# On the asymptotic geometry of finite-type $k$-surfaces in three-dimensional hyperbolic space.

Graham Smith

*Dedicated to the memory of Mark Sinclair-McGarvie*

**Abstract:** For $0 < k < 1$, a *finite-type $k$-surface* in 3-dimensional hyperbolic space is a complete, immersed surface of finite area and of constant extrinsic curvature equal to $k$. In [32], we showed that such surfaces have finite genus and finitely many cusp-like ends. Each of these cusps is asymptotic to an immersed cylinder of exponentially decaying radius about a complete geodesic and terminates at an ideal point which we call the *extremity* of the cusp. We show that every cusp of any finite-type $k$-surface has a well-defined axis, which we will call the *Steiner geodesic* of the cusp. One of its end-points is the extremity, and we will call the other, which constitutes new geometric data, the *Steiner point* of the cusp. We prove a new identity involving extremities and Steiner points in terms of Möbius invariant vector fields over the Riemann sphere.

We define two new functionals over the space of finite-type $k$-surfaces. The first, which will be called the *generalized volume*, is defined by the integral of a certain well-chosen form, and extends to the non-embedded case the concept of volume of the set bounded by the surface. The second, which will be called the *renormalized energy*, is related to the integral of the mean curvature of the surface, and is well-defined up to a choice of Busemann function. Upon describing natural parametrisations of the strata of the space of finite-type $k$-surfaces by open complex manifolds, we prove a new Schläfli-type formula relating the extremities and Steiner points to the first order variations of the generalized volume and the renormalized energy. In particular, Möbius invariance of this formula yields the aforementioned identity. We conclude by studying some applications of this identity and Schläfli-type formula.

## 1 - Introduction.

**1.1 - Overview.** Surfaces of constant extrinsic curvature in space forms have been natural objects of study since the publication of Gauss' famous Teorema Egregium in 1827, and the intriguing applications that they have found in recent years in such diverse fields as soliton theory, general relativity, hyperbolic geometry and Teichmüller theory[†] prove that, even almost two centuries after their introduction, they are still able to surprise and delight. In this paper, we study the asymptotic properties of a certain subclass of such surfaces in 3-dimensional hyperbolic space $\mathbb{H}^3$.

For $0 < k < 1$, a *finite-type $k$-surface* in $\mathbb{H}^3$ is a pair $(S, e)$, where $S$ is a smooth surface, and $e : S \to \mathbb{H}^3$ is a complete, smooth immersion of finite area and of constant extrinsic curvature equal to $k$. The space of reparametrisation equivalence classes of such surfaces identifies with the space of biholomorphism equivalence classes of marked ramified covers of the Riemann sphere $\hat{\mathbb{C}}$ (see Section 1.3). At this stage, it suffices to note that finite-type $k$-surfaces in $\mathbb{H}^3$ are topologically finite with cusp-like ends, as illustrated in Figure 1.1.1. Each cusp terminates in a well-defined point of $\partial_\infty\mathbb{H}^3 = \hat{\mathbb{C}}$, which we call its *extremity* and, since cusps are not necessarily embedded, each also has a well-defined positive integer-valued *winding number*.

In Section 1.2, we canonically associate to every cusp a complete geodesic in $\mathbb{H}^3$, which we view as the axis of the cusp, and which we call its *Steiner geodesic*. One of the end-points of this geodesic will be the extremity of the cusp, and we call the other the *Steiner point* of the cusp. A simple example of Steiner geodesics and Steiner points is illustrated in Figure 1.1.2.

---

[†] For recent applications of $k$-surfaces, the reader may consult [3], [6], [7], [19], [20], [25] and [33] as well as our review [9], written in collaboration with F. Fillastre. In the case where $k<0$, an attractive and modern discussion of the relationship between $k$-surfaces (there referred to as pseudospheres) and soliton theory is presented in [23].

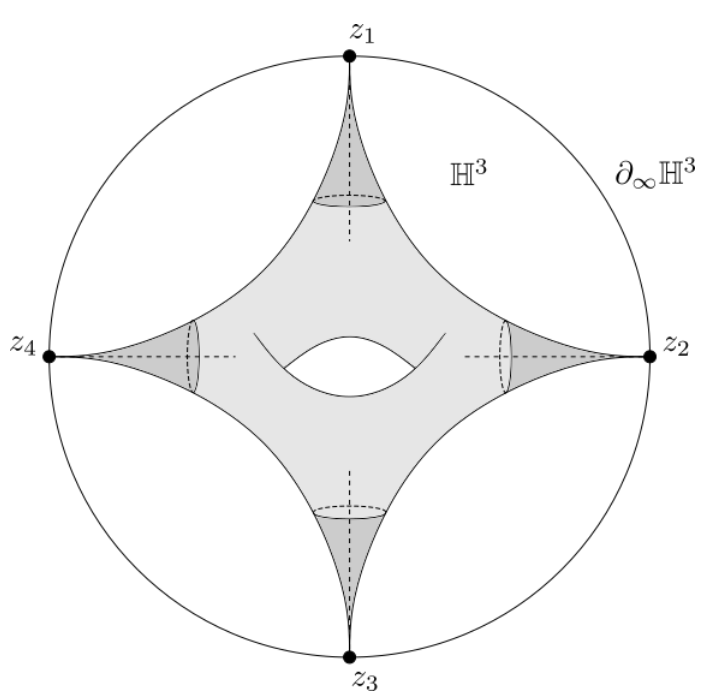

**Figure 1.1.1 - A typicial finite-type $k$-surface in $\mathbb{H}^3$ -** Such surfaces are topologically finite with cusp-like ends. Each cusp has a well-defined extremity in $\hat{\mathbb{C}}$ and a well-defined positive integer-valued winding number.

Our aim is to prove the following identity relating the extremities, winding numbers and Steiner points of any finite-type $k$-surface. For each $x \neq y \in \hat{\mathbb{C}}$, we define the vector field $\mathrm{K}[x,y]$ over $\hat{\mathbb{C}}$ by

$$\mathrm{K}[x,y](z) := \frac{(z-x)(z-y)}{(x-y)}\partial_z. \tag{1.1}$$

This family is Möbius-invariant, as is seen upon noting that $\mathrm{K}[x,y]$ is the Killing vector field of the stabiliser of $\{x,y\}$ in the Möbius group $\mathrm{SO}(3,1)$.

**Theorem 1.1.1**

*For every finite-type $k$-surface $(S,e)$ in $\mathbb{H}^3$ with extremities $z_1,\cdots,z_n$, winding numbers $m_1,\cdots,m_n$, and Steiner points $\zeta_1,\cdots,\zeta_n$,*

$$\sum_{i=1}^{n} m_i K[z_i,\zeta_i] = 0. \tag{1.2}$$

**Remark 1.1.1.** Theorem 1.1.1 is proven in Section 5.5.

Theorem 1.1.1 follows immediately from a Möbius invariant Schläfli-type formula that we now describe. We will see that, even in the non-embedded case, every finite-type $k$-surface $(S,e)$ bounds a well-defined, real-valued *generalized volume* $\mathrm{Vol}[e]$. In addition, given any Busemann function $h:\mathbb{H}^3\to\mathbb{R}$ with centre $h_\infty \in \hat{\mathbb{C}}$, we will say that $(S,e)$ is *$h$-admissable* whenever its extremities lie in $\hat{\mathbb{C}}\setminus\{h_\infty\}$ and, in this case, we will construct a well-defined renormalized energy $\hat{\mathrm{E}}[e;h]$ of the surface. We then define, for every $h$-admissable finite-type $k$-surface,

$$\mathcal{E}_k[e;h] := \hat{\mathrm{E}}[e;h] - 2(1+k)\mathrm{Vol}[e]. \tag{1.3}$$

For all $x\neq y\in\hat{\mathbb{C}}$, we define the 1-form $\lambda[x,y]$ over $\hat{\mathbb{C}}\setminus\{x,y\}$ by

$$\lambda[x,y](z) := \frac{(y-x)dz}{(y-z)(x-z)}. \tag{1.4}$$

This family is also Möbius-invariant, as is seen upon noting that, for all suitable $x$, $y$, $z$ and $\xi$

$$\lambda[x,y](z) = \partial_w[z,x,y,w]|_{w=z}\cdot\xi, \tag{1.5}$$

where $[\cdot,\cdot,\cdot,\cdot]$ here denotes the cross ratio. Finally, we note that the space of finite-type $k$-surfaces is stratified by finite-dimensional complex manifolds, with each stratum locally parametrised by the extremities of the cusps. In particular, given a finite-type $k$-surface $(S,e)$ with extremities $z_1,\cdots,z_n$, the tangent space of $\hat{\mathbb{C}}^n$ at $(z_1,\cdots,z_n)$ naturally identifies with the tangent space of the stratum at $(S,e)$. We obtain the following Schläfli-type formula.

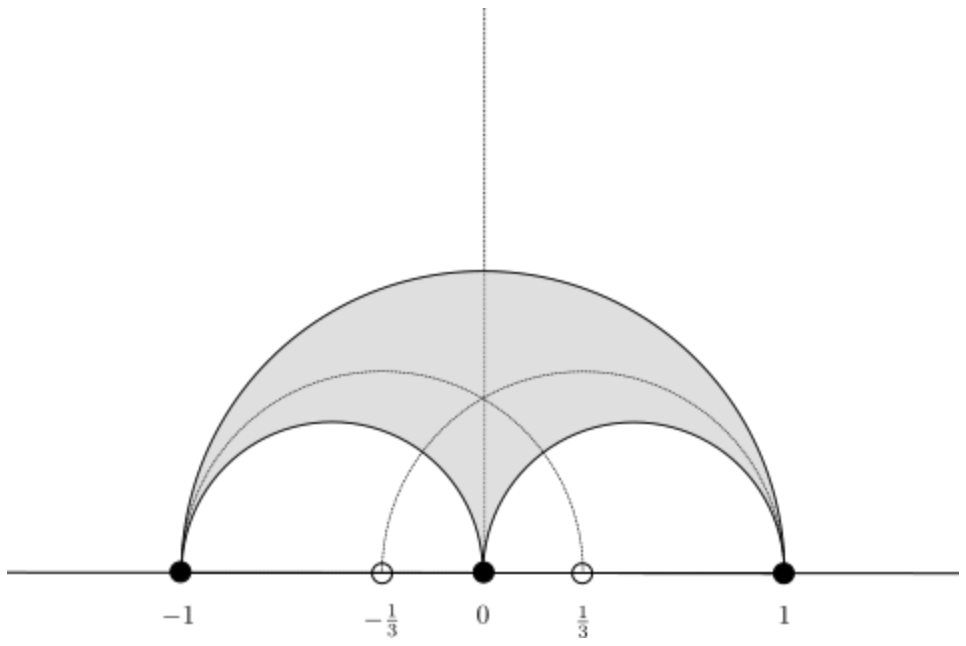

**Figure 1.1.2 - Steiner geodesics and Steiner points -** The embedded $k$-surface with three extremities at $-1$, 0 and 1 is a smooth fattening of the convex hull of these three points in $\mathbb{H}^3$. In this case, by symmetry, the Steiner geodesics are independent of $k$ and have Steiner points at $1/3$, $\infty$ and $-1/3$ respectively.

**Theorem 1.1.2**

*Let $h : \mathbb{H}^3 \to \mathbb{R}$ be a Busemann function with centre $h_\infty \in \hat{\mathbb{C}}$. For every $h$-admissable finite-type $k$-surface $(S, e)$ in $\mathbb{H}^3$ with extremities $z_1, \cdots, z_n$, winding numbers $m_1, \cdots, m_n$, and Steiner points $\zeta_1, \cdots, \zeta_n$, and for every tangent vector $\underline{\xi} = (\xi_1, \cdots, \xi_n)$ to $\hat{\mathbb{C}}^n$ at $(z_1, \cdots, z_n)$,*

$$D\mathcal{E}_k[e; h] \cdot \underline{\xi} = 4\pi \sum_{i=1}^{n} m_i Re(\lambda[h_\infty, \zeta_i](z_i) \cdot \xi_i). \tag{1.6}$$

**Remark 1.1.2.** Theorem 1.1.2 follows immediately from Theorem 5.5.5 and the subsequent remark. In particular, Theorem 1.1.1 follows immediately upon applying (1.6) to the stabiliser of $h$ in $\mathrm{SO}(3, 1)$.

**Remark 1.1.3.** The classical Schläfli formula is a staple of convex geometry which describes the first-order variation of the volume of a polyhedron in terms of the lengths of its edges and the first-order variations of its angles (see [1]). This formula has been extended in recent decades in various ways. For example, in [5], Bonahon adapted it to the case of equivariant pleated surfaces in hyperbolic space. More recently, in [15], Krasnov & Schlenker developed a version for first-order variations of smooth equivariant surfaces (c.f. also [24] and [26]). It is this smooth version that is closest in spirit to Theorem 1.1.2.

**1.2 - Finite-type $k$-surfaces, Steiner geodesics and Steiner points.** The remainder of this introduction will be devoted to detailing the results outlined in Section 1.1. In this section, we address the geometry of finite-type $k$-surfaces in $\mathbb{H}^3$. We refer the reader to [30] for proofs of the assertions that follow.

We first describe explicitely the ends of finite-type $k$-surfaces. Let $z \in \partial_\infty \mathbb{H}^3$ be an ideal point and let $h$ be a Busemann function of $\mathbb{H}^3$ centred on this point (c.f. [2]). For all $t \in \mathbb{R}$, let

$$H_t := h^{-1}(\{t\}) \tag{1.7}$$

denote the horosphere at height $t$ centred on $z$, and note that this surface is intrinsically euclidean. For a positive integer $m$, a *$k$-end* of *winding order* $m$ with *extremity* $z$ is defined to be a smooth immersion $e : S^1 \times [0, \infty[ \to \mathbb{H}^3$ of finite area and of constant extrinsic curvature equal to $k$ such that, for all $y$, $e(\cdot, y)$ is an immersed curve in $H_{-y}$ of total curvature equal to $2\pi m$. Up to rescaling, the metric that $e$ induces over $S^1 \times [0, \infty[$ is that of a hyperbolic cusp so that, in particular, the length of $e(\cdot, y)$ tends exponentially to zero as $y$ tends to infinity. Furthermore, for every unit-speed geodesic $\gamma : \mathbb{R} \to \mathbb{H}^3$ terminating at $z$ and parametrised by the value of $(-h)$,

$$d(\gamma(y), e(x, y)) = \mathrm{O}(e^{-y\sqrt{1-k}}) \tag{1.8}$$

as $y$ tends to infinity. Heuristically, $e$ wraps, ever more tightly, $m$-times around a complete geodesic. We will provide a complete description of the asymptotic geometry of $k$-ends in Chapter 4.

Now let $(S, e)$ be a finite-type $k$-surface. By Huber's Theorem, $S$ has finite genus and finitely many ends. We may therefore suppose that $S = \overline{S} \setminus P$, where $\overline{S}$ is a compact surface and $P := \{p_1, \cdots, p_n\}$ is a finite subset of $\overline{S}$. The immersion $e$ extends uniquely to a continuous function $\overline{e}:\overline{S}\to\mathbb{H}^3\cup\partial_\infty\mathbb{H}^3$ which maps every point of $P$ to an ideal point in $\partial_\infty\mathbb{H}^3$. For all $1\leq i\leq n$, the $i$'th *extremity* $\mathrm{z}_i[e]$ of $(S,e)$ is defined by

$$\mathrm{z}_i[e] := \overline{e}(p_i). \tag{1.9}$$

For each $i$, let $h_i:\mathbb{H}^3\to\mathbb{R}$ be a Busemann function centred on $\mathrm{z}_i[e]$, and denote

$$S_0 := \left\{p\in S\ |\ (h_i\circ e)(p)\geq 0\ \forall i\right\}. \tag{1.10}$$

Upon modifying $h_1,\cdots,h_n$ if necessary, $S_0$ may be taken to be compact with smooth boundary and the complement in $S$ of its interior may be taken to consist of $n$ connected components $S_1,\cdots,S_n$ such that, for each $i$, $(S_i,e)$ is a reparametrisation of some $k$-end with extremity $\mathrm{z}_i[e]$. In particular, for each $i$, the $i$'th end of $(S,e)$ has a well-defined winding number, which we denote by $m_i$. This completes our description of the geometry of finite-type $k$-surfaces in $\mathbb{H}^3$.

We now describe the construction of Steiner geodesics and Steiner points. We first recall some results from the theory of planar curves. Given a locally strictly convex, immersed, closed curve $\phi:S^1\to\mathbb{R}^2$ with geodesic curvature $\kappa:S^1\to\mathbb{R}$, unit normal vector field $\nu:S^1\to S^1$, and winding number $m$, its *Steiner curvature centroid* is defined by

$$s(\phi) := \frac{1}{m\pi}\int_{S^1}\langle\phi,\nu\rangle\nu\kappa\mathrm{dl}, \tag{1.11}$$

where dl here denotes the length element of $\phi$. The Steiner curvature centroid is also expressed in terms of the support function of $\phi$ as follows. Suppose first that $\phi$ has unit winding number, so that, by strict convexity, $\nu$ is a diffeomorphism. Recall (see [29]) that the *support function* of $\phi$ is defined by

$$\psi := \langle\phi,\nu\rangle\circ\nu^{-1}, \tag{1.12}$$

so that, upon parametrising $S^1$ by the interval $]-\pi,\pi[$ in the natural manner, we obtain

$$s(\phi) = \frac{1}{\pi}\int_{-\pi}^{\pi}\psi(\theta)(\cos(\theta),\sin(\theta))d\theta. \tag{1.13}$$

That is, $s(\phi)$ is the first Fourier mode of the support function. We leave the reader to determine the straightforward extension of this formula to the case of general winding numbers.

Let $z\in\partial_\infty\mathbb{H}^3$ be an ideal point, let $e:S^1\times[0,\infty[\to\mathbb{R}$ be a $k$-end with extremity $z$, and let $h:\mathbb{H}^3\to\mathbb{R}$ be a Busemann function centred on $z$. Since horospheres are intrinsically euclidean, locally strictly convex closed curves in horospheres have well-defined Steiner curvature centroids. We thus define $s:[0,\infty[\to\mathbb{H}^3$ such that, for all $y$, $s(y)$ is the Steiner curvature centroid in the horosphere $H_{-y}$ of the immersed curve $e(\cdot,y)$.

**Theorem & Definition 1.2.1**

*There exists a unique unit-speed geodesic $\gamma:\mathbb{R}\to\mathbb{H}^3$ such that*

$$d(\gamma(y),s(y)) = O(e^{-y\sqrt{4-3k}}) \tag{1.14}$$

*as $y$ tends to infinity. We call $\gamma$ the Steiner geodesic of $e$, and we call its end-point at $-\infty$ the Steiner point of $e$.*

**Remark 1.2.4.** Theorem & Definition 1.2.1 is proven in Section 4.2, below.

**Remark 1.2.5.** Recall that any two distinct unit-speed geodesics terminating at $z$ and parametrised by the value of $(-h)$ are separated by a distance proportional to $e^{-y}$ so that, since $0<k<1$, the geodesic constructed in Theorem 1.2.1 is indeed unique.

**Remark 1.2.6.** Although the main interest for us of Theorem & Definition 1.2.1 is the construction of a canonical axis with noteworthy properties, the Steiner curvature centroid is itself a remarkable object which has attracted the attention of convex geometers since the middle of the nineteenth century. The reader may consult Section 14.4 of [11], Section 10.5 of [13], and the introduction of [27] for overviews of its history from different perspectives. Interestingly, an entertaining exercise (see [16]) shows that, in the case of a convex polygon, the Steiner curvature centroid coincides with the centre of mass of the system obtained by placing a mass equal to the magnitude of the exterior angle at each vertex so that, in the case of triangles, it coincides with the Kimberling centre $X_{1115}$.*

---

* It was another pleasant surprise in writing this paper to learn of the existence of an entire theory devoted to the study of triangle

**1.3 - Finite-type $k$-surfaces and marked ramified covers.** We now describe in detail the topology and geometry of the space of finite-type $k$-surfaces in $\mathbb{H}^3$. We show, in particular, how this space identifies with the space of marked ramified covers of the Riemann sphere $\hat{\mathbb{C}}$. We refer the reader to [30] and [32] for proofs of the assertions that follow.

For $0 < k < 1$, we denote the space of finite-type $k$-surfaces in $\mathbb{H}^3$ by $\hat{\mathcal{S}}_k$. Two finite-type $k$-surfaces $(S, e)$ and $(S', e')$ will be considered to be equivalent whenever there exists a diffeomorphism $\alpha : S \to S'$ such that $e = e' \circ \alpha$. The quotient space of $\hat{\mathcal{S}}_k$ by this equivalence relation will be denoted by $\mathcal{S}_k$. We often identify a finite-type $k$-surface $(S, e)$ with its equivalence class in $\mathcal{S}_k$.

We define a *marked ramified cover* of the Riemann sphere $\hat{\mathbb{C}}$ to be a triple $(\overline{S}, P, \phi)$, where $\overline{S}$ is a compact Riemann surface, $P$ is a finite subset of $\overline{S}$, and $\phi : \overline{S} \to \hat{\mathbb{C}}$ is a non-constant holomorphic map with ramification points contained in $P$. Elements of $P$ will be called *generalized ramification points*. The space of marked ramified covers of $\hat{\mathbb{C}}$ will be denoted by $\hat{\mathcal{R}}$. Two marked ramified covers $(\overline{S}, P, \phi)$ and $(\overline{S}', P', \phi')$ will be considered to be equivalent whenever there exists a conformal diffeomorphism $\alpha : \overline{S} \to \overline{S}'$ such that $P = \alpha^{-1}(P')$ and $\phi = \phi' \circ \alpha$. The quotient space of $\hat{\mathcal{R}}$ by this equivalence relation will be denoted by $\mathcal{R}$. We often identify a marked ramified cover $(\overline{S}, P, \phi)$ with its equivalence class in $\mathcal{R}$.

We construct a canonical bijection from $\mathcal{S}_k$ to $\mathcal{R}$ as follows. Let $(S, e)$ be a finite-type $k$-surface. Let $\hat{e} : S \to \mathrm{U}\mathbb{H}^3$ denote the unit normal vector field over $e$. Let $n : \mathrm{U}\mathbb{H}^3 \to \partial_\infty \mathbb{H}^3$ denote the *horizon map* defined such that, for every unit vector $\xi_x \in \mathrm{U}\mathbb{H}^3$,

$$n(\xi_x) := \gamma(+\infty), \tag{1.15}$$

where $\gamma : \mathbb{R} \to \mathbb{H}^3$ denotes the unique geodesic whose derivative at zero is $\xi_x$. The function

$$\phi_e := n \circ \hat{e} \tag{1.16}$$

extends to a ramified cover of $\hat{\mathbb{C}} = \partial_\infty \mathbb{H}^3$ by $\overline{S}$ whose ramification points are elements of $P$, and the map

$$\Phi_k : \mathcal{S}_k \to \mathcal{R}; (S, e) \mapsto (\overline{S}, P, \phi_e) \tag{1.17}$$

is the desired bijection (see [32]).

The spaces $\mathcal{S}_k$ and $\mathcal{R}$ also carry natural topologies with respect to which $\Phi_k$ is a homeomorphism. They are constructed as follows. Let $\overline{S}$ be a compact surface. Let $U$ be a set of finite subsets of $\overline{S}$ which is open in the Hausdorff topology. Let $V$ be an open subset of $C^0(\overline{S}, \mathbb{H}^3 \cup \partial_\infty \mathbb{H}^3)$. We define the subset $\Omega_s(\overline{S}, U, V)$ of $\hat{\mathcal{S}}_k$ by

$$\Omega_s(\overline{S}, U, V) := \left\{ (\overline{S} \setminus P, e) \mid \#P < \infty,\ P \in U,\ \overline{e} \in V \right\}. \tag{1.18}$$

We furnish $\hat{\mathcal{S}}_k$ with the topology generated by all sets of this form. The sequence $(\overline{S} \setminus P_m, e_m)_{m \in \mathbb{N}}$ converges to $(\overline{S} \setminus P_\infty, e_\infty)$ with respect to this topology if and only if $(P_m)_{m \in \mathbb{N}}$ converges to $P_\infty$ in the Hausdorff sense and $(\overline{e}_m)_{m \in \mathbb{N}}$ converges uniformly to $\overline{e}_\infty$. Now let $W$ be an open subset of $C^0(\overline{S}, \hat{\mathbb{C}})$. We define the subset $\Omega_r(\overline{S}, U, W)$ of $\hat{\mathcal{R}}$ by

$$\Omega_r(\overline{S}, U, W) := \left\{ (\overline{S}, P, \phi) \mid \#P < \infty,\ P \in U,\ \phi \in W \right\}. \tag{1.19}$$

We furnish $\hat{\mathcal{R}}$ with the topology generated by all sets of this form. As before, the sequence $(\overline{S}, P_m, \phi_m)_{m \in \mathbb{N}}$ converges to $(\overline{S}, P_\infty, \phi_\infty)$ with respect to this topology if and only if $(P_m)_{m \in \mathbb{N}}$ converges to $P_\infty$ in the Hausdorff sense and $(\phi_m)_{m \in \mathbb{N}}$ converges uniformly to $\phi_\infty$. Finally, $\mathcal{S}_k$ and $\mathcal{R}$ are furnished with the induced quotient topologies, and it is with respect to these topologies that the bijection $\Phi_k$ defined in (1.17) is homeomorphic.

We conclude by describing the stratified holomorphic structures of these spaces. Let $(S, e)$ be an element of $\mathcal{S}_k$ with $n$ ends. A nearby element $(S', e')$ of $\mathcal{S}_k$ lies on the same stratum whenever it has the same number of ends. Since the numbers and orders of ends can only be varied over a continuous family in $\mathcal{S}_k$ by splitting or coalescing existing ends, it follows that the unordered vector $(m_1, \cdots, m_n)$ of winding orders is constant over every stratum. In Section 5.2, we show that each stratum naturally has the structure of a smooth centres. The interested reader may consult [14] for further information on this hidden treasure of modern mathematics.

complex manifold, of dimension equal to the number of ends, which is locally conformally parametrised by the extremities of these ends in $\hat{\mathbb{C}}$.

The strata of $\mathcal{R}$ are defined similarly. Given an element $(\overline{S},P,\phi)$ of $\mathcal{R}$, a nearby element $(\overline{S}',P',\phi')$ lies on the same stratum whenever $P'$ has the same cardinality as $P$. This means that a continuous family in $\mathcal{R}$ lies on a given stratum whenever no generalized ramification points of order 1 are added or removed and no generalized ramification points split or coalesce. Every stratum of $\mathcal{R}$ likewise has the structure of a smooth complex manifold, of dimension equal to the cardinality of the generalized ramification set, which is locally conformally parametrised by the images of the generalized ramification points. It follows trivially from these definitions that $\Phi_k$ restricts to a conformal diffeomorphism from strata of $\mathcal{S}_k$ into strata of $\mathcal{R}$.

**1.4 - Area, generalized volume and renormalized energy.** Finally, we introduce three natural geometric functionals over $\mathcal{S}_k$ which are smooth over each stratum. These will be the functionals which contribute to the Schläfli-like formula of Theorem 1.1.2. The first functional, which we include for completeness, is the area

$$\mathrm{Area}[e]:=\int_S \mathrm{dArea}[e]. \tag{1.20}$$

It is of little geometric interest since, by elementary hyperbolic surface theory, it is constant over every stratum. Indeed, for all $(S,e)\in\mathcal{S}_k$,

$$\mathrm{Area}[e]=-\frac{2\pi\chi[S]}{(1-k)}, \tag{1.21}$$

where $\chi[S]$ denotes the Euler characteristic of $S$.

The second functional generalizes to the case of immersions the concept of volume bounded by an embedding. We first construct in Section 5.3 a natural family $(\alpha_z)_{z\in\partial_\infty\mathbb{H}^3}$ of primitives of the volume form of $\mathbb{H}^3$ parametrised by ideal points in $\hat{\mathbb{C}}$. The *generalized volume* of a finite-type $k$-surface $(S,e)$ is then defined by

$$\mathrm{Vol}[e]:=\int_S e^*\alpha_z, \tag{1.22}$$

for some ideal point $z$. We verify that this functional is finite and independent of the ideal point chosen, that it is smooth over every stratum of $\mathcal{S}_k$ and, whenever $e$ is embedded, that it coincides with the volume of the convex set that this embedding bounds.

The third functional is what we choose to call the renormalized energy, and is defined as follows. Let $h:\mathbb{H}^3\to\mathbb{R}$ be a Busemann function centred on $h_\infty\in\hat{\mathbb{C}}$. Let $(S,e)$ be a finite-type $k$-surface whose extremities $\mathrm{z}_1[e],\cdots,\mathrm{z}_n[e]$ all lie in $\hat{\mathbb{C}}\setminus\{h_\infty\}$. For each $i$, let $h_i:\mathbb{H}^3\to\mathbb{R}$ be a Busemann function centred on $\mathrm{z}_i[e]$ normalized such that the horosphere $h_i^{-1}(\{0\})$ is tangent to $h^{-1}(\{0\})$ at some point. For all $T\in\mathbb{R}$, define

$$\hat{\mathrm{E}}_T[e;h]:=\int_{S_T}\mathrm{H}[e]\mathrm{dArea}[e], \tag{1.23}$$

where

$$S_T:=\left\{p\in S\ |\ (h_i\circ e)(p)\geq T\ \forall i\right\}, \tag{1.24}$$

and $\mathrm{H}[e]$ denotes the mean curvature of $e$. The *renormalized energy* of $(S,e)$ with respect to $h$ is defined by

$$\hat{\mathrm{E}}[e;h]:=\underset{T\rightarrow-\infty}{\mathrm{Lim}}\hat{\mathrm{E}}_T[e;h]+\sum_{i=1}^n 2\pi m_iT, \tag{1.25}$$

where $m_1,\cdots,m_n$ are the winding orders of the ends of $e$. In Section 5.3, we show that this limit exists and defines a function over an open, dense subset of $\mathcal{S}_k$ which is smooth over every stratum. A different Busemann function will yield another renormalized energy, defined over a different open, dense subset, which differs from the first over every stratum by a certain quadratic function of the extremities. It is the derivatives of the generalised volume and renormalized energy which are studied in Theorem & Definition 1.1.2.

Note that the renormalized energy arises in a natural manner from the geometry of the immersion $\hat{e}:S\to\mathrm{U}\mathbb{H}^3$ already studied in Section 1.3. Indeed, the area form of this immersion with respect to the Sasaki metric of $\mathrm{U}\mathbb{H}^3$ is

$$\mathrm{d}\hat{\mathrm{E}}[e]=\frac{1}{k}\mathrm{H}[e]\mathrm{dArea}[e]. \tag{1.26}$$

Since $\hat{e}$ is asymptotic over every end of $S$ to a finite cover of a cylinder in $\mathrm{U}\mathbb{H}^3$, the integral of this form grows linearly with the absolute value of $T$ as $T$ tends to minus infinity, from which convergence in (1.25) follows. In [18] (see also [17]), Labourie showed that $\hat{e}$ is pseudo-holomorphic with respect to a suitable almost complex structure over $\mathrm{U}\mathbb{H}^3$ and, using this property, derived many valuable results. Since the area of a pseudoholomorphic curve is usually interpreted as an energy (see, for example, [22]), this justifies our terminology.

**1.5 - Applications I - lagrangian embeddings.** We conclude this introduction with some applications of Theorems 1.1.1 and 1.1.2. We first show how the extremities and Steiner points yield canonical lagrangian immersions of the strata of $\mathcal{S}_k$ into certain open Kähler manifolds. Indeed, let $\Omega$ denote the complement of the diagonal in $\hat{\mathbb{C}}\times\hat{\mathbb{C}}$, that is

$$\Omega := \{(z,\zeta)\ |\ z\neq\zeta\}\,. \tag{1.27}$$

We define the complex-valued symplectic form $\omega$ over this set by

$$\omega := \frac{1}{(\zeta-z)^2}d\zeta\wedge dz. \tag{1.28}$$

Note that the family of 1-forms $\lambda$ introduced in (1.4) satisfies, for all $h_\infty$,

$$d\lambda[h_\infty,\cdot] = \omega. \tag{1.29}$$

The Schläfli formula (1.6) thus immediately yields the following result.

**Theorem 1.5.1**

*Let $X$ be an $n$-dimensional stratum of $\mathcal{S}_k$. Let $\underline{m}:=(m_1,\cdots,m_n)$ be such that, for all $i$, $m_i$ is the winding number of the $i$'th end of any element of $X$. The function $(\mathrm{z}_1,\zeta_1,\cdots,\mathrm{z}_n,\zeta_n)$ defines a smooth, lagrangian immersion from $X$ into $(\Omega^n,\omega_{\underline{m}})$, where the symplectic form $\omega_{\underline{m}}$ is given by*

$$\omega_{\underline{m}} = Re(m_1\omega_1\oplus\cdots\oplus m_n\omega_n). \tag{1.30}$$

From a physical perspective, the extremities and Steiner points may be considered as observable quantities over $X$. We may thus interpret Theorem 1.5.1 in at least two distinct ways. On the one hand, in analogy with classical thermodynamics, this result means that, for each $i$, the $i$'th extremity and the $i$'th Steiner point together constitute a pair of conjugate variables over the stratum $X$. On the other hand, since $X$ is locally parametrised by its extremities, the locally defined function $\sigma$ given by

$$\sigma(\mathrm{z}_1,\cdots,\mathrm{z}_n) := (\zeta_1,\cdots,\zeta_n) \tag{1.31}$$

may be considered as the map of scattering through finite-type $k$-surfaces. We recall from linear scattering theory that the symmetry of the scattering matrix corresponds to the physical reversibility of the process being studied. In the non-linear case, this corresponds to the lagrangian property of the scattering map, and this is precisely the property established in Theorem 1.5.1.

**1.6 - Applications II - the geometry of Steiner points.** As a second application, we determine the Steiner points of certain elementary finite-type $k$-surfaces. At this stage, it will be convenient to choose an explicit upper half-space parametrisation of $\mathbb{H}^3$, namely

$$\mathbb{H}^3 := \left\{(x,y,z)^t\ |\ z>0\right\}, \tag{1.32}$$

with metric given by

$$g_{ij} := \frac{1}{z^2}\delta_{ij}. \tag{1.33}$$

In this parametrisation, the ideal boundary $\partial_\infty\mathbb{H}^3$ of $\mathbb{H}^3$ naturally identifies with the extended complex plane $\mathbb{C}\cup\{\infty\}$. Let $(S,e)$ be a finite-type $k$-surface, none of whose extremities $\mathrm{z}_1[e],\cdots,\mathrm{z}_n[e]$ lie at infinity. For each $i$, let $m_i$ denote the winding number of the $i$'th end, let $\zeta_i[e]$ denote its Steiner point, and define its *Steiner vector* $\mathrm{c}_i[e]$ by

$$\mathrm{c}_i[e] := \frac{1}{\overline{\zeta}_i[e] - \overline{\mathrm{z}}_i[e]}. \tag{1.34}$$

Observe that this vector is always finite, since the Steiner point of any end is trivially distinct from its extremity. Upon expanding (1.2) as a quadratic function of $z$ and analysing its coefficients, we obtain

$$\sum_{i=1}^{n} m_i \mathrm{c}_i[e]\overline{\mathrm{z}}_i[e] = -\frac{1}{2}\sum_{i=1}^{n} m_i. \tag{1.35}$$

It is this form of (1.2) that will be of most use to us in this section.

Recall from Section 1.3 that finite-type $k$-surfaces are defined via their marked ramified covers. Consider now the case where $\overline{S} = \hat{\mathbb{C}}$ and $\phi$ is the identity map. For any finite subset $P$ of $\hat{\mathbb{C}}$, $(S, e) := \Phi_k^{-1}(\overline{S}, P, \phi)$ is an embedded surface bounding a convex set in $\mathbb{H}^3$ whose intersection with $\partial_\infty \mathbb{H}^3$ is $P$. In particular, the set $P$ coincides with the set $\{z_1, \cdots, z_n\}$ of extremities of $(S, e)$. Suppose now that $P := \{z_1, \cdots, z_n\}$ where, for each $m$,

$$z_m := e^{\frac{2\pi i m}{n}}.$$

$P$ is then symmetric under reflection about the unit circle in $\mathbb{C}$ as well as under reflection about the real line generated by $z_i$ for all $i$. It follows that, for each $i$, the $i$'th Steiner point of $(S, e)$ is

$$\zeta_i[e] = -z_i,$$

and the $i$'th Steiner vector is therefore

$$\mathrm{c}_i[e] = -\frac{1}{2} z_i.$$

The case where $n = 5$ is illustrated in Figure 1.6.3.

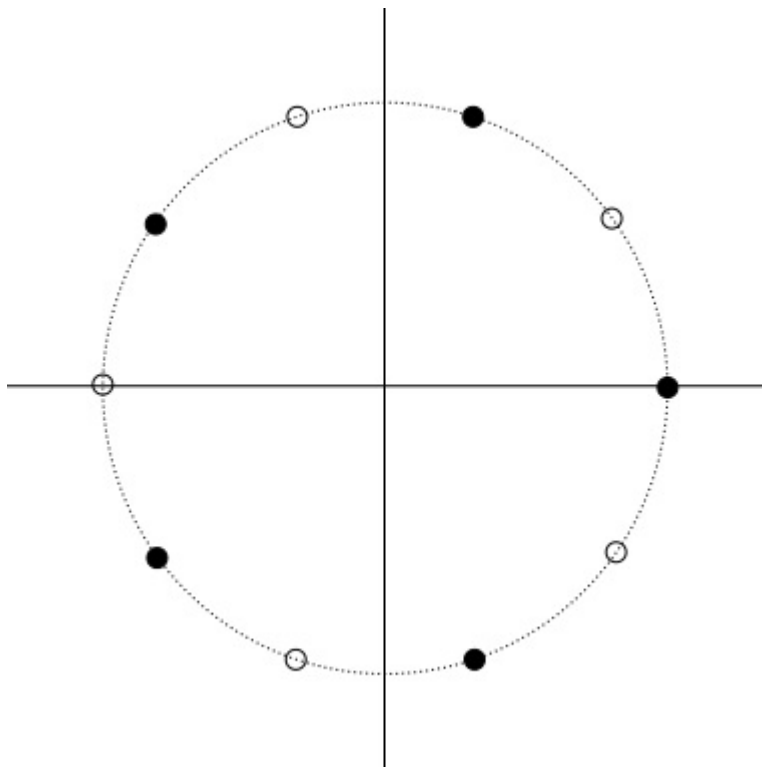

**Figure 1.6.3 - Steiner points I -** The extremities are shown in black and the Steiner points are shown in white. When the extremities of an embedded $k$-surface are evenly distributed along the unit circle in $\mathbb{C}$, the Steiner point of each end is the antipodal point of its extremity on the unit circle. Figure 1.1.2 is in fact obtained from the case of three points evenly distributed along a circle upon applying a suitable Möbius transformation.

Suppose now that $P := \{z_0, z_1, \cdots, z_n\}$, where $z_0 = 0$ and, for all $1 \leq i \leq n$, $z_i$ is as before. By symmetry again, the Steiner point of $(S, e)$ at 0 is

$$\zeta_0[e] = \infty,$$

and the corresponding Steiner vector is

$$\mathrm{c}_0[e] = 0.$$

Likewise, there exists a real number $a$ such that, for all $1 \leq i \leq n$, the Steiner point of $(S, e)$ at $z_i$ is

$$\zeta_i[e] = az_i.$$

However, symmetry alone is not sufficient to determine the value of $a$. Instead, using (1.35), we show that the Steiner vector of $(S, e)$ at $z_i$ is

$$\mathrm{c}_i[e] = -\frac{(n+1)}{2n} z_i,$$

so that the corresponding Steiner point is

$$\zeta_i[e] = \frac{(1-n)}{(1+n)} z_i.$$

The case where $n = 5$ is illustrated in Figure 1.6.4.

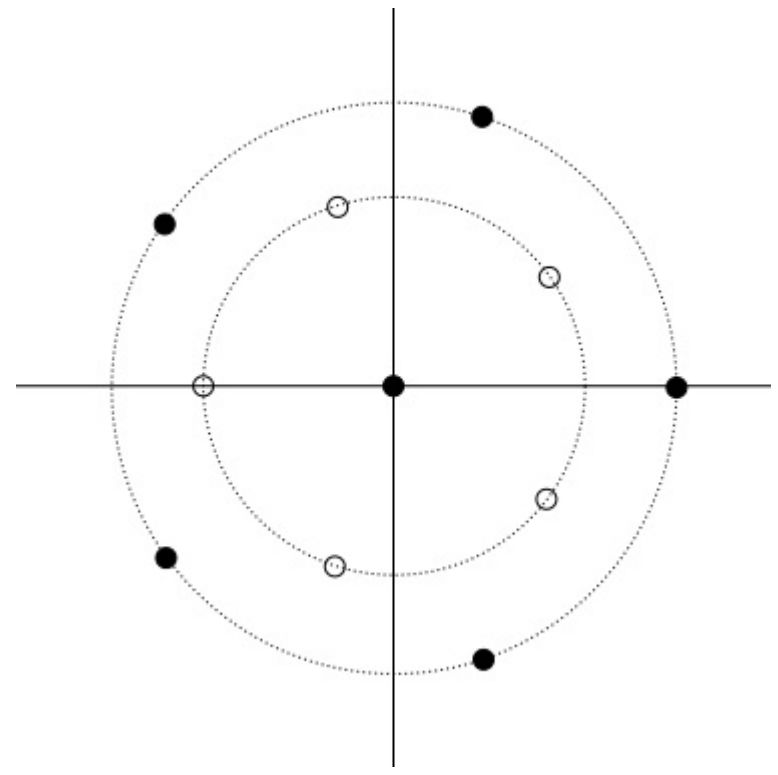

**Figure 1.6.4 - Steiner points II -** As before, the extremities are shown in black and the Steiner points are shown in white. The extra extremity at the origin displaces the other Stiener points towards the centre. In the case of 5 extremities evenly distributed along the unit circle, the Steiner points lie along the circle of radius 2/3 about the origin.

Finally, we construct a non-trivial covering of $\hat{\mathbb{C}}$ with a large number of symmetries as follows. Let $n$, $m_0$ and $m_1$ be positive integers such that

$$\frac{1}{m_0} + \frac{n}{m_1} \in \mathbb{Z}. \tag{1.36}$$

With $z_0, z_1, \cdots, z_n$ as before, let $\overline{S}$ be the Riemann surface of the function

$$f(z) = z^{\frac{1}{m_0}} \Pi_{i=1}^n (z - z_i)^{\frac{1}{m_1}},$$

and let $\phi : \overline{S} \to \hat{\mathbb{C}}$ denote the canonical projection. Condition (1.36) ensures that $\overline{S}$ is obtained from $\hat{\mathbb{C}}$ by branch cuts joining each $z_i$ by a radial line to the origin. In particular, the point at infinity is not the image of any ramification point of $\phi$. For each $i$, let $p_i \in \overline{S}$ be the unique preimage of $z_i$. By symmetry, the Steiner point of $(S, e)$ at $p_0$ is

$$\zeta_0[e] = \infty,$$

and the corresponding Steiner vector is

$$\mathrm{c}_0[e] = 0.$$

Likewise, by symmetry together with (1.35), for each $1 \leq i \leq n$, the Steiner vector of $(S, e)$ at $p_i$ is

$$\mathrm{c}_i[e] = -\frac{(m_0 + nm_1)}{2nm_1} z_i,$$

so that the corresponding Steiner point is

$$\zeta_i[e] = \frac{(m_0 - nm_1)}{(m_0 + nm_1)} z_i.$$

**1.7 - Notation.** Throughout this paper, we work with the upper half-space parametrisation of $\mathbb{H}^3$. The euclidean metric and norm will be denoted respectively by $\langle\cdot,\cdot\rangle_e$ and $\|\cdot\|_e$, and the hyperbolic metric and norm will be denoted respectively by $\langle\cdot,\cdot\rangle_g$ and $\|\cdot\|_g$. When describing asymptotic relations, we write $f(t)=\mathrm{O}(g(t))$ when $f(t)/g(t)$ remains bounded as $t$ tends to infinity, and we write $f(t)=\mathrm{o}(g(t))$ when $f(t)/g(t)$ tends to zero as $t$ tends to infinity.

**1.8 - Structure.** The paper is structured as follows.

**Section 2 :** We define Weinstein coordinates of the unitary bundle of hyperbolic space. These coordinates are well-adapted to the study of immersed surfaces that are asymptotic to cylinders over geodesic rays. We express the objects of classical differential geometry in terms of these coordinates. We conclude by determining a formula for the extrinsic curvature of an immersed surface.

**Section 3 :** We develop a theory of asymptotic series for solutions of a certain class of non-linear partial differential equations which includes $k$-ends. This is the most technical part of the paper and, since we are not aware of whether this problem has already been addressed in the literature, we study it in far greater detail than is necessary for our purposes.

**Section 4 :** We apply the analysis of Section 3 to determine asymptotic estimates of various geometric properties of $k$-ends, such as mean curvature, radius, area form, etc. In particular, we associate a well-defined axis to every $k$-end, thus proving Theorem 1.2.1.

**Section 5 :** We show how the immersions in each stratum vary smoothly with the extremities. We introduce the generalized volume and the renormalized energy. We study the derivatives of these functionals along the strata, and show that they satisfy the Schläfli formula of (1.6). We conclude by applying this Schläfli formula to the Killing fields of the isometry group of hyperbolic space, thereby proving the relations of Theorem 1.1.1.

**Appendix A :** For the reader's convenience, we review the theory of composition operators over Hölder spaces.

**1.9 - Acknowledgements.** The author is grateful to François Fillastre for inspiring conversations that led to the study of this problem; to the anonymous reviewer for very helpful comments made to earlier drafts of this paper; and to Giulio Belletti, Vladimir Fock and François Labourie for insightful comments leading to the present, much improved, draft of the paper.

## 2 - Geometry in Weinstein coordinates.

**2.1 - Overview.** Let $M:=M^{2m+1}$ be a contact manifold with contact form $\alpha$. An $m$-dimensional immersed submanifold $L:=L^m$ is said to be *legendrian* whenever the restrictions of $\alpha$ and $d\alpha$ to $L$ vanish. The total space of $T^*L\oplus\mathbb{R}$ also carries a natural contact structure given by

$$\tilde{\alpha}=dt-\lambda, \tag{2.1}$$

where $\lambda$ here denotes the canonical Liouville form of $T^*L$. Weinstein's theorem for contact manifolds then affirms the existence of a neighbourhood $\Omega$ of $L$ in $M$ contactomorphic to a neighbourhood of the zero section in $T^*L\oplus\mathbb{R}$ (see Theorem 3.4.13 of [21] and the preceding discussion). We call such parametrisations *Weinstein coordinates* of $M$ about $L$, and it follows from the proof of Weinstein's theorem that such coordinate systems are far from unique.

Consider now the case where $M:=\mathrm{U}\mathbb{H}^3$ is the total space of the unit sphere bundle over hyperbolic space, which is well-known to carry a natural contact structure. When working with hyperbolic space, it is useful to be aware of a number of coordinate systems, since different systems highlight different features of its geometry. In this spirit, Weinstein coordinates of $\mathrm{U}\mathbb{H}^3$ play a useful role in the study of those immersed surfaces in $\mathbb{H}^3$ which, like $k$-ends, are asymptotic to narrow cylinders about complete geodesics. Indeed, the unit normal bundle $\mathrm{N}\Gamma$ over a complete geodesic $\Gamma$ in $\mathbb{H}^3$ is an embedded lagrangian surface in $\mathrm{U}\mathbb{H}^3$, and there is a natural choice of Weinstein coordinates about this surface which respects the symmetries of $\mathbb{H}^3$ that preserve $\Gamma$. In this section, we will explicitly describe these coordinates, and we will study how the geometric properties of certain classes of surfaces are described in this framework. We underline that the explicit coordinates described here are not actually indispensable to what follows, but serve to greatly simplify our presentation. That said, it turns out that they do, in fact, possess a number of surprising properties which we believe warrant further study.

**2.2 - Weinstein coordinates of the unitary bundle.** We identify $\mathbb{H}^3$ with the upper half-space in $\mathbb{R}^3$, that is

$$\mathbb{H}^3 := \left\{(x,y,z)^t \mid z>0\right\}. \tag{2.2}$$

Recall that the hyperbolic metric over this space is given by

$$g_{ij} := \frac{1}{z^2}\delta_{ij}. \tag{2.3}$$

The total space of the tangent bundle of $\mathrm{T}\mathbb{H}^3$ identifies with an open subset of $\mathbb{R}^3\times\mathbb{R}^3$ in the natural manner. The total space of the unit tangent bundle then identifies with a codimension 1 submanifold of this product, namely

$$\mathrm{U}\mathbb{H}^3 := \left\{(x,y,z,u,v,w)^t \mid z>0,\ u^2+v^2+w^2=z^2\right\}.$$

Consider the Liouville form defined over $\mathrm{TU}\mathbb{H}^3$ by

$$\lambda := \frac{1}{z^2}(udx+vdy+wdz).$$

This form, which is invariant under the action of isometries of $\mathbb{H}^3$, defines a contact structure over $\mathrm{U}\mathbb{H}^3$.

Consider now the complete geodesic

$$\Gamma_{0,\infty} := \left\{(0,0,z)^t \mid z>0\right\},$$

and let $\mathrm{N}\Gamma_{0,\infty}$ denote the bundle of unit, normal vectors over this geodesic, that is

$$\mathrm{N}\Gamma_{0,\infty} := \left\{(0,0,z,u,v,0)^t \mid z>0,\ u^2+v^2=z^2\right\}.$$

Since $\mathrm{N}\Gamma_{0,\infty}$ is an embedded lagrangian submanifold of $\mathrm{U}\mathbb{H}^3$, Weinstein coordinates about this surface may be constructed. First, let $\mathrm{T}^*\mathrm{N}\Gamma_{0,\infty}$ denote the total space of its cotangent bundle. If $\tilde{\lambda}$ denotes the canonical Liouville form of $\mathrm{T}^*\mathrm{N}\Gamma_{0,\infty}$, then the form

$$dt-\tilde{\lambda}$$

defines a contact structure over the product $\mathrm{T}^*\mathrm{N}\Gamma_{0,\infty}\times\mathbb{R}$. Since the universal cover of $\mathrm{N}\Gamma_{0,\infty}$ is isometric to $\mathbb{R}^2$, $\mathrm{T}^*\mathrm{N}\Gamma_{0,\infty}\times\mathbb{R}$ naturally identifies with a quotient of $\mathbb{R}^5$. An explicit system $\Phi:\mathbb{R}^5\rightarrow\mathrm{U}\mathbb{H}^3$ of Weinstein coordinates about $\Gamma_{0,\infty}$ is then given by

$$\Phi(x,y,u,v,t) := \left(e^yt\mathrm{cos}(x)-e^yu\mathrm{sin}(x), e^yt\mathrm{sin}(x)+e^yu\mathrm{cos}(x), e^y, -\mathcal{C}e^y\mathrm{cos}(x), -\mathcal{C}e^y\mathrm{sin}(x), \mathcal{S}e^y\right)^t, \tag{2.4}$$

where

$$\begin{aligned}\mathcal{C} &:= \frac{1}{\sqrt{1+(t+v)^2}}\ \text{and}\\ \mathcal{S} &:= \frac{(t+v)}{\sqrt{1+(t+v)^2}}.\end{aligned} \tag{2.5}$$

Indeed, direct computation yields

$$\Phi^*\lambda = \frac{-1}{\sqrt{1+(t+v)^2}}(dt-udx-vdy). \tag{2.6}$$

This system of Weinstein coordinates is equivariant with respect to the group of isometries of $\mathbb{H}^3$ which preserve the point at infinity. Indeed, for all $\xi,\eta\in\mathbb{R}$ and for all $(a,b)^t\in\mathbb{R}^2$,

$$\begin{aligned}\mathrm{R}[\xi]_*\Phi(x,y,u,v,t) &= \Phi(x+\xi,y,u,v,t),\\ \mathrm{D}[\eta]_*\Phi(x,y,u,v,t) &= \Phi(x,y+\eta,u,v,t)\text{ and}\\ \mathrm{T}[a,b]_*\Phi(x,y,u,v,t) &= \Phi(x,y,u+\sigma_x(x,y),v+\sigma_y(x,y),t+\sigma(x,y)),\end{aligned}\tag{2.7}$$

where the hyperbolic isometries $\mathrm{R}[\xi]$, $\mathrm{D}[\eta]$ and $\mathrm{T}[a,b]$ are defined by

$$\begin{aligned}\mathrm{R}[\xi](x,y,z)^t &:= (\cos(\xi)x-\sin(\xi)y,\sin(\xi)y+\cos(\xi)x,z)^t,\\ \mathrm{D}[\eta](x,y,z)^t &:= (\eta x,\eta y,\eta z)^t\text{ and}\\ \mathrm{T}[a,b](x,y,z)^t &:= (x+a,y+b,z)^t,\end{aligned}\tag{2.8}$$

and the function $\sigma$ is defined by

$$\sigma(x,y) := \sigma[a,b](x,y) := ae^{-y}\cos(x)+be^{-y}\sin(x).\tag{2.9}$$

Finally, it will be convenient to introduce the variable

$$\theta := \arctan(t+v),\tag{2.10}$$

which is none other than the angle that the vector $\Phi(x,y,u,v,t)$ makes with the horizontal horosphere at height $e^y$. In particular

$$\mathcal{C}=\cos(\theta)\text{ and }\mathcal{S}=\sin(\theta),\tag{2.11}$$

justifying our notation. In addition, we denote

$$\mathcal{T} := \tan(\theta)=t+v.\tag{2.12}$$

These three abbreviations will be used frequently throughout the sequel.

**2.3 - The horizontal and vertical subbundles.** Let $\mathrm{WU}\mathbb{H}^3\subseteq\mathrm{TU}\mathbb{H}^3$ denote the contact distribution of $\mathrm{U}\mathbb{H}^3$. Let $\nabla$ denote the Levi-Civita covariant derivative of $\mathbb{H}^3$. We recall (see, for example, [31]) that $\mathrm{WU}\mathbb{H}^3$ decomposes as

$$\mathrm{WU}\mathbb{H}^3=\mathrm{HU}\mathbb{H}^3\oplus\mathrm{VU}\mathbb{H}^3,\tag{2.13}$$

where $\mathrm{HU}\mathbb{H}^3$ denotes the intersection of $\mathrm{WU}\mathbb{H}^3$ with the horizontal subbundle of $\nabla$ and $\mathrm{VU}\mathbb{H}^3$ denotes the vertical subbundle of $\mathrm{TU}\mathbb{H}^3$. Let

$$\Gamma := \nabla - D$$

denote the Christoffel symbol of $\nabla$, where $D$ here denotes the standard derivative of $\mathbb{R}^3$. We view $\Gamma$ as a symmetric bilinear form taking values in $\mathbb{R}^3$. By the Koszul formula, with respect to the standard basis $(\partial_x,\partial_y,\partial_z)$ of $\mathbb{R}^3$, it can be written as

$$\Gamma(x,y,z)=\frac{1}{z}\begin{pmatrix}\partial_z & 0 & -\partial_x\\ 0 & \partial_z & -\partial_y\\ -\partial_x & -\partial_y & -\partial_z\end{pmatrix}.\tag{2.14}$$

Consider now a point $(\underline{x},\underline{u})^t\in\mathrm{U}\mathbb{H}^3$. Given a tangent vector $\underline{\xi}$ of $\mathbb{H}^3$ at $\underline{x}$, its horizontal and vertical lifts to $\mathrm{T}_{(\underline{x},\underline{u})}\mathrm{U}\mathbb{H}^3$ are given by

$$\begin{aligned}[\underline{\xi},0]_{(\underline{x},\underline{u})} &:= (\underline{\xi},-\Gamma(\underline{x})(\underline{u},\underline{\xi}))^t\text{ and}\\ [0,\underline{\xi}]_{(\underline{x},\underline{u})} &:= (0,\underline{\xi})^t.\end{aligned}\tag{2.15}$$

The fibres over $(\underline{x},\underline{u})$ of the horizontal and vertical subspaces of $\mathrm{WU}\mathbb{H}^3$ are then given by

$$\begin{aligned}\mathrm{H}_{(\underline{x},\underline{u})}\mathrm{U}\mathbb{H}^3 &:= \left\{[\underline{\xi},0]_{(\underline{x},\underline{u})}\ \middle|\ \langle\underline{\xi},\underline{u}\rangle=0\right\}\text{ and}\\ \mathrm{V}_{(\underline{x},\underline{u})}\mathrm{U}\mathbb{H}^3 &:= \left\{[0,\underline{\xi}]_{(\underline{x},\underline{u})}\ \middle|\ \langle\underline{\xi},\underline{u}\rangle=0\right\}.\end{aligned}\tag{2.16}$$

In particular, there is a canonical bundle involution $\iota$ of $\mathrm{WU}\mathbb{H}^3$ defined such that, for all $(\underline{x},\underline{u})^t$ and for all $\underline{\xi}$,

$$
\begin{aligned}
\iota_{(\underline{x},\underline{u})}[\underline{\xi},0]_{(\underline{x},\underline{u})} &:= [0,\underline{\xi}]_{(\underline{x},\underline{u})}\ \text{and}\\
\iota_{(\underline{x},\underline{u})}[0,\underline{\xi}]_{(\underline{x},\underline{u})} &:= [\underline{\xi},0]_{(\underline{x},\underline{u})}.
\end{aligned}
\tag{2.17}
$$

We now determine $\mathrm{HU}\mathbb{H}^3$, $\mathrm{VU}\mathbb{H}^3$ and $\iota$ in the Weinstein coordinates defined in the previous section. By equivariance, we may suppose that $x=y=0$. Consider now the vector fields

$$
\begin{aligned}
\hat{\partial}_x &:= \partial_x + u\partial_t - u\partial_v - t\partial_u,\\
\hat{\partial}_y &:= \partial_y + v\partial_t - v\partial_v - u\partial_u,\\
\hat{\partial}_u &:= \partial_u\ \text{and}\\
\hat{\partial}_v &:= \mathcal{C}^2\partial_v.
\end{aligned}
\tag{2.18}
$$

We verify by inspection that these vector fields span $\Phi^*\mathrm{WU}\mathbb{H}^3$ and direct computation yields

$$
\begin{aligned}
\Phi_*(\mathcal{C}\hat{\partial}_u - \mathcal{S}\hat{\partial}_x) &= [(0,\mathcal{C},0),(0,0,0)],\\
\Phi_*(\mathcal{C}\hat{\partial}_y + \mathcal{S}\hat{\partial}_v) &= [(\mathcal{S},0,\mathcal{C}),(0,0,0)],\\
\Phi_*(-\hat{\partial}_x) &= [(0,0,0),(0,\mathcal{C},0)]\ \text{and}\\
\Phi_*(\hat{\partial}_v) &= [(0,0,0),(\mathcal{S},0,\mathcal{C})].
\end{aligned}
\tag{2.19}
$$

It follows that

$$
\begin{aligned}
\Phi^*\mathrm{HU}\mathbb{H}^3 &= \langle \mathcal{C}\hat{\partial}_u - \mathcal{S}\hat{\partial}_x, \mathcal{C}\hat{\partial}_y + \mathcal{S}\hat{\partial}_v\rangle,\\
\Phi^*\mathrm{VU}\mathbb{H}^3 &= \langle \hat{\partial}_v, \hat{\partial}_x\rangle,
\end{aligned}
\tag{2.20}
$$

and, with respect to the basis $(\hat{\partial}_u,\hat{\partial}_y,\hat{\partial}_x,\hat{\partial}_v)$,

$$
\Phi^*\iota = \begin{pmatrix} -\mathcal{S} & 0 & -\mathcal{C} & 0\\ 0 & -\mathcal{S} & 0 & \mathcal{C}\\ -\mathcal{C} & 0 & \mathcal{S} & 0\\ 0 & \mathcal{C} & 0 & \mathcal{S}\end{pmatrix}.
\tag{2.21}
$$

**2.4 - The geometry of legendrian immersions.** Let $\Omega$ be an open subset of $\mathbb{R}^2$. Let $u:\Omega\to\mathbb{R}$ be a smooth function. Define

$$
\hat{\Phi}[u] = \Phi\circ\hat{u},
\tag{2.22}
$$

where

$$
\hat{u}(x,y) := (x,y,u_x(x,y),u_y(x,y),u(x,t))^t.
\tag{2.23}
$$

Since $\hat{u}$ is a legendrian graph, $\hat{\Phi}[u]$ is a legendrian immersion. Furthermore, every immersed legendrian surface in $\mathrm{U}\mathbb{H}^3$ that is sufficiently close to $\mathrm{N}\Gamma_{0,\infty}$ is everywhere locally the image of such an immersion. Define also

$$
\Phi[u] := \pi\circ\hat{\Phi}[u],
\tag{2.24}
$$

where $\pi:\mathrm{U}\mathbb{H}^3\to\mathbb{H}^3$ is the canonical projection, so that

$$
\Phi[u] = (e^y u\cos(x) - e^y u_x\sin(x), e^y u\sin(x) + e^y u_x\cos(x), e^y).
\tag{2.25}
$$

We now review the elementary geometry of the map $\Phi[u]$. With respect to the bases $(\partial_x,\partial_y)$ of the domain and $(\hat{\partial}_u,\hat{\partial}_y,\hat{\partial}_x,\hat{\partial}_v)$ of the codomain, the derivative of $\hat{u}$ is given by

$$
D\hat{u} = \begin{pmatrix} M\\ N\end{pmatrix},
\tag{2.26}
$$

where

$$\begin{aligned} M &:= \begin{pmatrix} u_{xx}+u & u_{xy}+u_x \\ 0 & 1 \end{pmatrix} \text{ and} \\ N &:= \mathcal{C}^{-2}\begin{pmatrix} \mathcal{C}^2 & 0 \\ u_{yx}+u_x & u_{yy}+u_y \end{pmatrix}. \end{aligned} \tag{2.27}$$

It follows by (2.19) that

$$\begin{aligned} \Phi[u]_*\partial_x &= (0, u_{xx}+u, 0) \text{ and} \\ \Phi[u]_*\partial_y &= (\mathcal{T}, u_{xy}+u_x, 1). \end{aligned} \tag{2.28}$$

The first fundamental form of $\Phi[u]$ is thus given by

$$\mathrm{I}[u] := \begin{pmatrix} (u_{xx}+u)^2 & (u_{xx}+u)(u_{xy}+u_x) \\ (u_{xx}+u)(u_{xy}+u_x) & (u_{xy}+u_x)^2+\mathcal{C}^{-2} \end{pmatrix}, \tag{2.29}$$

its area form is

$$\mathrm{dArea}[u] = \mathcal{C}^{-1}(u_{xx}+u)dxdy, \tag{2.30}$$

and the length element that it induces over every horizontal curve is

$$\mathrm{dl}[u] = (u_{xx}+u)dx. \tag{2.31}$$

In particular, $\Phi[u]$ is an immersion if and only if

$$u_{xx}+u \neq 0, \tag{2.32}$$

The unit normal vector field $\mathrm{N}[u]$ over $\Phi[u]$ is, up to sign, simply the second component of $\hat{\Phi}[u]$, so that

$$\mathrm{N}[u] := e^y(\mathcal{C}\cos(x), \mathcal{C}\sin(x), -\mathcal{S})^t. \tag{2.33}$$

For all $y$, the restriction of $\Phi[u]$ to $\mathbb{R}\times\{y\}$ is the intersection of this immersion with the horizontal horosphere at height $e^y$. Let $\Phi^y[u]$ denote this restriction. We observe in passing that, when this immersion is locally strictly convex, $u(\cdot,y)$ is none other than the support function of this restriction (c.f. Section 1.2). The unit conormal vector field over this curve is

$$\nu[u] := e^y(\mathcal{S}\cos(x), \mathcal{S}\sin(x), \mathcal{C})^t. \tag{2.34}$$

Indeed, this vector field is orthogonal to the normal $\mathrm{N}[u]$ and, denoting by $\wedge_h$ the cross product of $\mathrm{T}\mathbb{H}^3$ compatible with its orientation, we see that

$$\mathrm{T}[u] := \mathrm{N}[u]\wedge_h \nu[u] = e^y(\sin(x), -\cos(x), 0)^t \tag{2.35}$$

is tangent to the curve. Using (2.19), we also obtain the useful formula

$$\nu[u] = \mathcal{C}\Phi_*\partial_y - \mathcal{C}\frac{(u_{xy}+u_x)}{(u_{xx}+u)}\Phi_*\partial_x. \tag{2.36}$$

The triplet $(\mathrm{N}[u], \nu[u], \mathrm{T}[u])$ defines an orthonormal frame over the immersion $\Phi[u]$. Furthermore, viewing $\Phi[u]$ itself as a vector field over this immersion, we obtain

$$\|\Phi[u]\|_g^2 = 1+u^2+u_x^2, \tag{2.37}$$

and

$$\begin{aligned} \langle\Phi[u], \mathrm{N}[u]\rangle_g &= -\mathcal{C}u_y, \\ \langle\Phi[u], \nu[u]\rangle_g &= \mathcal{C}+\mathcal{S}u \text{ and} \\ \langle\Phi[u], \mathrm{T}[u]\rangle_g &= -u_x. \end{aligned} \tag{2.38}$$

**2.5 - Curvatures of legendrian immersions.** Let $\pi_H$ denote the projection onto $\Phi^*\mathrm{HU}\mathbb{H}^3$ along $\Phi^*\mathrm{VU}\mathbb{H}^3$ and let $\pi_V$ denote the projection onto $\Phi^*\mathrm{VU}\mathbb{H}^3$ along $\Phi^*\mathrm{HU}\mathbb{H}^3$. Observe that

$$\begin{aligned}(\Phi^*\iota)\circ\pi_H &= \pi_V\circ(\Phi^*\iota)\ \text{and}\\ (\Phi^*\iota)\circ\pi_V &= \pi_H\circ(\Phi^*\iota).\end{aligned} \tag{2.39}$$

**Lemma 2.5.1**

*With respect to the basis $(\partial_x,\partial_y)$ of $\mathbb{R}^2$, the shape operator of $\Phi[u]$ is*

$$A[u] := \begin{pmatrix}-\mathcal{S} & 0\\ 0 & -\mathcal{S}\end{pmatrix} + M^{-1}\begin{pmatrix}-\mathcal{C} & 0\\ 0 & \mathcal{C}\end{pmatrix}N, \tag{2.40}$$

*where $M$ and $N$ are the matrices defined in* (2.27).

**Proof:** Indeed, with respect to the basis $(\partial_x,\partial_y)$ of $\mathbb{R}^2$ and the bases $(\hat{\partial}_u,\hat{\partial}_y)$ of $\Phi^*\mathrm{HU}\mathbb{H}^3$ and $(\hat{\partial}_x,\hat{\partial}_v)$ of $\Phi^*\mathrm{VU}\mathbb{H}^3$,

$$\begin{aligned}\pi_H\circ D\hat{u} &= M\ \text{and}\\ \pi_V\circ D\hat{u} &= N,\end{aligned}$$

so that,

$$\Phi^*\iota\circ\pi_V\circ D\hat{u} = \pi_H\circ\Phi^*\iota\circ D\hat{u} = \begin{pmatrix}-\mathcal{S} & 0\\ 0 & -\mathcal{S}\end{pmatrix}M + \begin{pmatrix}-\mathcal{C} & 0\\ 0 & \mathcal{C}\end{pmatrix}N.$$

Since $\mathrm{HU}\mathbb{H}^3$ is the horizontal bundle of the Levi-Civita covariant derivative of $\mathbb{H}^3$, the image of $D\hat{u}$ coincides with the graph of $A$, that is

$$(\pi_H\circ D\hat{u})\circ A = \Phi^*\iota\circ\pi_V\circ D\hat{u},$$

and the result follows. $\square$

Upon taking the trace and the determinant of (2.40), we obtain

**Lemma 2.5.2**

*The mean and extrinsic curvatures of $\Phi[u]$ are*

$$H[u] := \frac{\mathcal{C}}{(u_{xx}+u)} + 2\mathcal{S} - \frac{1}{\mathcal{C}(u_{xx}+u)}\big((u_{xx}+u)(u_{yy}+u_y)-(u_{xy}+u_x)^2\big),\ \textit{and} \tag{2.41}$$

$$K[u] := \mathcal{S}^2 + \frac{1}{(u_{xx}+u)}\big(\mathcal{S}\mathcal{C}+\mathcal{T}(u_{xy}+u_x)^2-\mathcal{T}(u_{xx}+u)(u_{yy}+u_y)-(u_{yy}+u_y)\big). \tag{2.42}$$

In particular, it follows from (2.42) that $\Phi[u]$ has constant extrinsic curvature equal to $k$ if and only if

$$ku_{xx}+u_{yy}-(1-k)u = F(u,Du,D^2u), \tag{2.43}$$

where $F$ is an analytic function of its arguments vanishing up to order 2 at $(0,0,0)$.

Finally, recall the restriction $\Phi^y[u]$ of $\Phi[u]$ to $\mathbb{R}\times\{y\}$ defined in the preceding section. Let $\kappa^y[u]$ denote its geodesic curvature with respect to the unit normal $\nu[u]$. Although it will only be of secondary importance to our work, we also show

**Lemma 2.5.3**

*The geodesic curvature $\kappa^y[u]$ of $\Phi^y[u]$ is*

$$\kappa^y[u]=\frac{\mathcal{C}(u_y-u_{xx})}{(u+u_{xx})}. \tag{2.44}$$

**Proof:** Indeed, the upward pointing unit normal over the horosphere at height $e^y$ is

$$\mathrm{N}^H=(0,0,e^y)^t.$$

Since every horosphere is totally umbilic with unit curvature,

$$\langle\nabla_{\Phi_*\partial_x}\mathrm{N}^H(\Phi[u]),\partial_x\rangle_g=-\|\Phi[u]_*\partial_x\|_g^2=-(u+u_{xx})^2.$$

On the other hand, bearing in mind (2.29) and (2.40),

$$\langle\nabla_{\Phi_*\partial_x}\mathrm{N}[u],\Phi_*\partial_x\rangle_g=\mathrm{I}[u](\mathrm{A}[u]\cdot\partial_x,\partial_x)=\mathcal{C}(u+u_{xx})+\mathcal{S}(u+u_{xx})^2.$$

Since

$$\nu[u]=\mathcal{C}^{-1}\mathrm{N}^H(\Phi[u])+\mathcal{T}\mathrm{N}[u],$$

the preceeding relations yield

$$\langle\nabla_{\Phi_*\partial_x}\nu[u],\Phi[u]_*\partial_x\rangle_g=-\mathcal{C}(u+u_{xx})^2+\mathcal{C}\mathcal{T}(u+u_{xx})=\mathcal{C}(u_y-u_{xx})(u+u_{xx}),$$

and (2.44) follows upon dividing both sides by $\|\Phi[u]_*\partial_x\|_g^2$. $\square$

## 3 - Asymptotic analysis.

**3.1 - Overview.** Theorem 1.2.1 will follow immediately from an asymptotic analysis of $k$-ends in Weinstein coordinates. Since this analysis applies to solutions of general totally non-linear partial differential equations defined over cylinders, and since we are not aware of whether it has already been carried out in the literature, we study it in far greater detail than is actually required for our work. Consider the non-linear second order partial differential operator

$$P[u]:=\partial_x^2u+\partial_y^2u-a^2u-F(u,Du,D^2u), \tag{3.1}$$

defined over the space of twice differentiable functions $u:S^1\times[0,\infty[\rightarrow\mathbb{R}$, where $a$ is a real constant and $F$ is a smooth function of its arguments. In what follows, we will suppose that

$$F(0,0,0),DF(0,0,0)=0, \tag{3.2}$$

and we will study solutions of the problem

$$\mathrm{P}[u]=0. \tag{3.3}$$

When $F$ vanishes, such solutions are completely described via the classical technique of separation of variables. Our analysis consists of a perturbation of this technique to the case of non-trivial $F$.

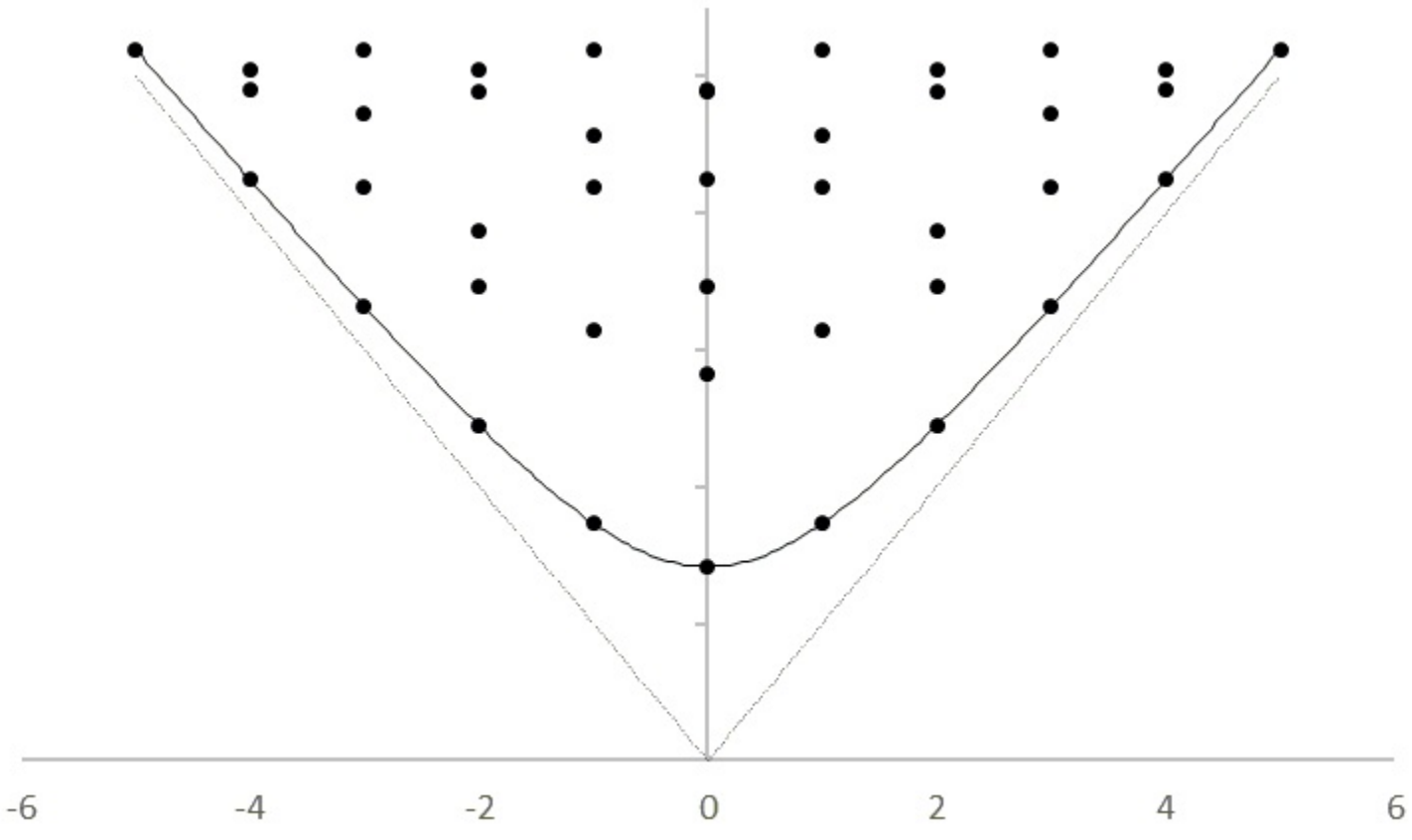

**Figure 3.1.5 - The index set -** The index set $\mathcal{M}$ is the subsemigroup of $\mathbb{R}^2$ generated by the set $\mathcal{M}_0$ consisting of those points of the hyperbola with integer $x$-coordinate.

Let $\mathcal{M}$ denote the subsemigroup of $\mathbb{R}\times\mathbb{R}$ generated by the set

$$\mathcal{M}_0 := \left\{(n,\sqrt{a^2+n^2})\ |\ n\in\mathbb{Z}\right\}, \tag{3.4}$$

as illustrated in Figure 3.1.5. Let $\mathcal{A}^0$ denote the vector space of all continuous functions $u: S^1\times[0,\infty[\rightarrow\mathbb{C}$ having the property that there exists a family $(a_{\lambda,\mu})_{(\lambda,\mu)\in\mathcal{M}}$ of complex constants indexed by $\mathcal{M}$ such that, for all $\omega>0$, there exists $C>0$ such that, for all $(x,y)$,

$$\left|u(x,y)-\sum_{(\lambda,\mu)\in\mathcal{M},\mu<\omega}a_{\lambda,\mu}e^{i\lambda x}e^{-\mu y}\right|\leq Ce^{-\omega y}. \tag{3.5}$$

Observe that, for all $(\lambda,\mu)$, the coefficient $a_{\lambda,\mu}$ is uniquely determined by $u$ and varies linearly with this function. When $u\in\mathcal{A}^0$, we write

$$u\sim\sum_{(\lambda,\mu)\in\mathcal{M}}a_{\lambda,\mu}e^{i\lambda x}e^{-\mu y}. \tag{3.6}$$

For all non-negative, integer $k$, let $\mathcal{A}^k$ denote the space of all $k$-times differentiable functions $u: S^1\times[0,\infty[\rightarrow\mathbb{R}$ all of whose derivatives up to and including order $k$ are elements of $\mathcal{A}^0$. Observe that when $u\in\mathcal{A}^k$ the asymptotic series of its derivatives are obtained by differentiating term by term the asymptotic series of $u$. We define

$$\mathcal{A} := \bigcap_{k\in\mathbb{N}}\mathcal{A}^k. \tag{3.7}$$

That is, $\mathcal{A}$ is the space of those functions $u: S^1\times[0,\infty[\rightarrow\mathbb{C}$ all of whose derivatives to all orders have asymptotic series of the form (3.6). We show in (3.39), below, that $\mathcal{A}$ is the intersection of a family of Banach spaces, from which it follows that it carries a natural Frechet structure given by the family of all norms of these Banach spaces. For all $(\lambda,\mu)\in\mathcal{M}$, we define the linear operator $\mathrm{a}_{\lambda,\mu}:\mathcal{A}\rightarrow\mathbb{C}$ such that, for all $u$,

$$\mathrm{a}_{\lambda,\mu}[u] := a_{\lambda,\mu}, \tag{3.8}$$

where $a_{\lambda,\mu}$ is the coefficient of $e^{i\lambda x}e^{-\mu y}$ in the asymptotic series (3.6) of $u$. It follows directly from the definition of the Frechet structure of $\mathcal{A}$ that, for all $(\lambda,\mu)$, this operator is continuous in the Frechet sense.

In order to study solutions of (3.3), we develop a straightforward calculus for determining the asymptotic series of sums, products and compositions. Indeed, the sum of two elements of $\mathcal{A}$ is trivially an element of $\mathcal{A}$ whose asymptotic series is given by the sums of their asymptotic series. In Lemma 3.4.4, below, and the subsequent remark, we show that the product operator defines a continuous bilinear map from $\mathcal{A}\oplus\mathcal{A}$ to $\mathcal{A}$ where the asymptotic series of the product of two elements of $\mathcal{A}$ is obtained by formal multiplication of the asymptotic series of each of these elements. In other words, $\mathcal{A}$ is a Frechet algebra. Likewise, in Lemma 3.5.4, below, and the subsequent remark, we show that, given any smooth function $\Phi$ defined in a neighhourhood of 0 such that $\Phi(0)=0$, the operator of composition by $\Phi$ defines a map from a neighbourhood of 0 in $\mathcal{A}$ to $\mathcal{A}$ which is smooth in the Frechet sense, and the asymptotic series of the image of any element of this neighbourhood is obtained by formally substituting the asymptotic series of this element into the MacLaurin series of $\Phi$.

Now let $\mathrm{R}:\mathcal{A}\rightarrow C^\infty(S^1)$ denote the operator of restriction onto $S^1\times\{0\}$. This linear operator is trivially continuous in the Frechet sense. The main result of this chapter is

**Theorem 3.1.1**

*If $F(0,0,0), DF(0,0,0)=0$, then the following assertions hold.*

*(A) If $u: S^1\times[0,\infty[\rightarrow\mathbb{R}$ is a smooth solution of* (3.3) *such that, for some $\lambda>0$,*

$$\begin{aligned}|u(x,t)| &= O(e^{-\lambda t}),\ \text{and}\\ \|D^k u(x,y)\| &= O(1)\ \forall k\geqslant 1,\end{aligned} \tag{3.9}$$

*then $u$ is an element of $\mathcal{A}$.*

*(B) There exists a neighbourhood $U$ of $0$ in $C^\infty(S^1)$ and an operator $S : U \to \mathcal{A}$ which is smooth in the Frechet sense such that, for all $v \in U$,*

*(1) $RS[v] = v$, and*

*(2) $PS[v] = 0$.*

*(C) Upon reducing $U$ if necessary, we may suppose that $S$ is unique.*

**Remark 3.1.7.** Part $(A)$ of Theorem 3.1.1 follows from Lemma 3.6.4, below. Parts $(B)$ and $(C)$ follow from Theorem 3.6.2, below, and the subsequent remark.

Theorem 3.1.1 applies to $k$-ends as follows. Let $e : S^1 \times [0,\infty[\to \mathbb{H}^3$ be a $k$-end which, for ease of presentation, we take to be of winding order 1. Up to reparametrisation, we may suppose that $e = \Phi[u]$ where $\Phi$ is as in (2.24) and $u : S^1 \times [0,\infty[\to \mathbb{R}$ solves (2.43) for some suitable analytic function $F$. In Lemma 3.6.5, below, we verify that $u$ satisfies the hypotheses of Part $(A)$ of Theorem 3.1.1 and therefore has a well-defined asymptotic series of the form (3.6). We will show in Chapter 4, below, that two well-chosen coefficients of this asymptotic series are eliminated upon applying a unique translation of the form (2.8). Expressed in geometric terms, this corresponds precisely to the existence of the Steiner geodesic, thus proving Theorem 1.2.1

However, Theorem 3.1.1 says a good deal more. Indeed, Part $(B)$ of this theorem tells us that every coefficient of the asymptotic series of $u$ locally depends only on the restriction of this function to the boundary curve $S^1 \times \{0\}$ and that, furthermore, this dependence is smooth. Consider now a finite type $k$-surface. Upon applying a suitable cut-off, as in (1.10), we see that each of its ends has a well-defined asymptotic series. Furthermore, the coefficients of these series, and thus, in particular, the Steiner geodesics and Steiner points, depend locally, and in a smooth manner, only on the compact part $S_0$ of this surface, which immediately yields the smoothness of the immersions studied in Theorem 1.5.1.

**3.2 - The one-dimensional linear problem.** We begin by studying the case where $u$ is constant in $x$. We first recall the formalism of weighted Hölder spaces. Let $E$ be a Banach space. For a weight $\omega \in \mathbb{R}$ and for all $(k,\alpha)$, define the *$\omega$-weighted $C^{k,\alpha}$-norm* for $k$-times differentiable functions $u : \mathbb{R} \to E$ by

$$\|u\|_{C^{k,\alpha}_\omega} := \|ue^{\omega\langle y\rangle}\|_{C^{k,\alpha}}, \tag{3.10}$$

where $\|\cdot\|_{C^{k,\alpha}}$ denotes the standard Hölder norm (see Appendix A) and

$$\langle y\rangle := \sqrt{1+y^2}. \tag{3.11}$$

For all $\omega$ and for all $(k,\alpha)$, the Banach space of $k$-times differentiable functions $u : \mathbb{R} \to E$ with finite $C^{k,\alpha}_\omega$-norm will be denoted by $C^{k,\alpha}_\omega(\mathbb{R}, E)$. For the sake of computations, we observe that the $C^{k,\alpha}_\omega$-norm is uniformly equivalent to

$$\|u\|'_{C^{k,\alpha}_\omega} := \sum_{i=0}^{k} \|e^{\omega|y|}D^i u\|_{C^0} + [e^{\omega|y|}D^l u]_\alpha. \tag{3.12}$$

In what follows, we will use without comment the more appropriate of these norms depending on the context in which we are working.

For all $\omega' \in \mathbb{R}$, let $\mu_{\omega'}$ denote the operator of multiplication by the function $e^{-\omega'\langle y\rangle}$. By definition, for all $\omega$, $\omega'$ and for all $(k,\alpha)$, $\mu_{\omega'}$ defines a linear isomorphism from $C^{k,\alpha}_\omega(\mathbb{R}, E)$ into $C^{k,\alpha}_{\omega+\omega'}(\mathbb{R}, E)$. In particular, $\mu_{-\omega}$ defines a linear isomorphism from $C^{k,\alpha}_\omega(\mathbb{R}, E)$ into $C^{k,\alpha}(\mathbb{R}, E)$. It follows that norm estimates for any given linear operator $L : C^{k+l,\alpha}_\omega(\mathbb{R}, E) \to C^{k,\alpha}_\omega(\mathbb{R}, E)$ are equivalent to norm estimates for the conjugate operator $\mu_{-\omega}L\mu_\omega$, viewed as a linear map from $C^{k+l,\alpha}(\mathbb{R}, E)$ into $C^{k,\alpha}(\mathbb{R}, E)$. In other words, the study of the analytic properties of families weighted functional norms is equivalent to the study of the analytic properties of corresponding families of conjugates of operators. Since it is often easier to study conjugated operators than it is to study functional norms, this perspective will often be used implicitly throughout the sequel.

For a real constant $a$, consider now the linear operator

$$\tilde{\mathrm{L}}_a u := \partial_y^2 u - a^2 u. \tag{3.13}$$

**Lemma 3.2.1**

*For all $a,\omega\in\mathbb{R}$, $\tilde{L}_a$ defines a bounded linear map from $C^{2,\alpha}_\omega(\mathbb{R},E)$ onto $C^{0,\alpha}_\omega(\mathbb{R},E)$.*

**Proof:** Indeed, we verify by inspection that $\mu_{-\omega}\tilde{\mathrm{L}}_a\mu_\omega$ is a second order, linear operator with coefficients bounded in $C^{0,\alpha}(\mathbb{R},E)$. It therefore defines a bounded linear map from $C^{0+2,\alpha}(\mathbb{R},E)$ into $C^{0,\alpha}(\mathbb{R},E)$, and the result follows. □

We use potential theory to study the invertibility properties of this operator over different function spaces. For all $a>0$, the Green's function of $\tilde{\mathrm{L}}_a$ is

$$\tilde{K}_a(y) := -\frac{1}{2a}e^{-a|y|}, \tag{3.14}$$

and its Green's operator is

$$\tilde{\mathrm{K}}_a[u](y) := \int_{-\infty}^{\infty}\tilde{K}_a(y-z)u(z)dz. \tag{3.15}$$

**Lemma 3.2.2**

*For all $|\omega|<a$, the operator $\tilde{K}_a$ defines a bounded linear map from $C^0_\omega(\mathbb{R},E)$ to itself.*

**Proof:** The exponential decay of $u$ ensures that the integral (3.15) exists and varies continuously with $y$. A straightforward calculation then yields, for all $u\in C^0_\omega(\mathbb{R},E)$ and for all $y\in\mathbb{R}$,

$$\|\tilde{\mathrm{K}}_a[u](y)e^{\omega|y|}\| \leq \frac{1}{(a^2-\omega^2)}\|u\|'_{C^0_\omega} + \frac{\omega}{a(a^2-\omega^2)}e^{-(a-\omega)|y|}\|u\|'_{C^0_\omega}.$$

Since $(a-\omega)>0$, the coefficient of the second term is bounded and the result follows. □

**Lemma 3.2.3**

*For all $|\omega|<a$, the operator $\tilde{K}_a$ defines a bounded linear map from $C^{0,\alpha}_\omega(\mathbb{R},E)$ into $C^{2,\alpha}_\omega(\mathbb{R},E)$ such that*

$$\tilde{\mathrm{L}}_a\tilde{K}_a = Id. \tag{3.16}$$

**Proof:** Indeed, for all $u$,

$$\tilde{\mathrm{K}}_a[u](y) = -\int_{-\infty}^{y}\frac{1}{2a}e^{a(z-y)}u(z)dz - \int_y^\infty\frac{1}{2a}e^{a(y-z)}u(z)dz.$$

Differentiating this equation under the integral yields, for all $u$,

$$\partial_y\tilde{\mathrm{K}}_a[u](y) = \int_{-\infty}^{y}\frac{1}{2}e^{a(z-y)}u(z)dz - \int_y^\infty\frac{1}{2}e^{a(y-z)}u(z)dz,$$

and differentiating a second time yields, for all $u$,

$$\partial_y^2\tilde{\mathrm{K}}_a[u](y) = u(y) + a^2\tilde{\mathrm{K}}_a[u](y).$$

In particular, bearing in mind Lemma 3.2.2, for all $u\in C^{0,\alpha}_\omega(\mathbb{R},E)$, $\tilde{\mathrm{K}}_a[u]\in C^{2,\alpha}_\omega(\mathbb{R},E)$, and

$$\tilde{\mathrm{L}}_a\tilde{\mathrm{K}}_a[u] = u.$$

However, it is straightforward to verify that there exists $A_1>0$ such that, for all $v\in C^{2,\alpha}_\omega(\mathbb{R},E)$,

$$\|v\|_{C^{2,\alpha}_\omega} \leq A_1\big(\|v\|_{C^0_\omega} + \|\tilde{\mathrm{L}}_a v\|_{C^{0,\alpha}_\omega}\big),$$

so that, by Lemma 3.2.2 again, there exists $A_2>0$ such that, for all $u\in C^{0,\alpha}_\omega(\mathbb{R},E)$,

$$\|\tilde{\mathrm{K}}_a u\|_{C^{2,\alpha}_\omega} \leq A_2\|u\|_{C^{0,\alpha}_\omega},$$

and the result follows. □

**Theorem 3.2.4**

*For all $0 \leq \omega < a$, $\tilde{L}_a$ defines a linear isomorphism from $C^{2,\alpha}_\omega(\mathbb{R},E)$ into $C^{0,\alpha}_\omega(\mathbb{R},E)$ with inverse $\tilde{K}_a$.*

**Proof:** Indeed, by Lemma 3.2.3, $\tilde{\mathrm{K}}_a$ defines a right inverse of $\tilde{\mathrm{L}}_a$. By the maximum principle, $\tilde{\mathrm{L}}_a$ is injective, and the result follows. $\square$

The preceding results adapt to the Dirichlet problem as follows. First, for all $\omega$ and for all $(k,\alpha)$, the Banach space $C^{k,\alpha}_\omega([0,\infty[,E)$ is defined in the natural manner and its closed subspace $C^{k,\alpha}_{\omega,0}([0,\infty[,E)$ is defined to consist of those functions which vanish at zero. We verify that the Green's operator of $\tilde{\mathrm{L}}_a$ for the Dirichlet problem is

$$\tilde{\mathrm{K}}_{a,0}[u](y) := \int_0^\infty \tilde{K}_a(y-z)u(z)dz - \int_{-\infty}^0 \tilde{K}_a(y-z)u(-z)dz. \tag{3.17}$$

Proceeding as before, we obtain

**Theorem 3.2.5**

*For all $0 \leq \omega < a$, $\tilde{L}_a$ defines a linear isomorphism from $C^{2,\alpha}_{\omega,0}([0,\infty[,E)$ into $C^{0,\alpha}_\omega([0,\infty[,E)$ with inverse $\tilde{K}_a$.*

We now consider the case where the weight $\omega$ is greater than $a$. This setting yields the richer structure underlying the asymptotic series described in Section 3.1. Consider a finite vector $\underline{\omega} := (\omega_0,\cdots,\omega_m)$ of real weights, where

$$a = \omega_0 < \cdots < \omega_m. \tag{3.18}$$

Define the spaces

$$\begin{aligned}
\tilde{\mathcal{A}}^{k,\alpha}_{\underline{\omega}} &:= \langle e^{-\omega_0 y},\cdots,e^{-\omega_{m-1}y}\rangle \oplus C^{k,\alpha}_{\omega_m}([0,\infty[),\\
\tilde{\mathcal{A}}^{k,\alpha}_{\underline{\omega},0} &:= \left\{ f \in \tilde{\mathcal{A}}^{k,\alpha}_{\underline{\omega}}\ \mid\ f(0)=0\right\}\ \text{and}\\
\tilde{\mathcal{A}}^{k,\alpha}_{\underline{\omega},*} &:= \langle e^{-\omega_1 y},\cdots,e^{-\omega_{m-1}y}\rangle \oplus C^{k,\alpha}_{\omega_m}([0,\infty[).
\end{aligned} \tag{3.19}$$

Observe that $\tilde{\mathrm{L}}_a$ maps $\tilde{\mathcal{A}}^{2,\alpha}_{\underline{\omega},0}$ into $\tilde{\mathcal{A}}^{0,\alpha}_{\underline{\omega},*}$.

**Theorem 3.2.6**

*For all $\underline{\omega}$ satisfying (3.18), $\tilde{K}_{a,0}$ defines a bounded linear map from $\tilde{\mathcal{A}}^{0,\alpha}_{\underline{\omega},*}$ into $\tilde{\mathcal{A}}^{2,\alpha}_{\underline{\omega},0}$ such that*

$$\tilde{\mathrm{L}}_a\tilde{K}_{a,0} = Id.$$

*In particular, $\tilde{L}_a$ defines a linear isomorphism from $\tilde{\mathcal{A}}^{2,\alpha}_{\underline{\omega},0}$ into $\tilde{\mathcal{A}}^{0,\alpha}_{\underline{\omega},*}$ with inverse $\tilde{K}_{a,0}$.*

**Proof:** Indeed, we verify that, for $1 \leq i \leq m-1$,

$$\tilde{\mathrm{K}}_{a,0}[e^{-\omega_i y}] = -\frac{1}{(\omega_i^2-a^2)}e^{-ay} + \frac{1}{(\omega_i^2-a^2)}e^{-\omega_i y}.$$

On the other hand, for $u \in C^{0,\alpha}_{\omega_m}([0,\infty[)$,

$$\tilde{\mathrm{K}}_{a,0}[u] = v_1 + v_2,$$

where

$$\begin{aligned}
v_1(y) &:= -\frac{1}{a}e^{-ay}\int_0^\infty \sinh(az)u(z)dz\ \text{and}\\
v_2(y) &:= -\frac{1}{a}\int_y^\infty \sinh(a(z-y))u(z)dz.
\end{aligned}$$

We verify that

$$\begin{aligned}
\left|-\frac{1}{a}\int_0^\infty \sinh(az)u(z)dz\right| &\leq \frac{1}{a(\omega_m-a)}\|u\|'_{C^0_{\omega_m}}\ \text{and}\\
\|v_2\|'_{C^0_{\omega_m}} &\leq \frac{1}{a(\omega_m-a)}\|u\|'_{C^0_{\omega_m}}.
\end{aligned}$$

Finally, since $\tilde{\mathrm{L}}_a v_2 = u$,

$$\|v_2\|_{C^{2,\alpha}_{\omega_m}} \leq C\|u\|_{C^{0,\alpha}_{\omega_m}}$$

for a suitable constant $C$, and the result follows. $\square$

**3.3 - The two-dimensional linear problem.** For all $m\in\mathbb{Z}$, define $\Pi_m:L^2(S^1)\rightarrow\mathbb{R}$ and $\mathrm{I}_m:\mathbb{R}\rightarrow L^2(S^1)$ by

$$\begin{aligned}\Pi_m[u]&:=\frac{1}{2\pi}\int_0^{2\pi}u(x)e^{-imx}dx\ \text{and}\\ \mathrm{I}_m[\lambda]&:=\lambda e^{imx}.\end{aligned}\tag{3.20}$$

The operators $\Pi_m[u]$ and $\mathrm{I}_m\Pi_m[u]$ yield respectively the $m$'th Fourier coefficient of $u$ and the orthogonal projection of $u$ onto the $m$'th Fourier mode with respect to the standard $L^2$ inner product of $S^1$. We denote also by $\Pi_m:L^2(\mathbb{S}^1\times\mathbb{R})\rightarrow L^2(\mathbb{R})$ and $\mathrm{I}_m:L^2(\mathbb{R})\rightarrow L^2(\mathbb{S}^1\times\mathbb{R})$ the natural extensions of these operators. For all $m$, denote

$$\begin{aligned}\mathrm{P}_m&:=\sum_{|n|<m}\mathrm{I}_n\Pi_n\ \text{and}\\ \mathrm{P}_m^\perp&:=\mathrm{Id}-\mathrm{P}_m.\end{aligned}\tag{3.21}$$

As before, let $E$ be a Banach space and, for all $\omega$ and for all $(k,\alpha)$, define the $\omega$-weighted $C^{k,\alpha}$-norm for $k$-times differentiable functions $u:S^1\times\mathbb{R}\rightarrow\mathbb{R}$ by

$$\|u\|_{C^{k,\alpha}_\omega}:=\|ue^{\omega\langle y\rangle}\|_{C^{k,\alpha}}.\tag{3.22}$$

Let $\underline{\omega}:=(\omega_0,\cdots,\omega_m)$ be a vector of $(m+1)$ real weights. For all $(k,\alpha)$, define the *$\underline{\omega}$-weighted $C^{k,\alpha}$-norm* by

$$\|u\|_{C^{k,\alpha}_{\underline{\omega}}}:=\sum_{|n|<m}\|\Pi_n[u]\|_{C^{k,\alpha}_{\omega_{|n|}}}+\|P_m^\perp[u]\|_{C^{k,\alpha}_{\omega_m}},\tag{3.23}$$

and let $C^{k,\alpha}_{\underline{\omega}}(S^1\times\mathbb{R}^1,E)$ denote the Banach space of $k$-times differentiable functions $u:S^1\times\mathbb{R}\rightarrow E$ for which this norm is finite.

Consider now the second-order, linear partial differential operator

$$\mathrm{L}_a u:=\partial_x^2u+\partial_y^2u-a^2u.\tag{3.24}$$

By classical Fourier analysis (see [4]), its Green's function is

$$K_a(x,y):=\sum_{m\in\mathbb{Z}}K_{a,m}(x,y),\tag{3.25}$$

where, for all $m$,

$$K_{a,m}(x,y):=-\frac{1}{4\pi\sqrt{m^2+a^2}}e^{imx}e^{-\sqrt{m^2+a^2}|y|},\tag{3.26}$$

and its Green's operator is

$$\mathrm{K}_a[u](x,y):=\int_0^{2\pi}\int_{-\infty}^{\infty}K_a(x-\xi,y-\eta)u(\xi,\eta)d\xi d\eta.\tag{3.27}$$

For all $m\geq0$, denote

$$K_{a,m}^\perp(x,y):=\sum_{|n|\geq m}K_{a,n}(x,y),\tag{3.28}$$

and denote by $\mathrm{K}_{a,m}^\perp$ the integral operator that it defines. Trivially,

$$\mathrm{K}_{a,m}=\mathrm{K}_{a,m}^\perp+\sum_{|n|<m}\mathrm{I}_n\tilde{\mathrm{K}}_{\sqrt{a^2+n^2}}\Pi_n,\tag{3.29}$$

where, for each $n$, $\tilde{\mathrm{K}}_{\sqrt{a^2+n^2}}$ is the operator defined in Section 3.2. It follows from elementary Fourier analysis that the function $\mathrm{K}_{a,m}^\perp$ is locally of class $L^2$ and therefore also locally of class $L^1$. In addition, since the sum (3.28) is close to being a geometric series, we obtain

**Lemma 3.3.1**

*For all $m$ and for all $Y>0$, there exists $B>0$ such that, for all $|y|\geq Y$,*

$$K^{\perp}_{a,m}(x,y)\leq Be^{-\sqrt{m^2+a^2}|y|}. \tag{3.30}$$

This in turn yields

**Lemma 3.3.2**

*For all $m$, and for all $|\omega|<\sqrt{m^2+a^2}$, $K^{\perp}_{a,m}$ defines a bounded linear map from $C^{0,\alpha}_{\omega}(S^1\times\mathbb{R},E)$ into $C^{2,\alpha}_{\omega}(S^1\times\mathbb{R},E)$ such that*

$$L_aK^{\perp}_{a,m}=Id-\sum_{|n|<m}I_n\Pi_n. \tag{3.31}$$

**Proof:** Indeed, using Lemma 3.3.1, we show as in Lemma 3.2.2 that $\mathrm{K}^{\perp}_{a,m}$ defines a bounded linear map from $C^0_{\omega}(S^1\times\mathbb{R},E)$ into $C^0_{\omega}(S^1\times\mathbb{R},E)$. Differentiating under the integral, we verify that, for all $u\in C^0_{\omega}(S^1\times\mathbb{R},E)$, $\mathrm{K}^{\perp}_{a,m}[u]$ is twice differentiable and satisfies

$$\mathrm{L}_a\mathrm{K}^{\perp}_{a,m}[u]=u-\sum_{|n|<m}\mathrm{I}_n\Pi_n[u].$$

By the interior Schauder estimates (see, for example, Chapter 6 of [10]), it follows that $\mathrm{K}^{\perp}_{a,m}$ defines a bounded linear map from $C^{0,\alpha}_{\omega}(S^1\times\mathbb{R},E)$ into $C^{2,\alpha}_{\omega}(S^1\times\mathbb{R},E)$. This completes the proof. $\square$

**Lemma 3.3.3**

*For all $\underline{\omega}$ such that, for all $0\leq i\leq m$,*

$$|\omega_i|<\sqrt{a^2+m^2}, \tag{3.32}$$

*$K_a$ defines a bounded linear map from $C^{0,\alpha}_{\underline{\omega}}(S^1\times\mathbb{R},E)$ into $C^{2,\alpha}_{\underline{\omega}}(S^1\times\mathbb{R},E)$ such that*

$$L_aK_a=Id.$$

**Proof:** This follows from Lemmas 3.2.3 and 3.3.2 together with (3.29). $\square$

This yields

**Theorem 3.3.4**

*For all $a>0$ and for all $\underline{\omega}$ such that, for all $0\leq i\leq m$,*

$$0\leq\omega_i<\sqrt{a^2+i^2}, \tag{3.33}$$

*$L_a$ defines a linear isomorphism from $C^{2,\alpha}_{\underline{\omega}}(S^1\times\mathbb{R},E)$ into $C^{0,\alpha}_{\underline{\omega}}(S^1\times\mathbb{R},E)$ with inverse $K_a$.*

**Proof:** Indeed, by Lemma 3.3.3, $\mathrm{K}_a$ defines a right inverse of $\mathrm{L}_a$. By the maximum principle, $\mathrm{L}_a$ is injective, and the result follows. $\square$

We now consider the Dirichlet problem. For all $\underline{\omega}$ and for all $(k,\alpha)$, the Banach space $C^{k,\alpha}_{\underline{\omega}}(S^1\times[0,\infty[,E)$ is defined in the natural manner and its closed subspace $C^{k,\alpha}_{\underline{\omega},0}(S^1\times[0,\infty[,E)$ is defined to consist of those functions which vanish along the boundary $S^1\times\{0\}$. We verify that the Green's operator of $\mathrm{L}_a$ for the Dirichlet problem is

$$\begin{aligned}\mathrm{K}_{a,0}[u](x,y):=&\int_0^{2\pi}\int_0^{\infty}K_a(x-\xi,y-\eta)u(\xi,\eta)d\xi d\eta\\ &-\int_0^{2\pi}\int_{-\infty}^{0}K_a(x-\xi,y-\eta)u(-\xi,\eta)d\xi d\eta.\end{aligned} \tag{3.34}$$

Proceeding as before, we obtain

**Theorem 3.3.5**

*For all $a > 0$ and for all $\underline{\omega}$ such that, for all $0 \leq i \leq m$,*

$$0 \leq \omega_i < \sqrt{a^2 + i^2}, \tag{3.35}$$

*$L_a$ defines a linear isomorphism from $C^{2,\alpha}_{\underline{\omega},0}(S^1 \times \mathbb{R}, E)$ into $C^{0,\alpha}_{\underline{\omega}}(S^1 \times \mathbb{R}, E)$ with inverse $K_{a,0}$.*

We now consider the case where the weight in each Fourier mode is permitted to be greater than the corresponding constant term. Recall the subsets $\mathcal{M}_0$ and $\mathcal{M}$ of $\mathbb{R}^2$ defined in Section 3.1. For $(\lambda, \mu) \in \mathcal{M}$, define

$$u_{(\lambda,\mu)}(x, y) := e^{i\lambda x} e^{-\mu y}. \tag{3.36}$$

For all $\omega > 0$ and for all $(k, \alpha)$, define

$$\begin{aligned} \mathcal{A}^{k,\alpha}_{\omega} &:= \langle u_{(\lambda,\mu)} \mid (\lambda, \mu) \in \mathcal{M},\ \mu < \omega \rangle \oplus C^{k,\alpha}_{\omega}(S^1 \times [0, \infty[) \\ \mathcal{A}^{k,\alpha}_{\omega,0} &:= \left\{ u \in \mathcal{A}^{k,\alpha}_{\omega} \mid u(x, 0) = 0\ \forall x \right\},\ \text{and} \\ \mathcal{A}^{k,\alpha}_{\omega,*} &:= \langle u_{(\lambda,\mu)} \mid (\lambda, \mu) \in \mathcal{M} \setminus \mathcal{M}_0,\ \mu < \omega \rangle \oplus C^{k,\alpha}_{\omega}(S^1 \times [0, \infty[). \end{aligned} \tag{3.37}$$

We now define

$$\mathcal{A}^{k,\alpha} := \bigcap_{\omega > 0} \mathcal{A}^{k,\alpha}_{\omega}, \tag{3.38}$$

and

$$\mathcal{A} := \bigcap_{k,\alpha} \mathcal{A}^{k,\alpha}. \tag{3.39}$$

Since $\mathcal{A}^{k,\alpha}$ and $\mathcal{A}$ are defined as intersections of families of Banach spaces, they carry natural Frechet structures given by the families of all norms of these spaces. Observe, in addition, that $\mathcal{A}$ is none other than the space defined in (3.7). In particular, for all $(k, \alpha)$, every function $u \in \mathcal{A}^{k,\alpha}$ has a unique asymptotic expansion of the form

$$u \sim \sum_{(\lambda,\mu) \in \mathcal{M}} a_{\lambda,\mu} u_{\lambda,\mu}, \tag{3.40}$$

where, for all $(\lambda, \mu) \in \mathcal{M}$, $a_{\lambda,\mu}$ is a complex coefficient. The derivatives of all such functions up to and including order $k$ also have unique asymptotic expansions of the same form, which are determined by differentiating (3.40) term by term. Finally, combining Theorem 3.2.6 and Lemma 3.3.2 yields

**Theorem 3.3.6**

*For all $a > 0$ and for all $\omega > 0$, $K_{a,0}$ defines a bounded linear map from $\mathcal{A}^{0,\alpha}_{\omega,*}$ into $\mathcal{A}^{2,\alpha}_{\omega,0}$ such that*

$$L_a K_{a,0} = Id.$$

*In particular, $L_a$ defines a linear isomorphism from $\mathcal{A}^{2,\alpha}_{\omega,0}$ into $\mathcal{A}^{0,\alpha}_{\omega,*}$ with inverse $K_{a,0}$.*

**3.4 - Products.** Let $E_1$, $E_2$ and $F$ be Banach spaces. Let $b : E_1 \oplus E_2 \to F$ be a bounded bilinear map, and define the operator

$$\mathrm{B}[u, v](x, y) := b(u(x, y), v(x, y)). \tag{3.41}$$

Let $X$ be a manifold locally isometric to $\mathbb{R}^d \times [0, \infty[$ for some $d \geq 0$. Recall (see Appendix A) that B defines a bounded bilinear map from $C^{k,\alpha}(X, E_1) \oplus C^{k,\alpha}(X, E_2)$ into $C^{k,\alpha}(X, F)$. We now extend this property to weighted spaces.

**Lemma 3.4.1**

*(1) If $\omega_1 \geq \omega_2$, then the canonical embedding $\mathrm{J}_{\omega_2,\omega_2} : C^{k,\alpha}_{\omega_1}(X, E) \to C^{k,\alpha}_{\omega_2}(X, E)$ is continuous.*

*(2) If $\omega_i \geq \omega$ for all $0 \leq i \leq m$, then the canonical embedding $\mathrm{J}_{\omega,\underline{\omega}} : C^{k,\alpha}_{\underline{\omega}}(X, E) \to C^{k,\alpha}_{\omega}(X, E)$ is continuous.*

*(3) If $\omega_i \leq \omega$ for all $0 \leq i \leq m$, then the canonical embedding $\mathrm{J}_{\underline{\omega},\omega} : C^{k,\alpha}_{\omega}(X, E) \to C^{k,\alpha}_{\underline{\omega}}(X, E)$ is continuous.*

**Proof:** It suffices to prove (1) as the proofs of (2) and (3) are almost identical. Since $e^{(\omega_2 - \omega_1)\langle x \rangle}$ is an element of $C^{k,\alpha}(\mathbb{R})$, the operator $\mu_{\omega_1 - \omega_2}$ defines a bounded linear map from $C^{k,\alpha}(X, E)$ to itself. Since $\mathrm{J}_{\omega_2,\omega_1} = \mu_{\omega_2} \mu_{\omega_1 - \omega_2} \mu_{-\omega_1}$, the result follows. $\square$

**Lemma 3.4.2**

*For all $\omega_1,\omega_2,\omega_3\in\mathbb{R}$ such that $\omega_1+\omega_2\geq\omega_3$, $B$ defines a bounded bilinear map from $C^{k,\alpha}_{\omega_1}(X,E)\oplus C^{k,\alpha}_{\omega_2}(X,E)$ into $C^{k,\alpha}_{\omega_3}(X,F)$.*

**Proof:** Indeed, $\mathrm{B}=\mathrm{J}_{\omega_3,(\omega_1+\omega_2)}\mu_{(\omega_1+\omega_2)}\mathrm{B}(\mu_{-\omega_1}\cdot,\mu_{-\omega_2}\cdot)$, and the result follows. □

**Lemma 3.4.3**

*If $\underline{\omega}$ is such that, for all $0\leq i,j\leq m$,*

$$\omega_i+\omega_j\geq\omega_{\min(i+j,m)},\tag{3.42}$$

*then $B$ defines a continuous bilinear map from $C^{k,\alpha}_{\underline{\omega}}(X,E_1)\oplus C^{k,\alpha}_{\underline{\omega}}(X,E_2)$ into $C^{k,\alpha}_{\underline{\omega}}(X,F)$.*

**Proof:** Indeed, for $u\in C^{k,\alpha}_{\underline{\omega}}(X,E_1)$ and $v\in C^{k,\alpha}_{\underline{\omega}}(X,E_2)$,

$$\mathrm{B}[u,v]=\mathrm{B}_1[u,v]+\mathrm{B}_2[u,v]+\mathrm{B}_3[u,v],$$

where

$$\begin{aligned}\mathrm{B}_1[u,v]&:=\sum_{|i|,|j|<m}\mathrm{I}_{i+j}\mathrm{B}[\Pi_iu,\Pi_jv],\\ \mathrm{B}_2[u,v]&:=\sum_{|i|<m}\mathrm{B}[\mathrm{I}_i\Pi_iu,\mathrm{P}^\perp_mv]+\sum_{|i|<m}\mathrm{B}[\mathrm{P}^\perp_mu,\mathrm{I}_i\Pi_iv],\ \text{and}\\ \mathrm{B}_3[u,v]&:=\mathrm{B}[\mathrm{P}^\perp_mu,\mathrm{P}^\perp_mv].\end{aligned}$$

By definition, for all $|i|<m$, $\Pi_i$ defines a continuous linear map from $C^{k,\alpha}_{\underline{\omega}}(X,E_i)$ into $C^{k,\alpha}_{\omega_i}(\mathbb{R},E_i)$ and continuity of $\mathrm{B}_1$ follows by Lemma 3.4.2. By Item (2) of Lemma 3.4.1, for all $0\leq|i|<m$, $\mathrm{I}_i\Pi_i$ defines a continuous linear map from $C^{k,\alpha}_{\underline{\omega}}(X,E_i)$ into $C^{k,\alpha}_{\omega_0}(X,E_i)$. By definition, $\mathrm{P}^\perp_m$ defines a continuous linear map from $C^{k,\alpha}_{\underline{\omega}}(X,E_i)$ into $C^{k,\alpha}_{\omega_m}(X,E_i)$. By Lemma 3.4.2 and (3.42), B defines a continuous bilinear map from $C^{k,\alpha}_{\omega_0}(X,E_1)\oplus C^{k,\alpha}_{\omega_m}(X,E_2)$ into $C^{k,\alpha}_{\omega_m}(X,F)$, and continuity of $\mathrm{B}_2$ follows by Item (3) of Lemma 3.4.1. Finally, by Lemma 3.4.2, B defines a continuous bilinear map from $C^{k,\alpha}_{\omega_m}(X,E_1)\oplus C^{k,\alpha}_{\omega_m}(X,E_2)$ into $C^{k,\alpha}_{\omega_m}(X,F)$, and continuity of $\mathrm{B}_3$ follows by Item (3) of Lemma 3.4.1. This completes the proof. □

A similar reasoning yields

**Lemma 3.4.4**

*For all $\omega>0$, $B$ defines a continuous bilinear map from $\mathcal{A}^{k,\alpha}_\omega\oplus\mathcal{A}^{k,\alpha}_\omega$ into $\mathcal{A}^{k,\alpha}_{\omega',*}$, where*

$$\omega'=min(2\omega,\omega+a).\tag{3.43}$$

**Remark 3.4.8.** In particular, the operator B defines a bilinear map which is continuous in the Frechet sense from $\mathcal{A}\oplus\mathcal{A}$ into $\mathcal{A}$. Furthermore, the argument used in the proof of Lemma 3.4.3 shows that, for any two $u,v\in\mathcal{A}$, the asymptotic series of the product $uv$ is obtained by multiplying term by term the asymptotic series of each of $u$ and $v$.

**3.5 - Non-linear operators.** Let $E$ and $F$ now be finite-dimensional vector spaces. Let $\Omega$ be an open subset of $E$. Let $\Phi:\Omega\to F$ be a smooth function. Let $\mathrm{C}_\Phi$ denote the operator of composition by $\Phi$, that is

$$\mathrm{C}_\Phi[u]:=\Phi\circ u.\tag{3.44}$$

Let $\underline{\omega}:=(\omega_0,\cdots,\omega_m)$ be a finite vector of real weights such that, for all $0\leq i,j\leq m$,

$$\omega_i+\omega_j\geq\omega_{\min(i+j,m)}.$$

By Lemma 3.4.3, upon rescaling the norm of $C^{k,\alpha}_{\underline{\omega}}(\mathbb{R}\times S^1)$ if necessary, we may suppose that, for $u,v\in C^{k,\alpha}_{\underline{\omega}}(\mathbb{R}\times S^1)$,

$$\|uv\|_{C^{k,\alpha}_{\underline{\omega}}}\leq\|u\|_{C^{k,\alpha}_{\underline{\omega}}}\|v\|_{C^{k,\alpha}_{\underline{\omega}}}.$$

It then follows that if $\Phi$ is analytic with radius of convergence $R$ about 0, then $\mathrm{C}_\Phi$ is also analytic with the same radius of convergence about 0. This in itself would be sufficient for our purposes since the functions of interest to us are all analytic. However, for completeness, we consider also the case where $\Phi$ is an arbitrary smooth function. To this end, define

$$\mathcal{O}^{k,\alpha}_{\underline{\omega}}(S^1\times\mathbb{R},\Omega) := \bigcup_{K\subseteq\Omega}\mathcal{O}^{k,\alpha}_{\underline{\omega}}(S^1\times\mathbb{R},K), \tag{3.45}$$

where $K$ varies over all compact subsets of $\Omega$ and, for all such $K$,

$$\mathcal{O}^{k,\alpha}_{\underline{\omega}}(S^1\times\mathbb{R},K) := \left\{f\in C^{k,\alpha}_{\underline{\omega}}(S^1\times\mathbb{R})\ |\ \mathrm{Im}(f)\subseteq K\right\}. \tag{3.46}$$

Observe that, if $\omega_0\geq 0$, then this set is open in $C^{k,\alpha}_{\underline{\omega}}(S^1\times\mathbb{R},E)$. It is the natural domain over which $\mathrm{C}_\Phi$ is defined.

**Lemma 3.5.1**

*If $\omega_0>0$, then for all $(k,\alpha)$, $C_\Phi$ defines a continuous function from $\mathcal{O}^{k,\alpha}_{\underline{\omega}}(S^1\times\mathbb{R},\Omega)$ into $C^{k,\alpha}_{\underline{\omega}}(S^1\times\mathbb{R})$.*

**Proof:** It suffices to prove the case where $E=F=\mathbb{R}$. The general case is similar. Let $M>0$ be such that

$$M\omega_0>\omega_m.$$

There exists a polynomial $P$ of order $(M-1)$ and a smooth function $\Psi:\Omega\rightarrow\mathbb{R}$ such that, for all $x\in\Omega$,

$$\Phi(x)=P(x)+\Psi(x)x^M.$$

In particular, for all $u$,

$$\mathrm{C}_\Phi[u]=\mathrm{C}_P[u]+\mathrm{C}_\Psi[u]u^M. \tag{3.47}$$

Since $P(0)=\Phi(0)=0$, by Lemma 3.4.3, $\mathrm{C}_P$ defines a continuous function from $C^{k,\alpha}_{\underline{\omega}}(S^1\times\mathbb{R})$ to itself. Recall now (see Appendix A) that $\mathrm{C}_\Psi$ defines a continuous function from $C^{k,\alpha}(S^1\times\mathbb{R})$ to itself. Furthermore, by Lemma 3.4.2, multiplication defines a continuous function from $C^{k,\alpha}(S^1\times\mathbb{R})\oplus C^{k,\alpha}_{\omega_0}(S^1\times\mathbb{R})^M$ into $C^{k,\alpha}_{M\omega_0}(S^1\times\mathbb{R})$. Thus, since

$$\mathrm{C}_\Psi[u]u^M=\mathrm{J}_{\underline{\omega},M\omega_0}\left[(\mathrm{C}_\Psi\circ\mathrm{J}_{0,\underline{\omega}})[u])\mathrm{J}_{\omega_0,\underline{\omega}}[u]^M\right],$$

the result follows by Lemma 3.4.1. $\square$

**Lemma 3.5.2**

*If $\omega_0>0$, then for all $(k,\alpha)$, $C_\Phi$ defines a differentiable function from $\mathcal{O}^{k,\alpha}_{\underline{\omega}}(S^1\times\mathbb{R})$ into $C^{k,\alpha}_{\underline{\omega}}(S^1\times\mathbb{R})$ with derivative at $u\in\mathcal{O}^{k,\alpha}_{\underline{\omega}}(S^1\times\mathbb{R})$ given by*

$$\big(DC_\Phi[u]v\big)(x,y)=C_{D\Phi}[u](x,y)v(x,y). \tag{3.48}$$

**Proof:** Choose $u\in\mathcal{O}^{k,\alpha}_{\underline{\omega}}(S^1\times\mathbb{R})$. For sufficiently small $v\in C^{k,\alpha}_{\underline{\omega}}(S^1\times\mathbb{R})$, denote

$$\Psi[v]:=\mathrm{C}_{D\Phi}[u+v]-\mathrm{C}_{D\Phi}[u].$$

By the fundamental theorem of calculus, for all $(x,y)$,

$$\mathrm{C}_\Phi[u+v](x,y)-\mathrm{C}_\Phi[u](x,y)-\mathrm{C}_{D\Phi}[u](x,y)v(x,y)=\int_0^1\Psi[sv](x,y)dsv(x,y).$$

Choose $\epsilon>0$. By Lemma 3.5.1, $\Psi$ is continuous. Since $\Psi[0]=0$, there therefore exists $\delta>0$ such that, for $\|v\|_{C^{k,\alpha}_{\underline{\omega}}}<\delta$,

$$\|\Psi[v]\|_{C^{k,\alpha}_{\underline{\omega}}}<\epsilon.$$

Furthermore, by continuity, the function $s\mapsto\Psi[sv]$ is integrable as a function taking values in the Banach space $C^{k,\alpha}_{\underline{\omega}}(S^1\times\mathbb{R})$. It follows by convexity of the norm that, for $\|v\|_{C^{k,\alpha}_{\underline{\omega}}}<\delta$,

$$\left\|\int_0^1\Psi[sv]ds\right\|_{C^{k,\alpha}_{\underline{\omega}}}\leq\int_0^1\|\Psi[sv]\|_{C^{k,\alpha}_{\underline{\omega}}}ds<\epsilon.$$

Consequently, for $\|v\|_{C^{k,\alpha}_{\underline{\omega}}}<\delta$,

$$\left\|\mathrm{C}_\Phi[u+v]-\mathrm{C}_\Phi[u]-\mathrm{C}_{D\Phi}[u]v\right\|_{C^{k,\alpha}_{\underline{\omega}}}<\epsilon\|v\|_{C^{k,\alpha}_{\underline{\omega}}}.$$

Since $\epsilon>0$ is arbitrary, the result follows. □

Applying Lemma 3.5.2 inductively yields

**Theorem 3.5.3**

*If $\omega_0>0$ and if $\Phi[0]=0$ then, for all $(k,\alpha)$, $\mathrm{C}_\Phi$ defines a smooth function from $\mathcal{O}^{k,\alpha}_{\underline{\omega}}(S^1\times\mathbb{R},\Omega)$ into $C^{k,\alpha}_{\underline{\omega}}(\mathbb{R}\times S^1,\Omega)$.*

Finally, for $\omega>0$ and for all $(k,\alpha)$. Define

$$\mathcal{U}^{k,\alpha}_\omega(\Omega):=\underset{K\subseteq\Omega}{\cup}\mathcal{U}^{k,\alpha}_\omega(S^1\times\mathbb{R},K),\tag{3.49}$$

where $K$ varies over every compact subset of $\Omega$ and, for all such $K$,

$$\mathcal{U}^{k,\alpha}_\omega(S^1\times\mathbb{R},K):=\left\{u\in\mathcal{A}^{k,\alpha}_\omega\ \mid\ \mathrm{Im}(u)\subseteq K\right\}.\tag{3.50}$$

Repeating the proof of Theorem 3.5.3 yields

**Lemma 3.5.4**

*(1) If $\omega>0$ and if $\Phi(0)=0$, then $\mathrm{C}_\Phi$ defines a smooth function from $\mathcal{U}^{k,\alpha}_\omega(\Omega)$ into $\mathcal{A}^{k,\alpha}_\omega$.*

*(2) If, in addition $D\Phi(0)=0$, then $\mathrm{C}_\Phi$ defines a smooth function from $\mathcal{U}^{k,\alpha}_\omega(\Omega)$ into $\mathcal{A}^{k,\alpha}_{\omega+a,*}$.*

**Remark 3.5.9.** When $\Phi(0)=0$ and $D\Phi(0)=0$, the MacLaurin series of this function involves only terms of quadratic or higher order, so that the stronger result of (2) follows from Lemma 3.4.4.

**Remark 3.5.10.** In particular, when $\Phi(0)=0$, $\mathrm{C}_\Phi$ defines a map from a neighbourhood of 0 in $\mathcal{A}$ into $\mathcal{A}$ which is smooth in the Frechet sense. Furthermore, bearing in mind (3.47), we see that, for any $u\in\mathcal{A}$, the asymptotic series of $\mathrm{C}_\Phi[u]$ is determined by substituting the asymptotic series of $u$ formally into the MacLaurin series of $\Phi$.

**3.6 - The Dirichlet solution operator.** Recall that $\mathrm{R}:C^{k,\alpha}_{\underline{\omega}}(S^1\times[0,\infty[)\rightarrow C^{k,\alpha}(S^1)$ denotes the operator of restriction onto $S^1\times\{0\}$.

**Theorem 3.6.1**

*If $a>0$ and if $\underline{\omega}:=(\omega_0,\cdots,\omega_m)$ is such that, for all $0\leq i\leq m$,*

$$0<\omega_i<\sqrt{a^2+i^2},\tag{3.51}$$

*and, for all $0\leq i\leq m$,*

$$\omega_i+\omega_j\geq\omega_{\min(i+j,m)},\tag{3.52}$$

*then there exists a neighbourhood $U$ of 0 in $C^{2,\alpha}(S^1)$ and a smooth map $S:U\rightarrow C^{2,\alpha}_{\underline{\omega}}([0,\infty[\times S^1)$ such that, for all $v\in U$,*

*(1) $RS[v]=v$ and*

*(2) $PS[v]=0$.*

*Furthermore, upon reducing $U$ if necessary, $S$ is unique.*

**Remark 3.6.11.** We show in the usual manner that for all $k+\beta\geq 2+\alpha$, S maps $U\cap C^{k,\beta}(S^1)$ smoothly into $C^{k,\beta}_{\underline{\omega}}([0,\infty[\times S^1)$. In addition, the same function S maps $U\cap C^{k,\beta}(S^1)$ smoothly into $C^{k,\beta}_{\underline{\omega}'}([0,\infty[\times S^1)$ for any other $\underline{\omega}'$ satisfying the hypotheses of Theorem 3.6.1.

**Proof:** Indeed, by Theorem 3.3.5, $(\mathrm{R},\mathrm{L}_a)$ defines a linear isomorphism from $C^{2,\alpha}_{\underline{\omega}}(S^1\times[0,\infty[)$ into the product $C^{2,\alpha}(S^1)\times C^{0,\alpha}_{\underline{\omega}}(S^1\times[0,\infty[)$. By the inverse function theorem, there exists a neighbourhood $\tilde{U}$ of 0 in $C^{2,\alpha}(S^1)\times C^{0,\alpha}_{\underline{\omega}}(S^1\times[0,\infty[)$ and a smooth map $\tilde{\mathrm{S}}$ such that, for all $(v,w)\in\tilde{U}$, $(\mathrm{R}\tilde{\mathrm{S}}[v,w],\mathrm{P}\tilde{\mathrm{S}}[v,w])=(v,w)$. Existence follows upon setting $\mathrm{S}[v]:=\tilde{\mathrm{S}}[v,0]$ and uniqueness follows by the uniqueness part of the inverse function theorem. □

For all $\omega$ and for all $(k,\alpha)$, let $\mathcal{A}^{k,\alpha}_{\omega}$ be the Banach space defined in (3.38).

**Theorem 3.6.2**

*Let $U$ and $S$ be as in Theorem 3.6.1. For all $\omega>0$, $S$ defines a smooth map from $U$ into $\mathcal{A}^{2,\alpha}_{\omega}$.*

**Remark 3.6.12.** A suitable refinement of Theorem 3.6.1 shows that, for all $k+\beta\geq 2+\alpha$ and for all $\omega$, S defines a smooth map from $U\cap C^{k,\beta}(S^1)$ into $\mathcal{A}^{k,\beta}_{\omega}$, so that S defines a Frechet smooth map from $U\cap C^{\infty}(S^1)$ into $\mathcal{A}$. In particular, for smooth initial data, the asymptotic series constructed by Theorem 3.6.2 are differentiable to all orders.

**Proof:** Indeed, for $\omega<a$, the result follows by Theorem 3.6.1. Let $W\subseteq]0,\infty[$ be the set of all weights for which the assertion is true. Since, for $\omega'<\omega$, the canonical embedding $\mathrm{J}:\mathcal{A}^{2,\alpha}_{\omega}\to\mathcal{A}^{2,\alpha}_{\omega'}$ is a bounded linear map, it follows that $W$ is an interval with lower extremity 0. Let $\omega_0:=\sup(W)$ and suppose that $\omega_0<\infty$. Denote $\omega:=\omega_0-a/2$. By Item (2) of Lemma 3.5.4, for $u\in\mathcal{U}^{2,\alpha}_{\omega}(\Omega)$,

$$\mathrm{C}_F(\mathrm{S}[u],D\mathrm{S}[u],D^2\mathrm{S}[u])\in\mathcal{A}^{0,\alpha}_{\omega+a,*}.$$

By Theorem 3.3.6, $(\mathrm{R},\mathrm{L}_a)$ defines a linear isomorphism from $\mathcal{A}^{2,\alpha}_{\omega+a}$ into $C^{2,\alpha}(S^1)\oplus\mathcal{A}^{0,\alpha}_{\omega+a,*}$. However, for all $u$,

$$(\mathrm{R},\mathrm{L}_a)\mathrm{S}[u]=(0,\mathrm{C}_F(\mathrm{S}[u],D\mathrm{S}[u],D^2\mathrm{S}[u])),$$

so that S maps $U$ smoothly into $\mathcal{A}^{2,\alpha}_{\omega+a}$. This is absurd, by definition of $\omega_0$. It follows that $\omega_0=\infty$, and this completes the proof. □

It remains only to verify Part $(A)$ of Theorem 3.1.1. We first require the following technical result.

**Lemma 3.6.3**

*Let $u:S^1\times[0,\infty[\to\mathbb{C}$ be a smooth function whose derivatives to all orders are bounded over $S^1\times[0,\infty[$. For $\lambda>0$, if $u(x,y)=O(e^{-\lambda y})$ then, for all $k\geqslant 1$ and for all $\epsilon>0$,*

$$\|D^k u(x,y)\|=O(e^{-(\lambda-\epsilon)y}).\tag{3.53}$$

**Proof:** It suffices to prove the result for functions defined over $\mathbb{R}\times[0,\infty[$. Let $\chi\in C^{\infty}_0(\mathbb{R}^2)$ be such that $\chi=1$ over $B_1(0)$ and $\mathrm{Supp}(\chi)\subseteq B_2(0)$. For all $(\xi,\eta)$, define

$$u_{\xi,\eta}(x,y):=u(x,y)\chi(x-\xi,y-\eta).$$

Trivially, $\|u_{\xi,\eta}\|_{C^0}=\mathrm{O}(e^{-\lambda\eta})$ and, for all $k\geqslant 1$, $\|D^k u(\xi,\eta)\|_{C^0}=\mathrm{O}(1)$. By the interpolation inequalities of Chapter 6.8 of [10], for all $0<k<l$,

$$\|D^k u(x,y)\|\leq\|D^k u_{x,y}\|_{C^0}=\mathrm{O}\big(\|u_{x,y}\|^{(l-k)/l}_{C^0}\cdot\|D^l u_{x,y}\|^{k/l}_{C^0}\big)=\mathrm{O}(e^{-(l-k)\lambda y/l}),$$

and the result follows upon chosing $l$ sufficiently large. □

This yields

**Lemma 3.6.4**

*If $u : S^1 \times [0,\infty[\to \mathbb{C}$ is a solution of (3.3) such that, for some $\lambda > 0$,*

$$\begin{aligned} |u(x,y)| &= O(e^{-\lambda y}),\ \text{and} \\ \|D^k u(x,y)\| &= O(1)\ \forall k \geqslant 1, \end{aligned} \tag{3.54}$$

*then $u$ is an element of $\mathcal{A}$.*

**Proof:** Indeed, by Lemma 3.6.3, $u \in C^{2,\alpha}_\omega(S^1 \times [0,\infty[)$ for all $0 < \omega < \lambda$. Upon truncating $u$ to $S^1 \times [Y_0,\infty[$ and reparametrising $y$, we may suppose that $\|u\|_{C^{2,\alpha}_\omega}$ is as small as we wish so that, by Theorem 3.6.1, $u$ is in the image of S and, by uniqueness, $u = \mathrm{SR}[u]$. It then follows by Theorem 3.6.2 and the subsequent remark that $u = \mathrm{SR}[u]$ is an element of $\mathcal{A}$, and this completes the proof. □

We conclude by verifying that Theorem 3.1.1 applies in particular to $k$-ends.

**Lemma 3.6.5**

*If $u : S^1 \times [0,\infty[$ is a bounded solution of (2.43), then $u$ satisfies the hypotheses of Part (A) of Theorem 3.1.1.*

**Proof:** Indeed, by (1.8) and (2.25),

$$u^2 \leq u^2 + u_x^2 = \mathrm{O}(e^{-2\sqrt{1-k}y}).$$

By Theorem 1.2 of [30], for all $k \geqslant 1$,

$$\|D^k u(x,y)\| = \mathrm{O}(1),$$

and the result follows. □

## 4 - The asymptotic geometry of $k$-ends.

**4.1 - Overview.** By the results of the preceding section, every $k$-end has an asymptotic series which is well-defined up to a choice of Weinstein coordinates. We now show how this asymptotic series is used to study the geometry of the end. Consider first a $k$-end $e : S^1 \times [0,\infty[\to \mathbb{H}^3$ which, for ease of presentation, we take for the moment to be of order 1. Up to reparametrisation, we may suppose that $e = \Phi[u]$, where $\Phi$ is as in (2.25) and $u : S^1 \times [0,\infty[\to \mathbb{R}$ solves (2.43). Observe now that, for all $(a,b)^t \in \mathbb{R}^2$,

$$\mathrm{T}[a,b] \circ e = \Phi[u + \sigma[a,b]], \tag{4.1}$$

where $\mathrm{T}[a,b]$ and $\sigma[a,b]$ are respectively as in (2.8) and (2.9). In particular, $u + \sigma[a,b]$ is also a solution of (2.43) which we call the *translate* of $u$ by the vector $(a,b)^t$.

More generally, for $m \in \mathbb{N}$, we denote

$$mS^1 := \mathbb{R}/2\pi m\mathbb{Z}, \tag{4.2}$$

and we define an *abstract $k$-end of order $m$* to be a smooth, exponentially decaying solution $u : mS^1 \times [0,\infty[\to \mathbb{R}$ of the equation

$$ku_{xx} + u_{yy} - (1-k)u = F(u, Du, D^2u), \tag{4.3}$$

where $F$ is a smooth function of its arguments satisfying

(1) $F(0,0,0) = 0$, $DF(0,0,0) = 0$, and

(2) for all $u$ and for all $(a,b)^t \in \mathbb{R}^2$,

$$F(u + \sigma[a,b], Du + D\sigma[a,b], D^2u + D^2\sigma[a,b]) = F(u, Du, D^2u). \tag{4.4}$$

Observe that, upon substituting $x' := x/m$ and $y' := \sqrt{k}y/m$, we recover (3.1) with $a^2 = m^2(1-k)/k$, and the results of Section 3 may therefore be applied. Observe, furthermore, that, by the second condition, whenever $u$ is an abstract $k$-end, so too is the function

$$u + \sigma[a,b], \tag{4.5}$$

for all $(a,b)^t \in \mathbb{R}^2$. As before, we call (4.5) the *translate* of $u$ by the vector $(a,b)^t$.

For $m \in \mathbb{N}$, let $\mathcal{M}_m$ denote the subsemigroup of $\mathbb{R}^2$ generated by

$$\mathcal{M}_{m,0} := \left\{ \frac{1}{m}\big(n, \sqrt{n^2k + m^2(1-k)}\big) \;\middle|\; n \in \mathbb{Z} \right\}. \tag{4.6}$$

For every weight $\omega > 0$, define the Banach space $\mathcal{A}^{k,\alpha}_{m,\omega}$ as in (3.37) with $\mathcal{M}_m$ instead of $\mathcal{M}$ and define the Frechet space $\mathcal{A}_m$ by

$$\mathcal{A}_m := \bigcap_{\omega>0}\bigcap_{k,\alpha} \mathcal{A}^{k,\alpha}_{m,\omega}. \tag{4.7}$$

By Part $(A)$ of Theorem 3.1.1, every abstract $k$-end of order $m$ is an element of $\mathcal{A}_m$. Since $(0,\sqrt{1-k})$ and $(\pm 1, 1)$ are elements of $\mathcal{M}_{m,0}$ for all $m$, we define

$$\begin{aligned} \mathrm{r} &:= \mathrm{a}_{(0,\sqrt{1-k})},\ \text{and} \\ \mathrm{c} &:= (\mathrm{a}_{(1,1)} + \mathrm{a}_{(-1,1)}, i\mathrm{a}_{(1,1)} - i\mathrm{a}_{(-1,1)}), \end{aligned} \tag{4.8}$$

where, for all $(\lambda,\mu)$, $\mathrm{a}_{\lambda,\mu}$ is as in Section 3.1. We call r and c respectively the *radius* and *centroid* operators. In the next section, we will show that every abstract $k$-end has a unique translate whose centroid vanishes. With the preceding discussion in mind, this result, expressed in geometric terms, corresponds exactly to the existence of the Steiner geodesic, thus proving Theorem 1.2.1. The remaining sections of this chapter will then be devoted to deriving asymptotic expressions for various geometric quantities which will be of use in the sequel.

**4.2 - Properties of the radius and centroid operators.** We first observe that the radius and centroid operators are equivariant under rotations, dilatations and translations in the sense that, for every abstract $k$-end $u$, for all $\xi, \eta \in \mathbb{R}$ and for all $(a,b)^t \in \mathbb{R}^2$,

$$\begin{aligned} \mathrm{r}[u(\cdot+\xi,\cdot)] &= \mathrm{r}[u], \\ \mathrm{c}[u(\cdot+\xi,\cdot)] &= \mathrm{R}[\xi]\mathrm{c}[u], \\ \mathrm{r}[u(\cdot,\cdot+\eta)] &= e^{-\sqrt{1-k}\eta}\mathrm{r}[u], \\ \mathrm{c}[u(\cdot,\cdot+\eta)] &= e^{-\eta}\mathrm{c}[u], \\ \mathrm{r}[u+\sigma[a,b]] &= \mathrm{r}[u]\ \text{and} \\ \mathrm{c}[u+\sigma[a,b]] &= \mathrm{T}[a,b]\mathrm{c}[u], \end{aligned} \tag{4.9}$$

where the linear maps $\mathrm{R}[\xi]$ and $\mathrm{T}[a,b]$ are defined by

$$\begin{aligned} \mathrm{R}[\xi](x,y) &:= (\cos(\xi)x + \sin(\xi)y, -\sin(\xi)x + \cos(\xi)y)\ \text{and} \\ \mathrm{T}[a,b](x,y) &:= (x+a, y+b). \end{aligned} \tag{4.10}$$

Recall from Section 2.4 that, for all $u \in \mathcal{A}_m$, $\Phi[u]$ defines a smooth map from $mS^1 \times [0,\infty[$ into $\mathbb{H}^3$. Since

$$u + u_{xx} = \mathrm{r}[u]e^{-\sqrt{1-k}y} + \mathrm{o}(e^{-\sqrt{1-k}y}), \tag{4.11}$$

it follows from (2.32) that $\Phi[u]$ is immersed for sufficiently large $y$ provided that

$$\mathrm{r}[u] \neq 0. \tag{4.12}$$

This property will be important, above all, in the study of families of smooth deformations of $k$-surfaces (see the proof of Lemma 5.2.1, below). In particular, we will verify in Lemma 4.3.3, below, that every $k$-end has non-vanishing radius.

We now turn our attention to the geometric meaning of the centroid operator.

**Lemma 4.2.1**

*For all $u\in\mathcal{A}_m$,*

$$\int_{mS^1}\big(a\cos(x)+b\sin(x)\big)u(x,y)dx=e^{-y}m\pi\langle(a,b),c[u]\rangle_e+O(e^{-\sqrt{4-3k}y}),\tag{4.13}$$

*as $y$ tends to infinity. Furthermore, the coefficient of the remainder term is locally uniformly bounded as $u$ varies in $\mathcal{A}_m$.*

**Proof:** Indeed, for all $y$, the Fourier series of $u(\cdot,y)$ is

$$u(x,y)=\frac{1}{2}\alpha_0(y)+\sum_{n=1}^{\infty}\alpha_n(y)\cos\left(\frac{nx}{m}\right)+\sum_{n=1}^{\infty}\beta_n(y)\sin\left(\frac{nx}{m}\right).$$

Consider now the subset $X$ of $\mathcal{M}_m\setminus\mathcal{M}_{m,0}$ given by

$$X:=\{(\lambda,\mu)\ |\ \lambda=1\},$$

and order the elements of this set by their $y$-components. Observe that the least element of this set is obtained by adding two elements of $\mathcal{M}_0$ so that, by (4.6), its $y$-component is equal to

$$y_0:=\frac{1}{m}\sqrt{n^2k+m^2(1-k)}+\frac{1}{m}\sqrt{(m-n)^2k+m^2(1-k)},\tag{4.14}$$

for some $0\leqslant n\leqslant m$. Observe, however, that the right-hand side of (4.14), being the sum of two convex functions of $n$, is itself also convex. In addition, by symmetry, it attains its unique minimum at $n_0:=m/2$, so that

$$y_0\geqslant\sqrt{4-3k}.$$

Finally, since $u\in\mathcal{A}_m$, it follows that

$$\begin{aligned}\alpha_m(y)&=\mathrm{c}_1[u]e^{-y}+\mathrm{O}(e^{-\sqrt{4-3k}y})\ \text{and}\\ \beta_m(y)&=\mathrm{c}_2[u]e^{-y}+\mathrm{O}(e^{-\sqrt{4-3k}y}).\end{aligned}$$

Furthermore, the coefficients of the remainder terms are locally bounded as $u$ varies in $\mathcal{A}_m$. The result now follows by the Fourier integral formula. $\square$

In particular, this yields Theorem 1.2.1.

**Theorem 1.2.1**

*For every $k$-end $e:mS^1\times[0,\infty[\to\mathbb{H}^3$, there exists a unique unit-speed geodesic $\gamma:\mathbb{R}\to\mathbb{H}^3$ such that*

$$d(\gamma(y),s(y))=O(e^{-y\sqrt{4-3k}})\tag{4.15}$$

*as $y$ tends to infinity.*

**Proof:** Choose an upper half-space parametrisation of $\mathbb{H}^3$ such that the extremity of $e$ coincides with the point at $\infty$. Upon reparamatrising $e$, we may suppose that $e=\Phi[u]$ for some function $u:mS^1\times[0,\infty[\to\mathbb{R}$ solving (2.43). There exists a unique translate of $u$, namely the translate by $-\mathrm{c}[u]$, for which the first term on the right-hand side of (4.13) vanishes. However, by (1.13), for each $y$, the integral on the left-hand side of (4.13) is none other than the inner-product of the vector $(a,b)^t$ with the Steiner curvature centroid of the intersection $\Phi^y[u]$ of this end with the horosphere at height $y$. The vertical half-line passing through the point $(\mathrm{c}[u],1)^t$ is therefore the desired geodesic, and this completes the proof. $\square$

Lemma 4.2.1 also guarantees convergence of the integrals that will be studied in the sequel. We thus conclude this section by establishing elementary conditions for the centroid of a given element of $\mathcal{A}_m$ to vanish. To this end, define

$$\begin{aligned}\mathcal{A}_{m,*}&:=\{u\in\mathcal{A}_m\ |\ \mathrm{a}_{\lambda,\mu}[u]=0\ \forall(\lambda,\mu)\in\mathcal{M}_{m,0}\}\ \text{and}\\ \mathcal{A}_{m,c}&:=\{u\in\mathcal{A}_m\ |\ \mathrm{c}[u]=(0,0)\}.\end{aligned}\tag{4.16}$$

Observe that

$$\mathcal{A}_{m,*}\subseteq\mathcal{A}_{m,c}\subseteq\mathcal{A}_m,$$

and that both $\mathcal{A}_{m,*}$ and $\mathcal{A}_{m,c}$ are ideals of the multiplicative Frechet algebra $\mathcal{A}_m$ which are closed under the action of differentiation.

**Lemma 4.2.2**

*For all $u\in\mathcal{A}_m$,*

$$\begin{aligned} u+u_{xx}&\in\mathcal{A}_{m,c}\ \text{and}\\ u+u_y&\in\mathcal{A}_{m,c}.\end{aligned}\tag{4.17}$$

**Proof:** Indeed, the asymptotic series of $u_{xx}$ and $u_y$ are obtained by differentiating term by term the asymptotic series of $u$. The result follows. $\square$

**4.3 - The geometry of $k$-ends.** Let $u:mS^1\times[0,\infty[\rightarrow\mathbb{R}$ be an abstract $k$-end. By (4.11) we may suppose that $\Phi[u]$ defines a smooth immersion from $mS^1\times[0,\infty[$ into $\mathbb{H}^3$. We now use the notation of Sections 2.4 and 2.5.

**Lemma 4.3.1**

*The form*

$$H[u]dArea[u]-dxdy$$

*is integrable over $mS^1\times[0,\infty[$. Furthermore, its $L^1$-norm is locally uniformly bounded as $u$ varies in $\mathcal{A}_m$.*

**Proof:** Indeed, by (2.30) and (2.41),

$$\mathrm{H}[u]\mathrm{dArea}[u]-dxdy=fdxdy,$$

where $f\in\mathcal{A}_m$. Since

$$f=\mathrm{r}[f]e^{-\sqrt{1-k}y}+\mathrm{o}(e^{-\sqrt{1-k}y}),$$

the result follows. $\square$

**Lemma 4.3.2**

*The length of $\Phi_y[u]$ satisfies*

$$L[u](y)=2\pi m\mathrm{r}[u]e^{-\sqrt{1-k}y}+o(e^{-\sqrt{1-k}y}).\tag{4.18}$$

*Furthermore, the coefficient of the remainder term is locally uniformly bounded as $u$ varies in $\mathcal{A}_m$.*

**Proof:** Indeed, by (2.31) and (4.11) the length element of $\Phi_y[u]$ satisfies

$$\mathrm{dl}[u]=\big(\mathrm{r}[u]e^{-\sqrt{1-k}y}+\mathrm{o}(e^{-\sqrt{1-k}y})\big)dx,$$

and the result follows upon integrating this form over $mS^1$. $\square$

In particular, we obtain

**Lemma 4.3.3**

*Every $k$-end has non-zero radius.*

**Proof:** Choose a half-space parametrisation of $\mathbb{H}^3$. Let $e:mS^1\times[0,\infty[\rightarrow\mathbb{H}^3$ be a $k$-end with respect to the Busemann function $h(x,y,z):=-\log(z)$. Let $\mathrm{I}:=\mathrm{I}[e]$ denote the first fundamental form of $e$ and observe that this metric is complete and of constant curvature equal to $(k-1)$. Let $\tilde{h}$ be a Busemann function of this cusp centred at the unique point at infinity.

Consider now the curve

$$\gamma:[0,\infty[\rightarrow mS^1\times[0,\infty[;y\mapsto(x_0,y),$$

for some fixed point $x_0$ of $mS^1$. By (2.27) and Lemma 3.6.4, the norm of $\partial_t(e\circ\gamma)$ tends to 1 as $t$ tends to infinity. Likewise, the geodesic curvature of $(e\circ\gamma)$ tends to zero as $t$ tends to infinity. It follows that the geodesic curvature of $\gamma$ with respect to I also tends to zero as $t$ tends to infinity. By classical hyperbolic geometry, since $\gamma$ terminates at the unique point at infinity of $mS^1\times[0,\infty[$, for all $\epsilon>0$, upon suitably modifying $e$ and $\tilde{h}$, we may suppose that, for all $t>0$,

$$(1-\epsilon/2)y\leqslant-(\tilde{h}\circ e)(y)\leqslant(1+\epsilon/2)y.$$

We may likewise suppose that, for all $y$, the length of the curve $mS^1\times\{y\}$ with respect to I is less than $\epsilon/2$ so that, for all $(x,y)$,

$$(1-\epsilon)y\leqslant-\tilde{h}(x,y)\leqslant(1+\epsilon)y.$$

In particular, for all $y$, the length of $mS^1\times\{y\}$ with respect to I is bounded below by the length of the level set of $\tilde{h}$ at height $-(1+\epsilon)y$. Consequently, for all $y$,

$$\mathrm{L}[e|_{mS^1\times\{y\}}]\geqslant ce^{-\sqrt{1-k}(1+\epsilon)y},$$

for some constant $c>0$. Observe now that, since $\mathcal{M}_m$ is discrete, the second term on the right-hand side of (4.18) may be replaced with $\mathrm{o}(e^{-\sqrt{1-k}(1+\delta)y})$ for some $\delta>0$. The result now follows upon choosing $\epsilon<\delta$. □

**Lemma 4.3.4**

*When $u$ has strictly positive radius, the geodesic curvature of $\Phi^y[u]$ with respect to the normal $\nu[u]$ satisfies*

$$\kappa^y[u]=-\sqrt{1-k}+\mathrm{o}(1).\tag{4.19}$$

*Furthermore, the coefficient of the remainder term is locally uniformly bounded as $u$ varies in $\mathcal{A}_m$.*

**Remark 4.3.13.** In the case where $u$ solves the gaussian curvature equation (2.43), the intrinsic curvature of the immersion $\Phi[u]$ is constant and equal to $(k-1)$. In particular, horocircles in this surface have constant geodesic curvature equal to $-\sqrt{1-k}$ with respect to their inward-pointing normals. Since $\Phi^y[u]$ is the intersection of $\Phi[u]$ with the horizontal horosphere at height $e^y$, Lemma 4.3.4 confirms our expectation that the intersections of $k$-ends with horoballs are asymptotic to horodisks in the surface.

**Proof:** Indeed, by (2.44),

$$\kappa^y[u]=\frac{\mathcal{C}(u_y-u_{xx})}{(u+u_{xx})}=-\sqrt{1-k}+\mathrm{o}(1),$$

and the result follows. □

**4.4 - Killing vector fields.** The Möbius group $\mathrm{SO}(3,1)$ acts by orientation-preserving isometries on $\mathbb{H}^3$. Its Lie algebra $\mathfrak{so}(3,1)$ therefore defines a 6-dimensional family of vector fields over $\mathbb{H}^3$ whose flows preserve the metric. These vector fields are known as *Killing vector fields* of $\mathbb{H}^3$. In what follows, we will be particularly interested in the field

$$X_{a,b}:=\mathrm{M}_*\underline{a},$$

where $\mathrm{M}:\mathbb{H}^3\rightarrow\mathbb{H}^3$ is the orientation-reversing hyperbolic isometry given by

$$\mathrm{M}\underline{x}:=\frac{\underline{x}}{\|\underline{x}\|_e^2},\tag{4.20}$$

and

$$\underline{a}:=(a,b,0)^t$$

is a constant horizontal vector.

The field $X_{a,b}$ is given explicitly by

$$X_{a,b}(x)=\|\underline{x}\|_e^2\underline{a}-2\langle\underline{x},\underline{a}\rangle_e\underline{x}.\tag{4.21}$$

By (2.25), (2.33) and (2.34), for all $u\in\mathcal{A}_m$,

$$\begin{aligned}\langle\underline{a},\Phi[u]\rangle_g&=ae^{-y}u\mathrm{cos}(x)-ae^{-y}u_x\mathrm{sin}(x)+be^{-y}u_x\mathrm{cos}(x)+be^{-y}u\mathrm{sin}(x),\\ \langle\underline{a},\mathrm{N}[u]\rangle_g&=\mathcal{C}ae^{-y}\mathrm{cos}(x)+\mathcal{C}be^{-y}\mathrm{sin}(x)\ \text{and}\\ \langle\underline{a},\nu[u]\rangle_g&=\mathcal{S}ae^{-y}\mathrm{cos}(x)+\mathcal{S}be^{-y}\mathrm{sin}(x),\end{aligned}\tag{4.22}$$

so that, by (2.37) and (2.38),

$$\begin{aligned}\langle X_{a,b}(\Phi[u]),\mathrm{N}[u]\rangle_g &= a(1+f_1)e^y\cos(x)+af_2e^y\sin(x)\\ &\qquad +bf_3e^y\cos(x)+b(1+f_4)e^y\sin(x),\end{aligned}\tag{4.23}$$

where

$$f_1,f_2,f_3,f_4\in\mathcal{A}_m.$$

Likewise, by (2.34), (2.37), (2.38) and (4.17),

$$\begin{aligned}\langle X_{a,b}(\Phi[u]),\nu[u]\rangle_g &= -2a(u+g_1)e^y\cos(x)+2a(u_x+g_2)e^y\sin(x)\\ &\qquad -2b(u_x+g_3)e^t\cos(x)-2b(u+g_4)e^y\sin(x),\end{aligned}\tag{4.24}$$

where

$$g_1,g_2,g_3,g_4\in\mathcal{A}_{m,c}.$$

This yields

**Lemma 4.4.1**

*For every abstract $k$-end $u$, the limit*

$$\underset{T\to\infty}{Lim}\int_0^T\int_{mS^1}\langle X_{a,b}(\Phi[u]),N[u]\rangle_g dArea[u]$$

*converges as $T$ tends to infinity. Furthermore, this convergence is locally uniform as $u$ varies in $\mathcal{A}_m$.*

**Proof:** Indeed, since $\mathcal{A}_{m,c}$ is an ideal in $\mathcal{A}_m$, it follows by (2.30), (4.17) and (4.23) that

$$\langle X_{a,b}(\Phi[u]),\mathrm{N}[u]\rangle_g\mathrm{dArea}[u]=(f_1e^y\cos(x)+f_2e^y\sin(x))dxdy,$$

where $f_1,f_2\in\mathcal{A}_{m,c}$. The result now follows by Lemma 4.2.1. □

**Lemma 4.4.2**

*For every abstract $k$-end $u$,*

$$\int_{mS^1}H[u]\langle X_{a,b}(\Phi[u]),\nu[u]\rangle_g dl[u]=-4\pi m\langle(a,b),c[u]\rangle_e+O(e^{y-\sqrt{4-3k}y}).\tag{4.25}$$

*as $y$ tends to infinity. Furthermore, the coefficient of the remainder term is locally uniformly bounded as $u$ varies in $\mathcal{A}_m$.*

**Proof:** Since $\mathcal{A}_{m,c}$ is an ideal in $\mathcal{A}_m$, it follows by (2.31), (2.41) and (4.24) that

$$\begin{aligned}\mathrm{H}[u]\langle X_{a,b}(\Phi[u]),\nu[u]\rangle_g\mathrm{dl}[u] &= -2aue^y\cos(x)dx+2au_xe^y\sin(x)dx-2bu_xe^y\cos(x)dx\\ &\qquad -2bue^y\sin(x)dx+f_1e^y\cos(x)dx+f_2e^y\sin(x)dx,\end{aligned}$$

where $f_1,f_2\in\mathcal{A}_{m,c}$. Since

$$\begin{aligned}&\mathrm{c}_2[u_x]=-\mathrm{c}_1[u]\text{ and}\\ &\mathrm{c}_1[u_x]=\mathrm{c}_2[u],\end{aligned}$$

the result now follows by Lemma 4.2.1. □

Let $\partial_\nu$ denote the derivative in the direction of $\nu$.

**Lemma 4.4.3**

*For every abstract $k$-end $u$,*

$$\int_{mS^1}\partial_\nu\langle X_{a,b}(\Phi[u]),\nu[u]\rangle_g dl[u]=O(e^{y-\sqrt{4-3k}y}),\tag{4.26}$$

*as $y$ tends to infinity. Furthermore, the coefficient of the term on the right-hand side is locally uniformly bounded as $u$ varies in $\mathcal{A}_m$.*

**Proof:** By (2.31) and (2.36), for arbitrary $f:mS^1\times[0,\infty[\rightarrow\mathbb{R}$,

$$\partial_\nu f\mathrm{dl}[u]=\big(\mathcal{C}(u+u_{xx})f_x-\mathcal{C}(u_{xy}+u_x)f_y\big)dx.$$

Since $\mathcal{A}_{m,c}$ is an ideal which is closed under differentiation, it follows by (4.17) that

$$\partial_\nu f\mathrm{dl}[u]=(g_1f_x+g_2f_y)dx,$$

where $g_1,g_2\in\mathcal{A}_{m,c}$. It now follows by (4.24) that

$$\partial_\nu\langle X_{a,b}(\Phi[u]),\nu[u]\rangle_g\mathrm{dl}[u]=g_3e^y\mathrm{cos}(x)\mathrm{dx}+g_4e^y\mathrm{sin}(x)\mathrm{dx},$$

where $g_3,g_4\in\mathcal{A}_{m,c}$. The result now follows by Lemma 4.2.1. $\square$

Finally, consider the form $\alpha_\infty$ defined over $\mathbb{H}^3$ by

$$\alpha_\infty:=-\frac{1}{2z^2}dx\wedge dy.$$

The geometric significance of this form will become clear presently. For the moment, it will be sufficient to show

**Lemma 4.4.4**

*For every abstract $k$-end $u$,*

$$\int_{mS^1}i_{X_{a,b}}\alpha_\infty=O(e^{y-\sqrt{4-3k}y})\tag{4.27}$$

*as $y$ tends to infinity. Furthermore, the coefficient of the term on the right-hand side is locally uniformly bounded as $u$ varies in $\mathcal{A}_m$.*

**Proof:** Indeed, by (2.33), (2.34) and (5.14),

$$\begin{gathered}\alpha_\infty(\mathrm{N}[u],\mathrm{T}[u]):=1+f_1\ \text{and}\\ \alpha_\infty(\nu[u],\mathrm{T}[u]):=f_2,\end{gathered}$$

where $f_1\in\mathcal{A}_{m,*}$ and $f_2\in\mathcal{A}_m$. It follows by (4.23) and (4.24) that

$$\alpha_\infty(X_{a,b}(\Phi[u]),\mathrm{T}[u])=-(1+f_3)ae^y\mathrm{cos}(x)-(1+f_4)be^y\mathrm{sin}(x),$$

where $f_3,f_4\in\mathcal{A}_m$. Since $\mathcal{A}_{m,c}$ is an ideal in $\mathcal{A}_m$, it follows by (2.31) and (4.17) that

$$\alpha_\infty(X_{a,b}(\Phi[u]),\mathrm{T}[u])\mathrm{dl}[u]=f_5ae^y\mathrm{cos}(x)+f_6be^y\mathrm{sin}(x),$$

where $f_5,f_6\in\mathcal{A}_{m,c}$. The result now follows by Lemma 4.2.1. $\square$

## 5 - Area, generalized volume and renormalized energy.

**5.1 - Overview.** We now study the local geometry of the space $\mathcal{S}_k$ of $k$-surfaces. Recall from Section 1.4 that every stratum of $\mathcal{S}_k$ has a natural holomorphic structure given by the unordered vector of extremities. We first determine how variations of these extremities affect the surfaces. This is the content of Theorem 5.2.5: we show that smooth variations along the strata yield variations of the corresponding surfaces which are smooth in the space of immersions. In particular, given a point $(\underline{x},\underline{y})$ of the stratum, corresponding to the immersed surface $(S,e)$, every tangent vector $(\underline{a},\underline{b})$ of the stratum at this point corresponds to a vector field $\mathrm{X}[\underline{a},\underline{b}]$ over the immersion $e$ which is asymptotic over every end to some Killing field. The proof of Theorem 5.2.5 is a more or less standard application of the implicit function theorem for functions defined over Banach spaces. However, we draw the reader's attention to the small but important role played in the proof by the non-vanishing of the radii of the ends, proven in Lemma 4.3.3 (see Lemma 5.2.1, below). Finally, having established these preliminaries, our main results readily follow upon applying the asymptotic analysis developed in the preceding chapters. We first verify that area, generalized volume and renormalized energy are indeed well-defined and smooth over every stratum. Next, we verify an approximate Schläfli type formula over truncated surfaces. Finally, upon letting the locus of truncation tend to infinity, we obtain the Schläfli-type formula of Theorem 1.1.2 and the identity of Theorem 1.1.1 is obtained upon substituting the Killing fields of $\mathbb{H}^3$ into this formula.

**5.2 - Perturbations of finite-type $k$-surfaces.** Let $(S,e)$ be a finite-type $k$-surface in $\mathbb{H}^3$, let $n$ denote the number of ends of this surface and, for $1\leq i\leq n$, let $m_i$ denote the winding order of the $i$'th end. Choose an explicit upper half-space parametrisation of $\mathbb{H}^3$ as in Section 2.2 and suppose that none of the extremities $\mathrm{z}_1[e],\cdots,\mathrm{z}_n[e]$ of $(S,e)$ lie at infinity. For all $i$, let $h_i$ be a Busemann function of $\mathbb{H}^3$ centred at $\mathrm{z}_i[e]$. As in Section 1.4, for all $T\in\mathbb{R}$ and for all $i$, denote

$$\begin{aligned} S_T &:= \{x\in S\ |\ (h_i\circ e)(x)\geq T\ \forall i\}\,,\ \text{and}\\ S_{T,i} &:= \{x\in S\ |\ (h_i\circ e)(x)\leq T\}\,. \end{aligned} \tag{5.1}$$

Upon modifying $h_1,\cdots,h_n$ if necessary, we may assume that $\partial S_0$ is smooth and that the complement of its interior in $S$ consists of $n$ distinct $k$-ends, $S_{0,1},\cdots,S_{0,n}$. For each $i$, identify $S_{0,i}$ with $m_iS^1\times[0,\infty[$, let $\mathrm{M}_i$ be the hyperbolic isometry

$$\mathrm{M}_i\underline{x} := \frac{\underline{x}}{\|\underline{x}\|^2}+\mathrm{z}_i[e], \tag{5.2}$$

and let $u_i\in\mathcal{A}_{m_i}$ be an abstract $k$-end such that

$$e|_{S_{0,i}} = \mathrm{M}_i\circ\Phi[u_i], \tag{5.3}$$

where $\Phi$ is the operator defined in Section 2.4.

In order to use the asymptotic theory of the previous sections to determine the derivatives of the generalized volume and the renormalized energy, we will show in this section how nearby immersions in the stratum of $(S,e)$ in $\mathcal{S}_k$ depend smoothly on their extremities. To this end, we first construct an infinite-dimensional family of perturbations of $e$ in $C^\infty(S,\mathbb{H}^3)$. This construction is carried out in two stages. In the first, we construct a finite-dimensional family of perturbations which are large at infinity and, in the second, we extend this finite-dimensional family by an infinite-dimensional family of perturbations which are small at infinity. The finite-dimensional family is constructed as follows. Let $\tilde{\mathrm{e}}:\mathbb{R}^{2k}\rightarrow C^\infty(S,\mathbb{H}^3)$ be such that

(1) the function $\tilde{\mathrm{e}}[\underline{a},\underline{b}](p):\mathbb{R}^{2k}\times S\rightarrow\mathbb{H}^3$ is smooth,

(2) for all $p\in S$,

$$\tilde{\mathrm{e}}[0,0](p) = e(p),$$

(3) for all $(\underline{a},\underline{b})\in\mathbb{R}^{2k}$ and for all $p\in S_1$,

$$\tilde{\mathrm{e}}[\underline{a},\underline{b}](p) = e(p)\ \text{and}$$

(4) for all $1\leq i\leq n$, for all $(\underline{a},\underline{b})\in\mathbb{R}^{2k}$ and for all $p\in S_{0,i}$,

$$\tilde{\mathrm{e}}[\underline{a},\underline{b}](p) = \mathrm{T}[a_i,b_i]e(p),$$

where $\mathrm{T}[a_i,b_i]$ is the hyperbolic isometry defined in (2.8). For all $\delta>0$, let $B^{2k}_\delta(0)$ denote the ball of radius $\delta$ about the origin in $\mathbb{R}^{2k}$, and choose $\delta>0$ such that, for all $(\underline{a},\underline{b})\in B^{2k}_\delta(0)$, $\tilde{\mathrm{e}}[\underline{a},\underline{b}]$ is a complete, locally strictly convex immersion. The function $\tilde{\mathrm{e}}$ is the desired finite-dimensional family.

The infinite-dimensional extension of $\tilde{\mathrm{e}}$ is constructed as follows. Let $\chi:S\to[0,1]$ be a smooth function such that,

(1) for all $p\in S_0$,

$$\chi(p)=0,\ \text{ and}$$

(2) for all $1\leq i\leq n$ and for all $p\in S_{-1,i}$,

$$\chi(p)=1.$$

Define $\mathcal{G}:\mathbb{H}^3\times\mathbb{H}^3\times\mathbb{R}\to\mathbb{H}^3$ such that, for all $\underline{x},\underline{y}\in\mathbb{H}^3$ and for all $t\in\mathbb{R}$,

$$\mathcal{G}(\underline{x},\underline{y},t):=\gamma(t),$$

where $\gamma:\mathbb{R}\to\mathbb{H}^3$ is the unique geodesic such that $\gamma(0)=\underline{x}$ and $\gamma(1)=\underline{y}$. Let $\mathrm{N}:B^{2k}_\delta(0)\to C^\infty(S,\mathrm{U}\mathbb{H}^3)$ be such that, for all $(\underline{a},\underline{b})$ and for all $p\in S$, $\mathrm{N}[\underline{a},\underline{b}](p)$ is the outward-pointing unit normal vector of the immersion $\tilde{\mathrm{e}}[\underline{a},\underline{b}]$ at the point $p$. Define $\tilde{\mathrm{e}}:B^{2k}_\delta(0)\times C^0(S)\to C^0(S,\mathbb{H}^3)$ such that,

(1) for all $(\underline{a},\underline{b},v)\in B^{2k}_\delta(0)\times C^0(S)$ and for all $p\in S_0$,

$$\tilde{\mathrm{e}}[\underline{a},\underline{b},v](p):=\mathrm{Exp}(v(p)\mathrm{N}[\underline{a},\underline{b}](p))\ \text{and}\tag{5.4}$$

(2) for all $1\leq i\leq n$, for all $(\underline{a},\underline{b},v)\in B^{2k}_\delta(0)\times C^0(S)$ and for all $p\in S_{0,i}$,

$$\tilde{\mathrm{e}}[\underline{a},\underline{b},v](p):=\mathcal{G}(\mathrm{Exp}(v(p)\mathrm{N}[\underline{a},\underline{b}](p)),(\mathrm{M}_i\circ\Phi[u_i+v])(p),\chi(p)),\tag{5.5}$$

where, for each $i$, $\mathrm{M}_i$ and $u_i$ are the hyperbolic isometry and abstract $k$-end given by (5.2) and (5.3) respectively. This yields the desired infinite-dimensional extension.

Significantly, since $u_i$ decays rapidly over $S_{i,0}$ for all $i$, it is not clear that $\tilde{\mathrm{e}}[\underline{a},\underline{b},v]$ is an immersion even when $v$ itself has rapid decay. In particular, there is no neighbourhood of 0 in $B^{2k}_\delta(0)\times C^2(S)$ over which the extrinsic curvature operator can be meaningfully defined. For this reason, we introduce the operator

$$\mathrm{F}[\underline{a},\underline{b},v]:=\frac{\mathrm{K}[\underline{a},\underline{b},v]-k}{\mathrm{H}[\underline{a},\underline{b},v]},\tag{5.6}$$

where, for all $(\underline{a},\underline{b})\in\mathbb{R}^{2k}$ and for all $v\in C^2(S)$, $\mathrm{K}[\underline{a},\underline{b},v]$ and $\mathrm{H}[\underline{a},\underline{b},v]$ are respectively the extrinsic and mean curvature functions of the immersion $\tilde{\mathrm{e}}[\underline{a},\underline{b},v]$. It follows from (2.41) and (2.42) that this operator is well-defined for all sufficiently small $v$ even when $\tilde{\mathrm{e}}[\underline{a},\underline{b},v]$ is not an immersion.

For $\omega>0$ and for all $(k,\alpha)$, define the $C^{k,\alpha}_\omega$-norm of functions over $S$ by

$$\|u\|_{C^{k,\alpha}_\omega}=\|u|_{S_1}\|_{C^{k,\alpha}}+\sum_{i=1}^n\|u|_{S_{0,i}}\|_{C^{k,\alpha}_\omega},\tag{5.7}$$

and let $C^{k,\alpha}_\omega(S)$ denote the Banach space of $k$-times differentiable functions $u:S\to\mathbb{R}$ whose $C^{k,\alpha}_\omega$-norm is finite. Observe that, for all $(\underline{a},\underline{b})\in B^{2k}_\delta(0)$, $\mathrm{F}[\underline{a},\underline{b},0]$ is supported in $S_0$. It is then straightforward to show that, upon reducing $\delta$ if necessary, there exists neighbourhood $\mathcal{V}^{2,\alpha}_\omega(S)$ of zero in $C^{2,\alpha}_\omega(S)$ over which F defines a smooth function taking values in $C^{0,\alpha}_\omega(S)$.

**Lemma 5.2.1**

*After reducing $\delta$ and $\mathcal{V}^{2,\alpha}_\omega(S)$ if necessary, if $(\underline{a},\underline{b},v)\in B^{2k}_\delta(0)\times\mathcal{V}^{2,\alpha}_\omega(S)$ solves*

$$F[\underline{a},\underline{b},v]=0,$$

*then $\tilde{\mathrm{e}}[\underline{a},\underline{b},v]$ is a complete, locally strictly convex immersion of constant extrinsic curvature equal to $k$.*

**Proof:** It suffices to prove that these properties are satisfied over each end of $S$. However, for $1\leq i\leq n$, over $S_{-1,i}=m_iS^1\times[1,\infty[$, the function $u_i+v$ is an abstract $k$-end. In particular, it is an element of $\mathcal{A}_{m_i}$ so that, by (4.11),

$$(u_i+v)+(u_i+v)_{xx}=\mathrm{r}[u_i+v]e^{-\sqrt{1-k}y}+\mathrm{o}(e^{-\sqrt{1-k}y}),$$

where the coefficient of the remainder term is locally uniformly bounded as the $k$-end varies in $\mathcal{A}_m$. Since r is continuous, upon reducing $\delta$ and $\mathcal{V}^{2,\alpha}_\omega(S)$ if necessary, there exists $T<0$ such that, for every triple $(\underline{a},\underline{b},v)\in B^{2k}_\delta(0)\times\mathcal{V}^{2,\alpha}_\omega(S)$ which satisfies $\mathrm{F}[\underline{a},\underline{b},v]=0$, for all $1\leq i\leq n$, and for all $p\in S_{T,i}$,

$$(u_i+v)(p)+(u_i+v)_{xx}(p)>0,$$

so that, by (2.32), $\tilde{\mathrm{e}}[\underline{a},\underline{b},v]$ restricts to an immersion over $S_{T,i}$. Upon reducing $\delta$ and $\mathcal{V}^{2k}_\omega(S)$ further if necessary, we may then suppose that, for all such $(\underline{a},\underline{b},v)$, $\tilde{\mathrm{e}}[\underline{a},\underline{b},v]$ is an immersion over the whole of $S$. By a similar reasoning, for all such $(\underline{a},\underline{b},v)$, $\tilde{\mathrm{e}}[\underline{a},\underline{b},v]$ may also be taken to be locally strictly convex so that, by (2.41) and (2.42), $\mathrm{H}[\underline{a},\underline{b},v]$ and $K[\underline{a},\underline{b},v]$ are well-defined positive functions over $S$. In particular,

$$\mathrm{K}[\underline{a},\underline{b},v]=\mathrm{H}[\underline{a},\underline{b},v]\mathrm{F}[\underline{a},\underline{b},v]+k=k,$$

and the result follows. □

We now study the derivatives of F. Define the positive function $\phi:S\rightarrow]0,\infty[$ such that

(1) for all $p\in S_0$,

$$\phi(p)=1\text{ and}$$

(2) for all $1\leq i\leq n$, and for all $p\in S_{0,i}$,

$$\phi(p)=1+\chi(\mathcal{C}_i-1),$$

where

$$\mathcal{C}_i=\frac{1}{\sqrt{1+(u_i+u_{i,y})^2}}$$

is the function introduced in (2.5). Let $\mu_\phi$ denote the operator of multiplication by $\phi$. Since, for all $(k,\alpha)$, $\phi$ is an element of $C^{k,\alpha}(S)$, by Lemma 3.4.3, for all $\omega$, the operator $\mu_\phi$ defines a linear isomorphism from $C^{k,\alpha}_\omega(S)$ to itself.

**Lemma 5.2.2**

*The partial derivative of $F$ with respect to the third component at $(0,0,0)$ is given by*

$$D_3F[0,0,0]\cdot v=\frac{1}{H}J\mu_\phi v,$$

*where $H:=H[0,0,0]$ here denotes the mean curvature function of $e$ and $J$ denotes the Jacobi operator of extrinsic curvature for $e$.*

**Proof:** Indeed, choose $v\in C^2(S)$. Define the function $w:S\rightarrow\mathbb{R}$ by

$$w:=\left\langle\frac{\partial}{\partial t}\tilde{\mathrm{e}}[0,0,tv]\bigg|_{t=0},\mathrm{N}[e]\right\rangle_g.$$

By (5.4), over $S_0$,

$$w = v.$$

On the other hand, by (2.25) and (2.33), over each end

$$\left\langle \frac{\partial}{\partial t}\Phi[u+tv]\bigg|_{t=0}, \mathrm{N}[u]\right\rangle_g = \mathcal{C}v,$$

so that, by (5.5)

$$w = \mu_\phi v.$$

Since the extrinsic curvature of $(S,e)$ is constant, the first-order variation of K only depends on the normal component of the first-order variation of $e$, so that

$$\frac{\partial}{\partial t}\mathrm{K}[0,0,tv]\bigg|_{t=0} = \mathrm{J}w = \mathrm{J}\mu_\phi v.$$

The result now follows by the product rule, since the numerator in (5.6) vanishes at $(0,0,0)$. □

In Proposition 3.1.1 of [20] it is shown that the Jacobi operator of extrinsic curvature of $(S,e)$ is given by

$$\mathrm{J}v = \frac{1}{k}H(1-k)v - \mathrm{Tr}(A^{-1}\mathrm{Hess}(v)), \tag{5.8}$$

where $A := \mathrm{A}[e]$ is the shape operator of $e$ and $\mathrm{Hess}(v)$ is the hessian matrix of $v$ with respect to the metric $e^*g$.

**Lemma 5.2.3**

*The Jacobi operator of extrinsic curvature of $(S,e)$ satisfies*

$$\frac{k}{H}J\phi = (1-k)v - \hat{\Delta}v, \tag{5.9}$$

*where $\hat{\Delta}$ is the Laplace-Beltrami operator of the metric $I[e]+III[e]$.*

**Proof:** Indeed, the metric $\mathrm{I}[e]+\mathrm{III}[e]$ is given by

$$\hat{g} = (e^*g)((\mathrm{Id}+(1/k)A^2)\cdot,\cdot).$$

Since $e$ has constant extrinsic curvature equal to $k$, the Codazzi-Mainardi equations together with the Koszul formula yield

$$\hat{\Delta}v = \frac{k}{H}\mathrm{Tr}(A^{-1}\mathrm{Hess}(v)),$$

and the result follows. □

**Lemma 5.2.4**

*For $0<\omega<\sqrt{1-k}$, the operator $L := (1-k)-\hat{\Delta}$ defines a linear isomorphism from $C^{2,\alpha}_\omega(S)$ into $C^{0,\alpha}_\omega(S)$.*

**Proof:** The asymptotic properties of this operator over each end of $S$ are studied in Section 3.3. Together with the classical theory of elliptic operators (see [10]), these properties show that L defines a Fredholm map from $C^{2,\alpha}_\omega(S)$ into $C^{0,\alpha}_\omega(S)$. Since L is formally self-adjoint, it has Fredholm index equal to zero. Finally, by the maximum principle, L has trivial kernel in $C^{2,\alpha}_\omega(S)$, and invertibility follows. □

Lemmas 5.2.2 and 5.2.4 together with the implicit function theorem now yield

**Theorem 5.2.5**

*Upon reducing $\delta$ if necessary, there exists a smooth function $U : B^{2k}_\delta(0) \to C^{2,\alpha}_\omega(S)$ such that, for all $(\underline{a},\underline{b}) \in B^{2k}_\delta(0)$,*

$$K[\underline{a},\underline{b},U[\underline{a},\underline{b}]] = k.$$

*Furthermore, we may suppose that $U$ is unique.*

For all $(\underline{a},\underline{b}) \in B^{2k}_\delta(0)$, we define the immersion $\mathrm{e}[\underline{a},\underline{b}] : S \to \mathbb{H}^3$ by

$$\mathrm{e}[\underline{a},\underline{b}] := \tilde{\mathrm{e}}[\underline{a},\underline{b},\mathrm{U}[\underline{a},\underline{b}]]. \tag{5.10}$$

For all $(\underline{a},\underline{b}) \in \mathbb{R}^{2k}$, we define the vector field $\mathrm{X}[\underline{a},\underline{b}] : S \to \mathrm{T}\mathbb{H}^3$ by

$$\mathrm{X}[\underline{a},\underline{b}] := \frac{\partial}{\partial t}\mathrm{e}[t\underline{a},t\underline{b}]\bigg|_{t=0}, \tag{5.11}$$

and we define

$$\phi[\underline{a},\underline{b}] := \langle \mathrm{X}[\underline{a},\underline{b}],\mathrm{N}[e]\rangle_g. \tag{5.12}$$

The vector field $\mathrm{X}[\underline{a},\underline{b}]$ is the first-order variation of the immersion $e$ resulting from a first-order variation of the end point of the $i$'th end by the vector $(a_i,b_i)$. The function $\phi[\underline{a},\underline{b}]$ is the normal component of this first-order variation. By the preceding construction, for all $i$, over $S_{0,i}$,

$$\mathrm{X}[\underline{a},\underline{b}] = \phi_i\mathrm{N}[e] + \mathrm{M}_{i*}X_{a_i,b_i}\circ e, \tag{5.13}$$

where $\phi_i \in \mathcal{A}_{m_i}$, $\mathrm{M}_i$ is the hyperbolic isometry defined in (5.2) and $X_{a_i,b_i}$ is the Killing vector field defined in (4.21).

**5.3 - Area and generalized volume.** Let $(S,e)$ be a finite-type $k$-surface in $\mathbb{H}^3$. We continue to use the notation of Section 5.2. We first reprove the following elementary result of classical surface theory, which we believe provides a nice illustration of the estimates developed in the previous sections.

**Theorem 5.3.1**

*The area of $(S,e)$ is given by*

$$\mathit{Area}[e] = \frac{-2\pi\chi[S]}{(1-k)}.$$

*where $\chi[S]$ here denotes the Euler characteristic of $S$.*

**Proof:** Indeed, since $(S,e)$ has constant intrisic curvature equal to $(k-1)$, by the Gauss-Bonnet Theorem, for all $T$,

$$-(1-k)\mathrm{Area}[e|_{S_T}] + \int_{\partial S_T}\kappa[e]\mathrm{dl}[e] = 2\pi\chi[S],$$

where $\kappa$ here denotes the geodesic curvature of $\partial S_T$ with respect to the outward-pointing unit normal. However, by (4.18) and (4.19),

$$\underset{T\to\infty}{\mathrm{Lim}}\int_{\partial S_T}\kappa[e]\mathrm{dl}[e] = 0,$$

and the result follows. $\square$

The volume contained by $(S,e)$ is a slightly more subtle concept. Indeed, as $(S,e)$ is not necessarily embedded, it does not necessarily have a well-defined interior. However, since $\mathbb{H}^3$ is contractible, its volume form dVol is exact so that, by Stokes' Theorem, the volume contained within $(S,e)$ can be defined by integrating primitives of dVol over this surface. However, since $S$ itself is non-compact, there is no reason that two primitives should yield the same volume or even that an arbitrary primitive should be integrable over this surface. For this reason, we restrict attention to a special family of primitives. The *horospherical primitive* of dVol centred at infinity is defined by

$$\alpha_\infty := -\frac{1}{2z^2}dxdy. \tag{5.14}$$

Observe that, at each point, $\alpha_\infty$ is the pull-back under the orthogonal projection of $-(1/2)$ times the area form of the horizontal horosphere passing through that point. In particular, it is invariant under the action of those isometries of $\mathbb{H}^3$ which preserve the point at infinity. For any ideal point $w \in \partial\mathbb{H}^3$, the *horospherical primitive* of dVol centred at $w$ is now defined by

$$\alpha_w := \mathrm{M}^*\alpha_\infty,$$

where M is any isometry of $\mathbb{H}^3$ sending $w$ to $\infty$. By the preceding observation, for all $w$,

$$\|\alpha_w\|_g \leq \frac{1}{2}, \tag{5.15}$$

so that, since $(S, e)$ has finite area, the form $e^*\alpha_w$ has absolutely convergent integral over $S$. We now verify that this integral is independent of the horospherical primitive chosen. It suffices to compare $\alpha_\infty$ and $\alpha_0$.

**Lemma 5.3.2**

*The difference between $\alpha_\infty$ and $\alpha_0$ is given by*

$$\alpha_\infty - \alpha_0 = dlog(cosh(r))d\theta, \tag{5.16}$$

*where $r$ here denotes the distance in $\mathbb{H}^3$ from the geodesic $\Gamma_{0,\infty}$ and $\theta$ denotes the angle parameter of Fermi coordinates around this geodesic.*

**Remark 5.3.14.** In fact

$$d\log(\cosh(r))d\theta = \frac{1}{\cosh(r)}\beta, \tag{5.17}$$

where $\beta$ is the pull-back through the orthogonal projection of the area form of totally geodesic planes orthogonal to the geodesic $\Gamma_{0,\infty}$.

**Proof:** Indeed, let $\mathrm{M} : \mathbb{H}^3 \to \mathbb{H}^3$ be the isometry given by

$$\mathrm{M}\underline{x} = \frac{x}{\|\underline{x}\|^2}.$$

Since M reverses orientation

$$\alpha_0 = -\mathrm{M}^*\alpha_\infty = \frac{1}{2z^2}dxdy - \frac{1}{2z^2\rho^2}d(\rho^2)(xdy - ydx),$$

where

$$\rho^2 := x^2 + y^2 + z^2.$$

A straightforward calculation then yields

$$\alpha_\infty - \alpha_0 = \frac{1}{2}d\log\left(\frac{z^2}{\rho^2}\right) \wedge d\theta.$$

However, by elementary hyperbolic geometry

$$\cosh(r) = \frac{\rho}{z},$$

and the result follows. □

**Lemma 5.3.3**

*For all $z,w\in\partial_\infty\mathbb{H}^3$,*

$$\int_S e^*\alpha_z=\int_S e^*\alpha_w.$$

**Proof:** Indeed, there exists a constant $C>0$ such that

$$\|\mathrm{log}(\mathrm{cosh}(r))d\theta\|_g\leq C.$$

It follows that, for all $T$,

$$\left|\int_{S_T}e^*(\alpha_z-\alpha_w)\right|=\int_{\partial S_T}e^*\mathrm{log}(\mathrm{cosh}(r))d\theta\leq C\int_{\partial S_T}\mathrm{dl},$$

and the result now follows by Lemma 4.3.2. □

The *generalized volume* contained by $(S,e)$ is defined by

$$\mathrm{Vol}[e]:=\int_S e^*\alpha_z, \tag{5.18}$$

where $z$ is any ideal point of $\partial_\infty\mathbb{H}^3$ and $\alpha_z$ is the horospherical primitive of dVol centred at this point. Since the restriction of this integral to each end of $(S,e)$ varies smoothly with the end, it follows that Vol restricts to a smooth function over every stratum of $\mathcal{S}_k$. Finally, the reader may readily verify that when $(S,e)$ is embedded, $\mathrm{Vol}[e]$ coincides with the volume of the convex body in $\mathbb{H}^3$ that this embedding bounds.

**5.4 - Renormalized energy.** Let $(S,e)$ be a finite-type $k$-surface. As in Section 1.4, let $\hat{e}:S\rightarrow\mathrm{U}\mathbb{H}^3$ be the unit normal vector field over $S$, considered as an immersion in its own right in the total space of $\mathrm{U}\mathbb{H}^3$. The area form of the pull-back through this map of a suitable rescaling of the Sasaki metric is

$$\mathrm{d}\hat{\mathrm{E}}[e]:=\mathrm{H}[e]\mathrm{dArea}[e], \tag{5.19}$$

where $\mathrm{H}[e]:S\rightarrow\mathbb{R}$ denotes the mean curvature function of $e$. By Lemma 4.3.1, although the area of $S$ with respect to this form is infinite, the area of the truncated surface $S_T$ grows linearly with the absolute value of $T$ as $T$ tends to $-\infty$. The residue obtained upon subtracting this linear term yields a function over the space of $k$-surfaces which is well-defined up to a constant. More precisely, the *renormalized energy* of $(S,e)$ with respect to the Busemann functions $h_1,\cdots,h_n$ is defined by

$$\hat{\mathrm{E}}[e;h_1,\cdots,h_n]:=\underset{T\rightarrow-\infty}{\mathrm{Lim}}\int_{S_T}\mathrm{H}[e]\mathrm{dArea}[e]+2\pi T\sum_{i=1}^n m_i. \tag{5.20}$$

Trivially, if $h_1',\cdots,h_n'$ are other Busemann functions centred respectively on $\mathrm{z}_1[e],\cdots,\mathrm{z}_n[e]$, then

$$\hat{\mathrm{E}}[e;h_1',\cdots,h_n']-\hat{\mathrm{E}}[e;h_1,\cdots,h_n]=2\pi\sum_{i=1}^n m_i(h_i-h_i'),$$

which is constant since, for all $i$, $(h_i-h_i')$ is constant over $\mathbb{H}^3$.

The dependence of the renormalized energy on the Busemann functions is reduced as follows. Given another Busemann function $h_\infty$, the Busemann functions $h_1,\cdots,h_n$ are normalized such that, for all $i$, the horospheres $h_\infty^{-1}(\{0\})$ and $h_i^{-1}(\{0\})$ meet tangentially at a single point. Given that we are working in the upper half-space, it makes sense to choose

$$h_\infty(\underline{x})=-\mathrm{log}(z),$$

so that, for each $i$, the Busemann function $h_i$ is normalized by

$$h_i(\mathrm{z}_i[e],1)=0.$$

With these normalisations we set

$$\hat{\mathrm{E}}[e]:=\hat{\mathrm{E}}[e;h_1,\cdots,h_n], \tag{5.21}$$

so that $\hat{\mathrm{E}}[e]$ is uniquely defined a choice of Busemann function at $\infty$ or, equivalently, given an explicit upper half-space parametrisation of $\mathbb{H}^3$. As with the generalized volume, the renormalized energy defines a function over each end of $(S,e)$ which varies smoothly with the end, so that $\hat{\mathrm{E}}$ restricts to a smooth function over an open, dense subset of each stratum of $\mathcal{S}_k$.

**5.5 - The Schläfli formula.** Let $(S,e)$ be a finite-type $k$-surface. Let $X$ be the stratum of $\mathcal{S}_k$ in which it lies. Using the notation of Section 5.2, for all real $T$, define

$$\mathrm{Vol}_T[e] := \int_{S_T} e^*\alpha_\infty. \tag{5.22}$$

Using the local parametrisation of $X$ given by (5.10), we identify every tangent vector of this stratum at $(S,e)$ with a vector $(\underline{a},\underline{b})\in\mathbb{R}^{2k}$. Since we work with an explicit upper half-space parametrisation of $\mathbb{H}^3$, it will be useful to recall the definition of Steiner vectors (c.f. Section 1.6). For all $i$, let $\mathrm{z}_i[e]$ and $\zeta_i[e]$ denote respectively the $i$'th extremity and Steiner point of $(S,e)$, and define the $i$'th *Steiner vector* by

$$\mathrm{c}_i[e] := \frac{1}{\overline{\zeta}_i[e]-\overline{\mathrm{z}}_i[e]}. \tag{5.23}$$

**Lemma 5.5.1**

*For sufficiently large, negative $T$, the derivative of $\mathrm{Vol}_T$ at $(S,e)$ satisfies*

$$D\mathrm{Vol}_T[e]\cdot(\underline{a},\underline{b}) := \int_{S_T}\phi[\underline{a},\underline{b}]d\mathrm{Area}[e] + O(e^{\sqrt{1-k}T}). \tag{5.24}$$

**Proof:** Indeed, denoting $X := X[\underline{a},\underline{b}]$, we have

$$\begin{aligned}
D\mathrm{Vol}_T[e]\cdot(\underline{a},\underline{b}) &= \int_{S_T}\mathcal{L}_X\alpha_\infty\\
&= \int_{S_T} i_Xd\alpha_\infty + di_X\alpha_\infty\\
&= \int_{S_T}\langle X,\mathrm{N}[u]\rangle \mathrm{dArea}[e] + \int_{\partial S_T} i_X\alpha.
\end{aligned}$$

However, by (5.13), for $1\leq i\leq n$, over $S_{0,i}$,

$$X = \phi_i\mathrm{N}[e] + \mathrm{M}_{i*}X_i\circ e,$$

where $\phi_i\in\mathcal{A}_{m_i}$ and $X_i := X_{a_i,b_i}$ is the Killing vector field defined in (4.21). Since $\|\alpha\|_g = 1/2$,

$$\|i_{\phi_i\mathrm{N}[e]}\alpha\| \leq \frac{1}{2}\left|\phi_i\right|.$$

so that, since $\phi_i\in\mathcal{A}_{m_i}$,

$$\int_{\partial S_T} i_{\phi_i\mathrm{N}[e]}\alpha_\infty = \mathrm{O}(e^{\sqrt{1-k}T}).$$

On the other hand, by Lemma 4.4.4,

$$\int_{\partial S_T}(M_i\circ e)^* i_{\mathrm{M}_{i*}X_i}\alpha_\infty = \int_{\partial S_T} e^* i_{X_i}\alpha_{\mathrm{M}_i^{-1}\infty} = \mathrm{O}(e^{\sqrt{1-k}T}),$$

and the result follows. □

**Theorem 5.5.2**

*The derivative of Vol at $(S,e)$ is given by*

$$D\mathit{Vol}[e]\cdot(\underline{a},\underline{b}) = \underset{T\rightarrow-\infty}{\mathit{Lim}}\int_{S_T}\phi[\underline{a},\underline{b}]d\mathrm{Area}[e]. \tag{5.25}$$

**Proof:** By (5.13) and Lemma 4.4.1, the limit on the right-hand side of (5.25) exists and converges locally uniformly as $(S,e)$ varies along its stratum in $\mathcal{S}_k$. The result follows. □

For all $T$, define

$$\hat{\mathrm{E}}_T[e] := \int_{S_T}\mathrm{H}[e]\mathrm{dArea}[e]. \tag{5.26}$$

**Lemma 5.5.3**

*For sufficiently large, negative $T$, the derivative of $\hat{\mathrm{E}}_T$ at $(S,e)$ is given by*

$$\begin{aligned} D\hat{\mathrm{E}}_T[e]\cdot(\underline{a},\underline{b}) = &\int_{S_T} 2(1+k)\phi[\underline{a},\underline{b}]dArea[e] - \int_{\partial S_T}\partial_\nu\phi[\underline{a},\underline{b}]dl[e]\\ &-\sum_{i=1}^n\int_{\partial S_{T,i}}\mathcal{T}\phi_i H[e]dl[e]\\ &+\sum_{i=1}^n\int_{\partial S_{T,i}}\langle X_{a_i,b_i}(\Phi[e]),\nu[e]\rangle H[e]dl[e], \end{aligned}\tag{5.27}$$

*where, for each $1\leq i\leq n$, $\phi_i$ and $X_{a_i,b_i}$ are as in* (5.13).

**Proof:** Choose $(\underline{a},\underline{b})\in\mathbb{R}^{2k}$. For all sufficiently small $t$, let $e_t := e[t\underline{a},t\underline{b}]$ be as in (5.10). For $f:\mathbb{H}^3\times]-\epsilon,\epsilon[\rightarrow\mathbb{R}$ smooth, consider the function

$$\hat{\mathrm{E}}_f(t) := \int_S (f_t\circ e_t)\mathrm{H}[e_t]\mathrm{dArea}[e_t].$$

Since the integrand of $\hat{\mathrm{E}}_f$ is a smoothly varying family supported in a compact subset of $S$, the tangential component of $\mathrm{X}[\underline{a},\underline{b}]$ does not contribute to its derivative. We therefore work as if it were equal to zero, so that

$$\mathrm{X}[\underline{a},\underline{b}] = \phi[\underline{a},\underline{b}]\mathrm{N}[e].$$

The first-order variation of the area form is given by (c.f. [8]),

$$\frac{\partial}{\partial t}\mathrm{dArea}[e_t]\bigg|_{t=0} = \mathrm{H}[e]\phi[\underline{a},\underline{b}]\mathrm{dArea}[e].$$

The first-order variation of the mean curvature is given by (c.f. [8]),

$$\frac{\partial}{\partial t}\mathrm{H}[e_t]\bigg|_{t=0} = (2(1+k)-\mathrm{H}[e]^2)\phi[\underline{a},\underline{b}] - \Delta\phi[\underline{a},\underline{b}].$$

Finally, the first-order variation of $f_t\circ e_t$ is given by

$$\frac{\partial}{\partial t}(f_t\circ e_t)\bigg|_{t=0} = \frac{\partial f}{\partial t}\circ e + \langle\nabla f_0,\mathrm{N}[e]\rangle_g\phi[\underline{a},\underline{b}].$$

It follows by the product rule that

$$\begin{aligned} \frac{\partial}{\partial t}\hat{\mathrm{E}}_f\bigg|_{t=0} = &\ 2(1+k)\int_S(f_t\circ e)\phi[\underline{a},\underline{b}]\mathrm{dArea}[e] - \int_S(f_t\circ e)\Delta\phi[\underline{a},\underline{b}]\mathrm{dArea}[e]\\ &+\int_S\left(\frac{\partial f}{\partial t}\bigg|_{t=0}\circ e\right)\mathrm{H}[e]\mathrm{dArea}[e]\\ &+\int_S\langle\nabla^g f_0,\mathrm{N}[e]\rangle_g\phi[\underline{a},\underline{b}]\mathrm{H}[e]\mathrm{dArea}[e]. \end{aligned}\tag{5.28}$$

Now, for all $1\leq i\leq n$, let $z_{i,t}$ be the $i$'th extremity of $(S,e_t)$ and let $h_{i,t}:\mathbb{H}^3\rightarrow\mathbb{R}$ be a Busemann function centred on $z_{i,t}$ and normalized such that, for all $t$,

$$h_{i,t}(z_{i,t},1) = 0.$$

By definition, for all $1\leq i\leq n$,

$$\frac{\partial}{\partial t}z_{i,t}\bigg|_{t=0}=(a_i,b_i).$$

Define $f:\mathbb{H}^3\times]-\epsilon,\epsilon[\rightarrow\mathbb{R}$ by

$$f_t(x)=\begin{cases}1 & \text{if } h_{i,t}(x)\geq T\ \forall i \text{ and}\\ 0 & \text{otherwise.}\end{cases}$$

Observe that $f$ is an element of $\mathrm{BV}_{\mathrm{loc}}(\mathbb{H}^3\times]-\epsilon,\epsilon[)$, the space of functions of locally bounded total variation over $\mathbb{H}^3\times]-\epsilon,\epsilon[$ (c.f. [28]). By approximating $f$ by smooth functions, we find that (5.28) continues to hold with the derivatives of $f$ now being interpreted in the distributional sense. In particular,

$$\int_S\bigg(\frac{\partial f}{\partial t}\bigg|_{t=0}\circ e\bigg)\mathrm{H}[e]\mathrm{dArea}[e]=\sum_{i=1}^n\int_{\partial S_{T,i}}\mathcal{C}^{-1}\langle X_{a_i,b_i}(\Phi[e]),\nabla^g h_{i,0}\circ e\rangle_g\mathrm{H}[e]\mathrm{dl}[e],$$

and

$$\int_S\langle\nabla^g f_0,\mathrm{N}[e]\rangle_g\phi[\underline{a},\underline{b}]\mathrm{H}[e]\mathrm{dArea}[e]=-\sum_{i=1}^n\int_{\partial S_{T,i}}\mathcal{C}^{-1}\langle\mathrm{N}[e],\nabla^g h_{i,0}\circ e\rangle_g\phi[\underline{a},\underline{b}]\mathrm{H}[e]\mathrm{dl}[e],$$

where, for each $i$, $\nabla^g h_{i,0}$ denotes the gradient with respect to $g$ of the function $h_{i,0}$. Since

$$\nabla^g h_{i,0}\circ e=-\mathcal{C}\nu[e]+\mathcal{S}\mathrm{N}[e],$$

the result now follows upon substituting these relations into (5.28) and applying Stokes' Theorem. □

**Lemma 5.5.4**

*For all sufficiently large, negative $T$, the derivatives of $Vol_T$ and $\hat{\mathrm{E}}_T$ at $(S,e)$ are related by*

$$2(1+k)D\,Vol_T[e]\cdot(\underline{a},\underline{b})-D\hat{\mathrm{E}}_T[e]\cdot(\underline{a},\underline{b})=\sum_{i=1}^n 4\pi m_i\langle(a_i,b_i),\mathrm{c}_i[e]\rangle_e+O(e^{\sqrt{1-k}T}),\tag{5.29}$$

*where, for each $i$, $m_i$ denotes the winding order of the $i$'th end of $(S,e)$ and $\mathrm{c}_i[e]$ denotes its Steiner vector given by* (5.23).

**Proof:** It suffices to analyse the last three terms of (5.27) over each end. Choose $1\leq i\leq k$. By (5.13), over $S_{T,i}$,

$$X[\underline{a},\underline{b}]=\phi\mathrm{N}[e]+\mathrm{M}_{i*}X_{a_i,b_i}\circ e,$$

where $\phi\in\mathcal{A}_{m_i}$ and $X_{a_i,b_i}$ is the vector field given in (4.21). Since $\phi\in\mathcal{A}_{m_i}$, by (2.5), (2.31), (2.36) and (2.41),

$$\int_{\partial S_{T,i}}\partial_\nu\phi_i\mathrm{dl}[e]=\mathrm{O}(e^{2\sqrt{1-k}T})\text{ and}$$

$$\int_{\partial S_{T,i}}\mathrm{H}[e]\mathcal{S}\phi_i\mathrm{dl}[e]=\mathrm{O}(e^{2\sqrt{1-k}T}),$$

by Lemma 4.4.3,

$$\int_{\partial S_{T,i}}\partial_\nu\langle X_{a_i,b_i},\mathrm{N}[e]\rangle\mathrm{dl}[e]=\mathrm{O}(e^{-T+\sqrt{4-3k}T}),$$

and, by Lemma 4.4.2,

$$\int_{\partial S_{T,i}}H[e]\mathcal{C}\langle X,\nu[e]\rangle\mathrm{dl}[e]=-4\pi m_i\langle(a_i,b_i),\mathrm{c}_i[e]\rangle+\mathrm{O}(e^{-T+\sqrt{4-3k}T}).$$

The result follows. □

**Theorem 5.5.5**

*The derivatives of Vol and* $\hat{\mathrm{E}}$ *at* $(S,e)$ *are related by*

$$2(1+k)D\mathit{Vol}[e]\cdot(\underline{a},\underline{b})-D\hat{\mathrm{E}}[e]\cdot(\underline{a},\underline{b})=\sum_{i=1}^{n}4\pi m_i\langle(a_i,b_i),c_i[e]\rangle_e,\tag{5.30}$$

*where, for each* $i$, $m_i$ *denotes the winding order of the* $i$*'th end of* $(S,e)$ *and* $c_i[e]$ *denotes its Steiner vector given by* (5.23).

**Remark 5.5.15.** Theorem 1.1.2 follows immediately upon expressing (5.30) in Möbius-invariant form.

**Proof:** This follows immediately from (5.29) since the limit converges locally uniformly as $(S,e)$ varies along its stratum in $\mathcal{S}_k$. □

**Theorem 5.5.6**

*The extremities and Steiner vectors of* $(S,e)$ *satisfy*

$$\sum_{i=1}^{n}m_ic_i[e]=0,\tag{5.31}$$

$$\sum_{i=1}^{n}m_ic_i[e]\overline{\mathrm{z}}_i[e]=-\frac{1}{2}\sum_{i=1}^{n}m_i,\ \text{and}\tag{5.32}$$

$$\sum_{i=1}^{n}m_ic_i[e]\overline{\mathrm{z}}_i[e]^2=-\sum_{i=1}^{n}m_i\overline{\mathrm{z}}_i[e],\tag{5.33}$$

*where, for each* $i$, $m_i$ *denotes the winding order of the* $i$*'th extremity of* $(S,e)$.

**Remark 5.5.16.** Theorem 1.1.1 follows immediately upon expressing (5.31), (5.32) and (5.33) in Möbius invariant form.

**Proof:** These relations are obtained by applying Killing vector fields to the Schläfli formula (5.30). It suffices to prove the real part of (5.32), as the proofs of the remaining formulae are identical. Consider the family of hyperbolic isometries given by

$$\mathrm{M}_t\underline{x}=e^t\underline{x}.$$

The Killing vector field of this family is

$$X(\underline{x})=\underline{x}.$$

Consequently,

$$\left.\frac{\partial}{\partial t}\mathrm{z}_i[M_t\circ e]\right|_{t=0}=\mathrm{z}_i[e],$$

$$\left.\frac{\partial}{\partial t}\mathrm{Vol}[M_t\circ e]\right|_{t=0}=0\ \text{and}$$

$$\left.\frac{\partial}{\partial t}\hat{\mathrm{E}}[M_t\circ e]\right|_{t=0}=2\pi\sum_{i=1}^{n}m_i.$$

Substituting these values into (5.30) yields the real part of (5.32). The imaginary part is obtained in the same manner using rotations. Translations yield (5.31) and parabolic transformations about the origin yield (5.33). This completes the proof. □

## A - Smooth functions over Hölder spaces.

For the reader's convenience, we review the smoothness properties of composition operators over Hölder spaces. Although similar properties are studied in [12] and [34], it is not clear to us where a straightforward treatment of the difficulties particular to the Hölder space case may be found in the literature. Let $X$ be a metric space and let $E$ be a Banach space. For $\alpha\in[0,1]$, the $\alpha$-Hölder seminorm of a function $f:X\rightarrow E$ is defined by

$$[f]_\alpha:=\mathrm{Sup}_{x\neq y}\frac{\|f(x)-f(y)\|}{d(x,y)^\alpha}. \tag{A.1}$$

Observe that the 1-Hölder seminorm is the Lipschitz seminorm whilst the 0-Hölder seminorm is the total oscilation. The $C^{0,\alpha}$-norm is then defined by

$$\|f\|_{C^{0,\alpha}}:=\|f\|_{C^0}+[f]_\alpha. \tag{A.2}$$

More generally, when $X$ is a riemannian manifold, which for convenience we take to be locally isometric to $\mathbb{R}^m$ for some $m$, for all $(k,\alpha)$, the $C^{k,\alpha}$-norm of a $k$-times differentiable function $f:X\rightarrow E$ is defined by

$$\|f\|_{C^{k,\alpha}}:=\sum_{i=0}^k\|D^if\|_{C^0}+[D^kf]_\alpha. \tag{A.3}$$

Observe, in particular, that for all $k\geq 1$,

$$\|f\|_{C^{k,\alpha}}=\|f\|_{C^0}+\|Df\|_{C^{k-1,\alpha}}. \tag{A.4}$$

This recurrence relation will be useful for the induction arguments that we will invoke presently. For all $(k,\alpha)$, let $C^{k,\alpha}(X,E)$ denote the Banach space of functions with finite $C^{k,\alpha}$ norm. We readily obtain

**Lemma A.1.1**

*For $f,g\in C^{0,\alpha}(X,E)$,*

$$\|fg\|_{C^{0,\alpha}}\leq\|f\|_{C^{0,\alpha}}\|g\|_{C^{0,\alpha}}. \tag{A.5}$$

An induction argument, starting with Lemma A.1.1 and using (A.4), then yields

**Lemma A.1.2**

*For all $f,g\in C^{k,\alpha}(X,E)$,*

$$\|fg\|_{C^{k,\alpha}}\leq(2^{k+1}-1)\|f\|_{C^{k,\alpha}}\|g\|_{C^{k,\alpha}}. \tag{A.6}$$

Now let $\Omega$ be an open subset of $E$. Define

$$\mathcal{O}^{k,\alpha}(X,\Omega):=\underset{\epsilon>0}{\cup}\mathcal{O}^{k,\alpha}_\epsilon(X,\Omega), \tag{A.7}$$

where, for all $\epsilon>0$,

$$\mathcal{O}^{k,\alpha}_\epsilon(X,\Omega):=\left\{f\in C^{k,\alpha}(X,E)\ |\ d(f(x),\Omega^c)\geq\epsilon\ \forall x\right\}, \tag{A.8}$$

Observe that $\mathcal{O}^{k,\alpha}(X,\Omega)$ is an open subset of $C^{k,\alpha}(X,E)$. Given another Banach space $F$ and a suitably regular function $\Phi:\Omega\rightarrow F$, define the composition operator $\mathrm{C}_\Phi:\mathcal{O}^{k,\alpha}(X,\Omega)\rightarrow C^{k,\alpha}(X,F)$ by

$$\mathrm{C}_\Phi[f]:=\Phi\circ f. \tag{A.9}$$

Smoothness of composition operators over Hölder spaces is subtle for low regularity. First, we have

**Lemma A.1.3**

*If $\Phi\in C^{1,1}(\Omega,F)$, then $C_\Phi$ defines a continuous function from $\mathcal{O}^{0,\alpha}(X,\Omega)$ into $C^{0,\alpha}(X,F)$.*

**Proof:** Without loss of generality, we may suppose that $\Omega$ is convex. Then, for all $f, f+g\in\mathcal{O}^{0,\alpha}(X,\Omega)$, and for all $x\in X$,

$$\mathrm{C}_\Phi[f+g](x)-\mathrm{C}_\Phi[f](x)=\int_0^1\mathrm{C}_{D\Phi}[f+sg](x)g(x)dx.$$

Thus, for all $x,y\in X$,

$$\begin{aligned}&|\mathrm{C}_\Phi[f+g](x)-\mathrm{C}_\Phi[f](x)-\mathrm{C}_\Phi[f+g](y)+\mathrm{C}_\Phi[f](y)|\\ &\qquad=\left|\int_0^1\mathrm{C}_{D\Phi}[f+sg](x)g(x)-\mathrm{C}_{D\Phi}[f+sg](y)g(y)ds\right|\\ &\qquad\leq\int_0^1[\mathrm{C}_{D\Phi}[f+sg]g]_\alpha dsd(x,y)^\alpha\\ &\qquad\leq\int_0^1\big([D\Phi]_1[f+sg]_\alpha\|g\|_{C^0}+\|D\Phi\|_{C^0}[g]_\alpha\big)dsd(x,y)^\alpha,\end{aligned}$$

so that

$$[\mathrm{C}_\Phi[f+g]-\mathrm{C}_\Phi[f]]_\alpha\leq\|\Phi\|_{C^{1,1}}\big([f]_\alpha+[g]_\alpha\big)\|g\|_{C^{0,\alpha}}.$$

Since

$$\|\mathrm{C}_\Phi[f+g]-\mathrm{C}_\Phi[f]\|_{C^0}\leq\|\Phi\|_{C^{0,1}}\|g\|_{C^0},$$

continuity now follows. □

**Lemma A.1.4**

*If $\Phi\in C^{2,1}(\Omega,F)$, then $C_\Phi$ defines a continuously differentiable function from $\mathcal{O}^{0,\alpha}(X,\Omega)$ into $C^{0,\alpha}(X,F)$ with derivative given by*

$$(DC_\Phi[f]g)(x)=C_{D\Phi}[f](x)g(x).\tag{A.10}$$

**Proof:** Without loss of generality, we may suppose again that $\Omega$ is convex. Define $\Psi:\mathcal{O}^{0,\alpha}(X,\Omega)^2\rightarrow C^{0,\alpha}(X,F)$ by

$$\Psi[f,g]:=\mathrm{C}_{D\Phi}[g]-\mathrm{C}_{D\Phi}[f].$$

By Lemma A.1.3, $\Psi$ is continuous. By the fundamental theorem of calculus, for all $f,f+g\in\mathcal{O}^{0,\alpha}(X,\Omega)$ and for all $x\in X$,

$$\begin{aligned}\mathrm{C}_\Phi[f+g](x)-\mathrm{C}_\Phi[f](x)&=\int_0^1\mathrm{C}_{D\Phi}[f+sg](x)dsg(x)\\ &=\int_0^1\Psi[f,f+sg](x)dsg(x)+\mathrm{C}_{D\Phi}[f](x)g(x).\end{aligned}$$

By continuity, the curve $s\mapsto\Psi[f,f+sg]$ is integrable as a curve taking values in the Banach space $C^{0,\alpha}(X,F)$. Furthermore, by convexity of the norm,

$$\left\|\int_0^1\Psi[f,f+sg]ds\right\|_{C^{0,\alpha}}\leq\int_0^1\|\Psi[f,f+sg]\|_{C^{0,\alpha}}ds\leq\operatorname*{Sup}_{s\in[0,1]}\|\Psi[f,f+sg]\|_{C^{0,\alpha}}.$$

Since $\Psi$ vanishes when $g=f$, this tends to 0 as $g$ tends to $f$, and the result now follows by Lemma A.1.1. □

**Lemma A.1.5**

*If $\Phi \in C^{k+1,1}(\Omega, F)$, then $C_\Phi$ defines a continuous function from $\mathcal{O}^{k,\alpha}(X,\Omega)$ into the space $C^{k,\alpha}(X,F)$.*

**Proof:** We prove this by induction on $k$. By Lemma A.1.3, the result holds when $k=0$. Moreover, by the chain rule, for all $f \in \mathcal{O}^{k+1,\alpha}(X,\Omega)$,

$$D(\mathrm{C}_\Phi[f]) = \mathrm{C}_{D\Phi}[f]Df.$$

It follows by Lemma A.1.2 and the induction hypothesis that the operator $f \mapsto D\mathrm{C}_\Phi[f]$ defines a continuous function from $\mathcal{O}^{k+1,\alpha}(X,E)$ into $C^{k,\alpha}(X,\mathrm{Lin}(\mathbb{R}^n,F))$. Since $\mathrm{C}_\Phi$ trivially defines a continuous function from $\mathcal{O}^{k+1,\alpha}(X,E)$ into $C^0(X,F)$, the result now follows by (A.4). $\square$

In the same manner, we obtain

**Lemma A.1.6**

*If $\Phi \in C^{k+2,1}(\Omega, F)$, then $C_\Phi$ defines a continuously differentiable function from $\mathcal{O}^{k,\alpha}(X,\Omega)$ into $C^{k,\alpha}(X,F)$ with derivative given by* (A.10).

Finally, by induction on $l$, we obtain

**Lemma A.1.7**

*If $\Phi \in C^{k+l+1,1}(\Omega, F)$, then $C_\Phi$ defines a $C^l$ function from $\mathcal{O}^{k,\alpha}(X,\Omega)$ into $C^{k,\alpha}(X,F)$.*

We leave the reader to verify that when $k \geq 1$, a more straightforward argument shows that $\mathrm{C}_\Phi$ is of class $C^l$ whenever $\Phi$ is of class $C^{k+l,\beta}$ for some $\beta > \alpha$. Furthermore, this condition is sharp in the sense that there exist functions $\Phi$, of class $C^{k+l,\alpha}$, for which $\mathrm{C}_\Phi$ is not $l$-times continuously differentiable. However, when $k=0$, we do not know whether the hypotheses of Lemma A.1.3 may be relaxed or whether the technical arguments of this appendix may be bypassed. This merely formal concern is moot for geometric applications where we rarely, if ever, are concerned with functions which are not smooth. That said, when studying curvature problems, it is generally more elegant to work with functions of class $C^{2,\alpha}$, being the minimal Hölder regularity required to apply analytic techniques. The arguments presented in this appendix then become necessary when the curvature problems in question involve totally non-linear, second-order, partial differential operators, since the second derivatives of such functions are of class $C^{0,\alpha}$ and compositions of these derivatives by smooth functions could not otherwise be assumed to define smooth operations over function spaces. There is not a single situation we know of where this subtlety cannot be easily bypassed by working with functions of greater regularity. Nevertheless, it is clearly important to be formally correct.

## B - Bibliography.

[1] Alekseevskij D. V., Vinberg E. B., Solodovnikov A. S., *Geometry II: Spaces of constant curvature*, Encyclopaedia Math. Sci., **29**, Springer, (1993)

[2] Ballmann W., Schroeder V., Gromov M., *Manifolds of Nonpositive Curvature*, Progress in Mathematics, **61**, Birkhäuser, (1985)

[3] Barbot T., Béguin F., Zeghib A., Prescribing Gauss curvature of surfaces in 3-dimensional spacetimes: application to the Minkowski problem in the Minkowski space, *Ann. Inst. Fourier*, **61**, no. 2, (2011), 511–591

[4] Berger M., Gauduchon P., Mazet E., *Le Spectre dune Variété Riemannienne*, Lecture Notes in Mathematics, **194**, Springer, (1971)

[5] Bonahon F., A Schläfli-type formula for convex cores of hyperbolic 3-manifolds, *J. Diff. Geom.*, **50**, (1998), 25–58

[6] Bonsante F., Mondello G., Schlenker J. M., A cyclic extension of the earthquake flow I, *Geom. Topol.*, **17**, no. 1, (2013), 157–234

[7] Bonsante F., Mondello G., Schlenker J. M., A cyclic extension of the earthquake flow II, *Ann. Sci. Ec. Norm. Supér.*, **48**, no. 4, (2015), 811–859

[8] Colding T. H., Minicozzi W. P. II, *A course in minimal surfaces*, Graduate Studies in Mathematics, Amer. Math. Soc., (2011)

[9] Fillastre F., Smith G., Group actions and scattering problems in Teichmüller theory, in *Handbook of group actions IV*, Advanced Lectures in Mathematics, **40**, (2018), 359–417

[10] Gilbarg D., Trudinger N. S., *Elliptic partial differential equations of second order*, Classics in Mathematics, Springer-Verlag, Berlin, (2001)

[11] Grünbaum B., *Convex Polytopes*, Graduate Texts in Mathematics, **221**, Springer-Verlag, New York, (2003)

[12] Hamilton R. S., The inverse function theorem of Nash and Moser, *Bull. Amer. Math. Soc.*, **7**, no. 1, (1982), 65–222

[13] Honsberger R., *Episodes in nineteenth and twentieth century euclidean geometry*, New Mathematical Library, **37**, The Mathematical Association of America, (1995)

[14] Kimberling C., *Encyclopedia of Triangle Centers*, http://faculty.evansville.edu/ck6/encyclopedia

[15] Krasnov K., Schlenker J. M., A symplectic map between hyperbolic and complex Teichmller theory, *Duke Math. J.*, **150**, no. 2, (2009), 331–356

[16] Kubota T., Über die Schwerpunkte der konvexen geschlossenen Kurven und Flächen, *Tohoku Math. J.*, **14**, (1918), 20–27

[17] Labourie F., Exemples de courbes pseudo-holomorphes en géométrie riemannienne, in *Holomorphic curves in symplectic geometry*, Progress in Mathematics, **117**, Birkhäuser, (1994)

[18] Labourie F., Problèmes de Monge-Ampère, courbes pseudo-holomorphes et laminations, *Geom. Funct. Anal.*, **7**, (1997), 496–534

[19] Labourie F., Métriques prescrites sur le bord des variétés hyperboliques de dimension 3, *J. Diff. Geom.*, **35**, no. 3, (1992), 609–626

[20] Labourie F., Un Lemme de Morse pour les surfaces convexes, *Inventiones Mathematicae*, **141**, No 2, (2000), 239–297

[21] McDuff D., Salamon D., *Introduction to Symplectic Topology*, Oxford Mathematical Monographs, Oxford University Press, (1999)

[22] McDuff D., Salamon D., *J-holomorphic Curves and Symplectic Topology*, Colloquium Publications, **52**, Amer. Math. Soc., (2012)

[23] Rogers C., Schief W. K., *Bäcklund and Darboux Transformations: Geometry and Modern Applications in Soliton Theory*, Cambridge Texts in Applied Mathematics, **30**, Cambridge University Press, (2002)

[24] Schlenker J. M., Souam R., Higher Schläfli formulas and applications, *Comp. Mat.*, **135**, (2003), 1–24

[25] Schlenker J. M., Hyperbolic manifolds with convex boundary, *Inventiones math.*, **163**, (2006), 109–169

[26] Schlenker J. M., Souam R., Higher Schläfli formulas and applications II: vector-valued differential relations, *Int. Mat. Res. Not.*, (2008)

[27] Schneider R., On Steiner points of convex bodies, *Israel J. Math.*, **9**, no. 2, (1971), 241–249

[28] Simon L., *Lectures on geometric measure theory*, Proceedings of the Centre for Mathematical Analysis, **3**, Australian National University, (1983)

[29] Schneider R., *Convex bodies: the Brunn-Minkowski theory*, Encylopedia of Mathematics and its Applications, **44**, Cambridge University Press, Cambridge, (1993)

[30] Smith G., Pointed $k$-surfaces, *Bull. Soc. Math. France*, **134**, no. 4., (2006), 509557

[31] Smith G., Constant scalar curvature hypersurfaces in (3+1)-dimensional GHMC Minkowski spacetimes, *J. Geom. Phys.*, **128**

[32] Smith G., On an Enneper-Weierstrass-type representation of constant Gaussian curvature surfaces in 3-dimensional hyperbolic space, to appear in *Minimal surfaces: integrable systems and visualisation*

[33] Tamburelli A., Prescribing metrics on the boundary of anti-de Sitter 3-manifolds, *Int. Math. Res. Not.*, no. 5, (2018), 1281–1313

[34] Tayler M. E., *Partial differential equations III: non-linear equations*, Applied Mathematics Sciences, **117**, Springer, (2011)